\documentclass[11pt]{amsart}
\usepackage[margin= 1in]{geometry}

\usepackage{amssymb,amsmath,amsthm,amsfonts,amscd,latexsym, dsfont, color}
\usepackage[dvipsnames]{xcolor}

\usepackage{multirow}

\newcommand{\U}{\mathrm{U}}
\newcommand{\T}{\mathrm{T}}

\newcommand{\Thick}{\mathrm{Thick}}

\newcommand{\titsangle}{\angle_{\mathrm{Tits}}}
\newcommand{\image}{\mathrm{image}} 
 
\newcommand{\graph}{\mathrm{graph}}
\newcommand{\p}{\mathfrak{p}}
\newcommand{\solie}{\mathfrak{so}}
\usepackage{dynkin-diagrams}
\usepackage[title]{appendix}
\usepackage[american]{babel} 

\usepackage{array}
\usepackage[citecolor=ForestGreen,linkcolor=blue,urlcolor=black]{hyperref} 
 \usepackage{graphics}
\usepackage{wrapfig}
\usepackage{pst-node}

\usepackage{tikz-cd} 

\hypersetup{colorlinks}
\usepackage[normalem]{ulem}

\usepackage{setspace} 

\usepackage{tcolorbox}

\usepackage{footnote}

\usepackage{cancel}

\newcommand{\pr}{\mathrm{pr}}

\usepackage{MnSymbol}
\usepackage{enumitem}

\newcommand{\spann}{\mathrm{span}}

\newcommand{\V}{\mathcal{E}} 
\definecolor{bluegreen}{RGB}{5, 170, 204}

\counterwithin{equation}{section}
\definecolor{defnyellow}{RGB}{255,224,102} 
\definecolor{lightblue}{RGB}{170, 220, 255} 
\definecolor{darkblue}{RGB}{0, 125, 230}

\newcommand{\Gtwo}{\mathsf{G}_2}

\newcommand{\Gtwosplit}{\mathsf{G}_2'}

\newtcolorbox{bluebox}[1]{colback=lightblue,colframe=darkblue,fonttitle=\bfseries,title=#1}
\newtcolorbox{redbox}[1]{colback=red!5!white,colframe=red!75!black,fonttitle=\bfseries,title=#1}

\usepackage{pgfplots}

\usepgfplotslibrary{colormaps}
\pgfplotsset{
    compat=newest,
    colormap={mycolormap}{color=(lightgray) color=(white) color=(lightgray) } }

\definecolor{mydarkblue}{RGB}{37, 42, 200}

\definecolor{mygreen}{RGB}{0, 150, 50}

\newcommand{\vis}{\partial_{\mathrm{vis}}}

\newcommand{\Ein}{\mathrm{Ein}} 
 
\newcommand{\fraka}{\mathfrak{a}}

\newcommand{\SO}{\mathrm{SO}}

\newcommand{\PSL}{\mathrm{PSL}}
\newcommand{\Gr}{\mathrm{Gr}} 
\newcommand{\Stab}{\mathrm{Stab}} 
\usepackage{mathrsfs}
\usepackage[utf8]{inputenc}
\usepackage[T1]{fontenc}
\newcommand{\R}{\mathbb R}
\newcommand{\RP}{\mathbb R \mathbb P}
\newcommand{\CP}{\mathbb C \mathbb P}
\newcommand{\C}{\mathbb C}
\newcommand{\Z}{\mathbb Z}

\newcommand{\Ha}{\mathbb{H}} 

\newcommand{\g}{\mathfrak{g}}
\newcommand{\frakk}{\mathfrak{k}}

\newcommand{\X}{\mathbb{X}}
\newcommand{\K}{\mathcal{K}}
\newcommand{\End}{\mathsf{End}} 

\newcommand{\sphere}{\mathbb{S}}

\newcommand{\sllie}{\mathfrak{sl}}
\newcommand{\SL}{\mathrm{SL}}
\newcommand{\Sp}{\mathrm{Sp}}
\newcommand{\GL}{\mathrm{GL}}
\newcommand{\gl}{\mathfrak{gl}}

\newcommand{\tr}{\text{tr}}
\newcommand{\id}{\text{id}}

 \newcommand{\Hom}{\mathrm{Hom}}

\newcommand{\delbar}{\overline{\partial}}

\newcommand{\Iso}{\mathrm{Iso}} 
\usepackage{mathtools}

\newcommand{\Pho}{\mathrm{Pho}}

\newcommand{\rank}{\text{rank}}

\newcommand{\Flag}{\mathrm{Flag}}

\newcommand{\diag}{\mathsf{diag}}

\newcommand{\Hit}{\mathrm{Hit}}

\usepackage{thmtools, thm-restate}
\newcommand{\Diff}{\mathrm{Diff}}

\theoremstyle{plain}

\newtheorem{theorem}{Theorem}[section]
\newtheorem{proposition}[theorem]{Proposition}
\newtheorem{lemma}[theorem]{Lemma}
\newtheorem{example}[theorem]{Example}
\newtheorem{corollary}[theorem]{Corollary}

\newtheorem{definition}[theorem]{Definition}

\newtheorem{remark}[theorem]{Remark}
\newcommand{\alignL}{\begin{flushleft}}
\newcommand{\alignLend}{\end{flushleft}}

\usepackage{soul}

\title{On Einstein Structures for $\SO_0(p,p+1)$-Surface Group Representations}
\author{Colin Davalo and Parker Evans}

\begin{document}

\begin{abstract}
Let $S$ be a closed surface of genus $g \geq 2$. 
We study the cocompact domain of discontinuity $\Omega_{\rho}$ in the Einstein universe $\Ein^{p-1,p}$ defined by Guichard-Wienhard and Kapovich-Leeb-Porti for a class of $p$-Anosov representations $\rho:\pi_1S \rightarrow \SO_0(p,p+1)$ including Hitchin representations, for $p \geq 3$. 
The quotient $M_{\rho} = \rho(\pi_1S)\backslash \Omega_{\rho}$ is abstractly known to be realizable as a fiber bundle over $S$, with unknown fiber of unique homotopy type $F_{\rho}$. We explicitly exhibit $M_{\rho}$ as a smooth $\mathfrak{F}_{\rho}$-fiber bundle over $S$, determining the diffeomorphism type of $M_{\rho}$ and the unique homotopy type  $F_{\rho}$. Surprisingly, in many situations the fiber bundle $\mathfrak{F}_{\rho} \rightarrow M_{\rho}\rightarrow S$ is trivial. 
\end{abstract}

\maketitle 

\tableofcontents

\section{Introduction}
 
In the seminal work \cite{Hit92}, Hitchin discovered surprising connected components of representations of surface groups $\pi_1S$ into higher rank Lie groups that now bear his name. Let $G$ be a split, real, simple Lie group. The \emph{$G$-Hitchin component} $\Hit(S,G)$ in the character variety $\chi(S,G)$ consists of all deformations of representations $\pi_1S \rightarrow G$ that factor through the \emph{principal embedding} $\iota_{pr}:\PSL(2,\R)\to G$ associated with the \emph{principal} $\mathfrak{sl}_2$-subalgebra $\mathfrak{s} < \g$ popularized by Kostant \cite{Kos59}. Hitchin showed the component $\Hit(S,G)$ is a ball in analogy with Teichm\"uller space when $G= \PSL(2,\R)$. Remarkably, all the representations in this component, called \emph{Hitchin representations}, are discrete and faithful. Indeed, they were shown to be \emph{Anosov}, a condition introduced by Labourie \cite{Lab06} specifically for Hitchin representations, but which has since become a topic of independent interest \cite{GW12,GGKW17, KLP17, BPS, DMS25}. Guichard and Wienhard introduced the notion of a $P$-Anosov representation $\Gamma \rightarrow G$ for a hyperbolic group $\Gamma$ and a parabolic subgroup $P<G$, which are representations $\rho:\Gamma \rightarrow G$ that are strongly discrete and faithful, and, in particular, admit a continuous, injective, and $\rho$-equivariant boundary map $\xi_\rho:\partial \Gamma\to G/P$ from the Gromov boundary $\partial\Gamma$ of $\Gamma$ to the flag manifold $G/P$ \cite{GW12}. In this paper, we will only consider  the case $\Gamma = \pi_1(S_{g})$, where $S_g$ is a closed surface of genus $g \geq 2$.

Anosov representations $\rho: \pi_1S \rightarrow G$ are holonomies of locally homogeneous geometric structures. Let us unravel this statement. 
 
Guichard-Wienhard \cite{GW12} and Kapovich-Leeb-Porti \cite{KLP18} have constructed cocompact domains of discontinuity $\Omega$ in flag manifolds $\mathcal{F}=G/Q$ for $(P-)$Anosov representations. Usually, the parabolic subgroup $Q$ is not the same as $P$. 
These domains $\Omega \subset \mathcal{F}$ are associated with the choice of a \emph{balanced Tits-Bruhat ideal}, which contains, in particular, the data of the flag manifold $\mathcal{F}$ in which the domain of discontinuity lies and the parabolic subgroup $P$ for which the Anosov condition holds. The quotient $M_{\rho} \coloneqq \rho(\pi_1S)\backslash \Omega_{\rho}$ of the domain defined via a balanced Tits-Bruhat ideal provides a compact manifold $M_{\rho}$ equipped with a locally homogeneous $(G,\mathcal{F}$)-geometric structure. Since the Anosov condition is  \emph{open} in $\chi(S,G)$, this means the mysterious manifold $M_{\rho}$ carries a large moduli space of $(G,\mathcal{F})$-structures that remain to be understood. 

While the constructions of Guichard-Wienhard and Kapovich-Leeb-Porti apply generally to Anosov representations of any hyperbolic group $\Gamma$, there is a topological structure theorem for the quotients $M_{\rho}$ only in a special case. Let us state this result for surface groups. Alessandrini-Maloni-Tholozan-Wienhard recently proved the following: for Hitchin representations, and more generally representations $\rho: \pi_1S\rightarrow G$ that, up to deformation in the space of $(P-)$Anosov representations, factor through $(\mathrm{P})\SL(2,\R)$, the quotient $M_{\rho}$ is diffeomorphic to a fiber bundle over $S$ \cite{AMTW25}. 
Working backwards, for $\iota:\SL(2,\R)\rightarrow G$ fixed, we shall call a representation $\rho=\iota \circ \rho_0$, for $\rho_0:\pi_1S\rightarrow \SL(2,\R)$ Fuchsian, an $\iota$-\emph{Fuchsian representation}. For $\iota = \iota_{\pr}$, we call an $\iota$-Fuchsian to be \emph{Fuchsian-Hitchin}. 

The topological structure theorem of \cite{AMTW25} is not true in general, even in rank one groups $G$. Indeed, in \cite{GLT88}, Gromov-Lawson-Thurston provide examples (in our modern language) of Anosov representations $\rho:\pi_1S\rightarrow \SO_0(4,1)$ whose quotients $M_{\rho}$, for $X = \mathbb{S}^3 = \partial_{\infty}\Ha^4$ are not realizable as fiber bundles over the original surface $S$. In this case, the cocompact domain of discontinuity $\Omega_{\rho}$ obtains the form $\Omega_{\rho} = \sphere^3 \backslash \Lambda$, where $\Lambda \cong \sphere^1$ is the limit set of the surface group $\rho(\pi_1S)$; this domain falls under the general constructions of \cite{GW12, KLP18}. 

Presently, we consider surface group representations that are Anosov deformations of $\iota$-Fuchsians, so that $M_{\rho}$ is known a priori to diffeomorphic to the total space $\mathscr{F}_{\rho}$ of a fiber bundle $\mathscr{F}_{\rho}\rightarrow S$ by \cite{AMTW25}. In this case, the domain $\Omega_{\rho}$ deformation retracts onto the fiber $\mathfrak{F}_{\rho}$ of $\mathscr{F}_{\rho} \rightarrow S$, hence we may speak of \emph{the homotopy type} of $\mathfrak{F}_{\rho}$ (well-defined among any such fibration). 
Here is the crux: the construction of the domain of discontinuity $\Omega_{\rho}$ does not by itself elucidate the structure of $M_{\rho}$. In particular, two essential questions remain open in general: 
\begin{enumerate}[noitemsep, label=(\roman*)]
    \item What is the homotopy type of $\mathfrak{F}_{\rho}$? 
    \item What is the diffeomorphism type of the quotient $M_{\rho}$? 
\end{enumerate}
We briefly remark that in the case of $\rho$ an $\iota$-Fuchsian, one can often produce an explicit fibration $\pi: M_{\rho} \rightarrow S$, allowing a refinement of the second question above: namely, \medskip 

(iii) What is the smooth isomorphism type of the fiber bundle $\pi: M_{\rho} \rightarrow S$? \medskip

To answer questions (i) and (ii) in the general case, it is sufficient to consider $\iota$-Fuchsian representations, as the homotopy type of $\mathfrak{F}_{\rho}$ and the diffeomorphism type of $M_{\rho}$ are invariant under Anosov deformations by a variant of the Ehresmann-Thurston principle. In particular, for Hitchin representations it is sufficient to consider a single Fuchsian-Hitchin representation. There has been a lot of work recently to investigate the $(G,X)$-manifolds associated to Anosov surface group representations via balanced thickenings, with a focus on these two questions \cite{Har26, AMTW25, Dav25, DE25a, RT25, ADL24, DS20, CTT19}. Other problems regarding these $(G,X)$-manifolds, such as providing a synthetic characterization of them, have also been studied in \cite{Dav26,NR25, GW08}.

A particular but already interesting instance of the \cite{GW12,KLP18}-domains are those associated with \emph{Tits metric thickenings}. Namely, in this case, $\Omega = \mathcal{F}\backslash K_{\Lambda}$, where the \emph{thickening} $K_{\Lambda}$ of the limit set $\Lambda=\image(\xi_{\rho})$ is given by a $(\frac{\pi}{2}$-)neighborhood of $\Lambda \subset G/P$ in the flag manifold $\mathcal{F}= G/Q$ with respect to the Tits angle metric. The domains considered in this paper will be of this form. Moreover, to the best of our knowledge, all concrete examples for which the topology of $\mathfrak{F}_{\rho}$ is known fall under this umbrella, as we recall in Section \ref{Sec:RelatedWork}.

\subsection{Main Results}

In this paper, we study the topology of the $(G,X)$-quotients of Anosov deformations of $\iota$-Fuchsians for a new infinite family, namely $(G,X)=(\SO_0(p,p+1), \Ein^{p-1,p})$, for $p \geq 3$, where $\mathcal{F}=\Ein^{p-1,p}$ is the \emph{Einstein universe} of isotropic lines in pseudo-Euclidean space $\R^{p,p+1}$. We find both (i) the homotopy type of $\mathfrak{F}_{\rho}$ and (ii) the diffeomorphism type of the manifold $M_{\rho}$ as a byproduct of our solution to  problem (iii). Namely, for $\iota$-Fuchsians, we determine the smooth isomorphism type of a specific fiber bundle $\pi:M_{\rho} \rightarrow S$ defined via the \emph{nearest point projection} as in \cite{Dav25}. Our results handle all $P_p$-Anosov (or \emph{$p$-Anosov}, for short) deformations of $\iota$-Fuchsian representations, where $\iota: \SL(2,\R)\rightarrow \SO_0(p,p+1)$  and $P_p<G$ is the stabilizer of an isotropic $p$-plane. 

In this case, the domain parabolic $Q<G$ is the stabilizer of an isotropic line, so $\mathcal{F} =G/Q = \Ein^{p-1,p}$. The relevant domain of discontinuity $\Omega_{\rho}  \subset \mathcal{F}$ is given by 
\begin{align}\label{Omega_Thickening_Intro}
     \Omega_{\rho} = \mathcal{F} \backslash \bigcup_{x \in \partial\pi_1S} K_{\xi(x)}. 
\end{align}
In this construction, $\Iso_p(\R^{p,p+1})=G/P_p$ is the space of isotropic $p$-planes in $\R^{p,p+1}$, the map $\xi: \partial\pi_1S \rightarrow \Iso_p(\R^{p,p+1})$ is the (unique) $\rho$-equivariant $p$-Anosov boundary map, and the thickening of a point $T \in \Iso_p(\R^{p,p+1})$ is given by 
\[ K_{T}= \left \{\ell \in \Ein^{p-1,p} \mid \ell \subset T \right\} .\]
Since $\titsangle(\ell, T) \leq \frac{\pi}{2} \iff \ell \subset T$ for $\ell \in \Ein^{p-1,p}$ and $T \in \Iso_p(\R^{p,p+1})$, the domain $\Omega_{\rho}$ is indeed constructed by Tits metric thickening. 

\subsubsection{Topology of the fiber}

For Hitchin representations, the case $p=2$ was treated by Collier-Tholozan-Toulisse in \cite{CTT19}, where the fiber is geometrically $\Ein^{1,0}$, a topological circle.  
We find the homotopy type of the fiber $\mathfrak{F}_{\rho}$ of $M_{\rho}$ in the remaining cases $p \geq 3$. 

\begin{theorem}
\label{IntroThm:EinFibers}
Let $p \geq 3$ be an integer and $\rho:\pi_1S \rightarrow \SO_0(p,p+1)$ a Hitchin representation. The manifold $M_{\rho}$ is diffeomorphic to an $\mathfrak{F}_{\rho}$-fiber bundle over $S$, where 
\begin{enumerate}[label=(\roman*), noitemsep]
    \item If $p$ is even, then $\mathfrak{F}_{\rho}$ is the unit tangent bundle $\mathrm{T}^1\mathbb{RP}^{p-1}$ of $\mathbb{RP}^{p-1}$. 
    \item If $p$ is odd, then $\mathfrak{F}_{\rho}$ is the space $\Ein^{p-1,p-2}$ of isotropic lines in $\R^{p,p-1}$.
\end{enumerate}
\end{theorem}

In order to prove these results, we first describe the fiber using the nearest point projection technique inside the symmetric space from \cite{Dav25}. This allows us to describe the fiber as a \emph{base of pencil}. Due to the simple nature of the embedding of $\Ein^{p-1,p}$ in the visual boundary $\vis\X$, where $\X$ denotes the $\SO_0(p,p+1)$-symmetric space, we can show the associated base of pencil is diffeomorphic to the total space of a sphere bundle $\mathbb{S}(E) \rightarrow \mathbb{RP}^{p-1}$, with fiber $\mathbb{S}^{p-2}$, associated to a vector bundle $E \rightarrow \RP^{p-1}$ of rank $p-1$. 
The topology of a base of pencil remains invariant under certain \emph{regular} deformations, and we use this freedom crucially. We carefully deform the original pencil $\mathcal{P}$ to a suitable ``simplification'' $\mathcal{P}_0$ via regular deformation. The associated vector bundle $E =E(\mathcal{P}_0)$ is then explicitly determined for the simplified pencil.
\medskip

We extend Theorem \ref{IntroThm:EinFibers} by considering non-principal representations of $\SL(2,\R)$. To describe our more general results, we need to recall some representation theory of $\sllie_2\R$. 

Recall that for any integer partition of $n$, there is a unique associated representation of $\sllie(2,\R) \rightarrow \sllie(n,\R)$ up to isomorphism.\footnote{There is one exception when all parts of the integer partition are even. See \cite[Theorem 9.3.3]{CM17}.}
Let $I$ denote an integer partition with one odd part and all even parts of even multiplicity.  
In this case, the associated representation $\sllie_2\R \rightarrow \sllie(2p+1)$ factors (essentially uniquely) through a copy of $\solie(p,p+1)$. Denote $\iota_I: \SL(2,\R)\rightarrow \SO_0(p,p+1)$ as the associated representation and $\Ha^2_I$ the sub-symmetric space of the image subgroup. The hyperbolic plane $\Ha^2_I$ is special in the sense that the Cartan projection $\mu: \T\X \rightarrow \overline{\mathfrak{a}}^+$ maps $\T\Ha^2_I$ to a subset of the model Weyl chamber $\overline{\mathfrak{a}}^+$ that avoids the walls (hyperplanes) defined by the short roots in the associated type $B_p$ root system. We call this property \emph{Einstein (Ein)-regularity}. Geometrically, this means $\iota_{I}$ will map Fuchsian representations in $\SL(2,\R)$ to $p$-Anosov representations in $\SO_0(p,p+1)$.

We show the homotopy type of $\mathfrak{F}_{\rho}$ is the same for any deformation $\rho$ of an $\iota_I$-Fuchsian inside the space of $p$-Anosov representations as in the Hitchin case. 
\begin{corollary}[Fibers for Deformations of $\iota$-Fuchsians]\label{Cor:FiberIotaFuchsians}
Let $p \geq 3$ be an integer and $\rho: \pi_1S\rightarrow \SO_0(p,p+1)$ a deformation of an $\iota_I$-Fuchsian in the space of $p$-Anosov representations. Then $M_{\rho}$ is diffeomorphic to an $\mathfrak{F}_{\rho}$-fiber bundle over $S$, where $\mathfrak{F}_{\rho}$ is as in Theorem \ref{IntroThm:EinFibers}. 
\end{corollary}

There is a corollary to the case $p=3$ in Theorem \ref{IntroThm:EinFibers}. Here, let $\Gtwosplit<\SO_0(3,4)$ denote the split real form of the exceptional Lie group $\Gtwo^\C$. Since $\Hit(S,\Gtwosplit) \hookrightarrow \Hit(S,\SO_0(3,4))$, we find the following corollary by topological invariance. 
\begin{corollary}\label{Cor:G2EinsteinFibers}
For $\rho: \pi_1S \rightarrow \Gtwosplit$ Hitchin, $M_{\rho}$ is diffeomorphic to an $\Ein^{2,1}$-fiber bundle over $S$. 
\end{corollary}

We use Corollary \ref{Cor:G2EinsteinFibers} as inspiration in \cite{DE25a}, where we solve the converse problem: given certain surface group representations $\rho: \pi_1S \rightarrow \Gtwosplit$, including, but not limited to $\Gtwosplit$-Hitchin representations, we construct an associated 5-manifold $N_{\rho}$, realized as an $\Ein^{2,1}$-fiber bundle over $S$, which carries a $(\Gtwosplit,\Ein^{2,3})$-structure whose holonomy $\pi_1N \rightarrow \Gtwosplit$ descends to $\pi_1S$ as the original representation $\rho$. The construction of $N_{\rho}$ is with differential geometry, rather than domains of discontinuity, but is shown to relate to the \cite{GW12, KLP18}-manifolds nevertheless. \medskip 

\subsubsection{Global Topology}

We now describe our results and technique to determine the global topology of these  quotients $M_{\rho}$. In general, the results of \cite{AMTW25} show that $M_{\rho}$ can be realized as the associated fiber bundle $\mathcal{P}\times_{\phi} \mathfrak{F}_{\rho}$ to a principal $\U(1)$-bundle $\mathcal{P}\rightarrow S$ via a circle action $\phi$ on some unknown manifold $\mathfrak{F}_{\rho}$. In the present circumstance, excluding $p=4k+2$, we avoid the associated bundle approach, as we can instead provide a more concrete description that clarifies which topological invariants of $\rho$ influence the global topology of $M_{\rho}$ and reveals global triviality of the fiber bundle $\mathfrak{F}_{\rho}\rightarrow M_{\rho}\rightarrow S$ in many situations.

The main idea is to use the base of pencil construction described earlier now inside the $\Ein^{p-1,p}$-fiber bundle 
$E=\tilde{S}\times_{\rho} \Ein^{p-1,p}$ over the surface $S$ associated to $\rho$. 
This approach allows us to assemble the pointwise constructions of bases of pencils into a codimension two fiber subbundle $\mathfrak{E}$ of $E$ that is diffeomorphic to $M_{\rho}$. By careful study of $\mathfrak{E}$, we obtain an explicit description of the quotient $M_{\rho}$, whose topology we now describe. \medskip

Let $\mathcal{K}$ be the total space of a degree $2g-2$ real vector bundle over $S$, which is diffeomorphic to the holomorphic cotangent line bundle $\K_{\Sigma} = (\T^{1,0}\Sigma)^*$ for any choice of Riemann surface $\Sigma=(S,J)$ on $S$. We remark that the description of the global topology inevitably varies in nature when $p$ is even and odd due to the discrepancy in the fibers from Corollary \ref{Cor:FiberIotaFuchsians}.

Now, let $p$ be odd. In this case, using that $\Ein^{k,\ell} \cong (\sphere^k\times \sphere^{\ell})/(-\id, -\id)$ is nearly a product, the description of $M_{\rho}$ is simpler. Given a $\Z_2$-graded vector bundle $E_1 \oplus E_2$, we shall similarly write $\Ein(E_1 \oplus E_2)$ for the associated fiber bundle $\sphere(E_1) \oplus \sphere(E_2)/\sim$, under the quotient $(p,u,v) \sim (p,-u,-v)$. 
For the following theorem statement, recall that smooth vector bundles $E \rightarrow S$ on a closed oriented surface $S$ are determined up to isomorphism by their first and second Stiefel-Whitney classes 
$w_1(E) \in H^1(S,\Z_2)$ and $w_2(E) \in H^2(S,\Z_2) \cong \Z_2$. It turns out the global topology of $M_{\rho}$ depends only on the parity of a single integer.

\begin{theorem}[Global Topology for $\iota$-Fuchsians, $p$ odd]\label{Thm:GlobalTopologyGeneral_OddIntro}
Let $p \geq 5$ be an odd integer and $\rho:\pi_1S\rightarrow \SO_0(p,p+1)$ a deformation of an $\iota_I$-Fuchsian representation in the space of $p$-Anosov representations. Suppose the unique odd part of the partition $I$ is $2q+1$. 

Then $M_\rho=\rho(\pi_1S)\backslash\Omega_\rho$ is diffeomorphic to the fiber bundle $\Ein(E_+\oplus E_-)$, where $E_+, E_-$ are orientable real vector bundles over $S$ of ranks $p$ and $p-1$ respectively, with second Stiefel-Whitney classes $w_2(\rho)\coloneqq (g-1)\frac{p-q}{2} \bmod 2 \in \Z_2$. 
\end{theorem}

\begin{corollary}\label{Cor:IntroOddCaseGeneral}
When $p\geq 5$ is odd, the quotient $M_{\rho}$  satisfies the following: 
\begin{enumerate}[label=(\roman*)]
    \item The global topology of $M_{\rho}$ depends only on $w_2(\rho)$. 
    \item When $\rho$ is an $\iota_I$-Fuchsian representation, the fiber bundle $\pi: M_{\rho}\rightarrow S$ defined via nearest point projection is non-trivial if and only if $w_2(\rho)\neq 0$.
\end{enumerate}
\end{corollary}

For clarity, we emphasize the topology in the Hitchin case, where $p=3$ is anomalous due to the presence of the circle bundle $\sphere(E_-)$. Here, denote $\varepsilon^k_{\R}\rightarrow S$ a trivial vector bundle of rank $k$.
\begin{corollary}\label{Cor:HitchinOdd}
Let $p \geq 3$ be an odd integer and $\rho:\pi_1S \rightarrow \SO_0(p,p+1)$ a Hitchin representation. The manifold $M_{\rho} =\rho(\pi_1S)\backslash \Omega_{\rho}$ is diffeomorphic to the following fiber bundle: 
\begin{enumerate}[noitemsep, label=(\roman*)]
    \item $\Ein(\varepsilon^3_{\R} \oplus \K^3)$ when $p = 3$.
    \item $S \times \Ein^{p-1,p}$ when $p \geq 5$. 
\end{enumerate}
\end{corollary}
The result of Corollary \ref{Cor:HitchinOdd} says that, excluding the exceptional case $p=3$, the fiber bundles for Hitchin representations are globally trivial in the odd case. Indeed, in the case of Fuchsian-Hitchin $\rho$, Corollary \ref{Cor:IntroOddCaseGeneral}(ii) applies and $\pi: M_{\rho} \rightarrow S$ is a trivial fiber bundle.  
\medskip 

When $p$ is even, we first give a general description of the fiber bundle obtained with fiber $\T^1\mathbb{RP}^{p-1}$. To this end, we note that there is a natural $\U(1)$-action on $\T^1 \mathbb{P}(V)$ for any real vector space $V$ of dimension $p$ equipped with an orthogonal complex structure $J$, given as follows. The tangent space at $\ell\in \mathbb{P}(V)$ decomposes orthogonally as $\T_{\ell}\mathbb{P}(V)=J\ell\oplus W$, where $W=(\C\ell)^\bot$ is preserved by $J$. Hence, for $\lambda\in \U(1)$ one can define ${\eta_0}_{|\ell}(\lambda)\in \GL( \T_\ell\mathbb{P}(V))$ so that $\eta_0(\lambda)$ acts trivially on $J\ell$ and acts by complex scalar multiplication on $W$. Since $J$ is orthogonal, $\eta_0$ preserves $\T^1\mathbb{P}(V)$ and the fiberwise actions produce a map $\eta_0: \U(1) \rightarrow \Diff(\T^1\RP^{p-1})$. 

Finally, for the following theorem statement, let $\mathcal{T}$ be the unit tangent bundle of $S$, seen as a principal $\U(1)$-bundle. We first give a description of $M_{\rho}$ as a certain kind of associated fiber bundle to a triple $(\mathcal{T},\mathcal{F}, \eta_0)$, where $\mathcal{F}$ is a $\T^1\RP^{p-1}$-fiber bundle over $S$. 
\begin{theorem}[Global Topology for $\iota$-Fuchsians, $p$ even]\label{Thm:GlobalTopologyGeneral_EvenIntro}
Let $p \geq 4$ be an even integer, $\rho: \pi_1S\rightarrow \SO_0(p,p+1)$ a deformation of an $\iota_I$-Fuchsian through $p$-Anosov representations, and $E_+ \rightarrow S$ an orientable rank $p$ vector bundle with $ w_2(E_+) = w_2(\rho)$, where $w_2(\rho)$ is as in Theorem \ref{Thm:GlobalTopologyGeneral_OddIntro}.
 
Then the manifold $M_{\rho} = \rho(\pi_1S)\backslash \Omega_{\rho}$ is diffeomorphic to the fiber bundle $\mathcal{T}\oplus\T^1\mathbb{P}(E_+)/\sim$, where $(x\cdot \lambda,v)\sim (x, \eta_0(\lambda)\cdot v)$ for all $\lambda \in \U(1)$ and $(x,v)\in \mathcal{T}\oplus\T^1\mathbb{P}(E_+)$. 
\end{theorem}

Recall the fiber bundle direct sum $\mathcal{T}\oplus\T^1\mathbb{P}(E_+)$ is defined so that for all $x\in S$, the fiber is given by $\big(\mathcal{T}\oplus\T^1\mathbb{P}(E_+) \big)|_x=\mathcal{T}|_x\times\T^1\mathbb{P}(E_+)|_x$.

\medskip

When $p$ is a multiple of four, we can simplify the description and obtain the following: 

\begin{theorem}
Let $p \geq 4$ be a multiple of four,
$S=S_g$ a closed surface of genus $g$, and $\rho: \pi_1S\rightarrow \SO_0(p,p+1)$ a deformation of an $\iota_I$-Fuchsian representation inside the space of $p$-Anosov representations. Then $M_{\rho} = \rho(\pi_1S)\backslash \Omega_{\rho}$ is diffeomorphic to $\T^1\mathbb{P}(E_+)$, where $E_+$ is as in Theorem \ref{Thm:GlobalTopologyGeneral_EvenIntro}.
\end{theorem}

\begin{corollary}\label{Cor:IntroEvenTopology}
When $p$ is a multiple of four, the smooth manifold $M_{\rho}$ satisfies the following:
\begin{enumerate}[noitemsep,label=(\roman*)]
    \item The global topology depends only on $w_2(\rho)$.
    \item If $\rho$ is an $\iota_I$-Fuchsian representation, the nearest point projection $\pi: M_{\rho} \rightarrow S$ defines a nontrivial fiber bundle if and only if $w_2(\rho) \neq 0$. 
    \item If $w_2(\rho)=\frac{p-q}{2}(g-1) \bmod 2 = 0$, then $M_{\rho}$ is diffeomorphic to $S \times \T^1\RP^{p-1}$.
\end{enumerate}
\end{corollary}

Thus, in the Hitchin case for $p=4k$, Corollary \ref{Cor:IntroEvenTopology}(iii) asserts that $M_{\rho}$ is a trivial product. 

\subsection{Related Work and Further Remarks}\label{Sec:RelatedWork}

The topology of $\mathfrak{F}$ and $M$ obtained for $G$-Hitchin representations in flag manifolds $G/Q$ is known explicitly only in some cases. Here, it is natural to also consider $G^\C$-\emph{quasi-Hitchin representations}, namely all $(P-)$Anosov deformations of a Fuchsian-Hitchin representation $\rho_0: \pi_1S\rightarrow \PSL(2,\R) \hookrightarrow G$ under the inclusion $G\hookrightarrow G^\C$ of $G$ into its complexification. 
This notion generalizes the classical case of \emph{quasi-Fuchsian} representations $\pi_1S \rightarrow \PSL(2,\C)$. 

Table \ref{Table:TopologyFibers} below summarizes all examples where the topology of the fiber is known. Here, we omit redundancies: for example, $X =\RP^{2n-1}$ is a flag manifold for both $G =\PSL(2n,\R)$ and $G= \mathrm{PSp}(2n,\R)$ and also $\Ein^{2,3}$ is a flag manifold for $G =\SO_0(3,4)$ and $G=\Gtwosplit$. \medskip 
\begin{table}[ht]
\renewcommand{\arraystretch}{1.3}
    \centering
    \begin{tabular}{|c|c|c|c|c|c|}\hline 
    $G$ & $G/P$ & $G/Q$ & $\mathfrak{F}$ & Reference(s) \\ \hline  
      $\SL(3,\C)$ & $\Flag(\C^3)$ & $\Flag(\C^3)$ & $(\sphere^2\times\sphere^2)\# (\sphere^2 \times \sphere^2)$ & \cite{Har26} \\ \hline  
    $\SO_0(p,p+1)$ ($p\geq 3$) & $\Iso_p(\R^{p,p+1})$ & $\Ein^{p-1,p}$ & $\begin{cases} \Ein^{p-1,p-2} & p \; \text{odd}\\ 
    \T^1\RP^{p-1} & p \; \text{even}\end{cases}$ & * \\ \hline 
    $ \Gtwosplit$ & $\Pho^\times$ & $\Pho^\times$ & $\RP^2\times \mathbb{S}^1$ & \cite{DE25a} \\ \hline
    $\SL(3,\R)$ & $\text{Flag}(\R^3)$ & $\text{Flag}(\R^3)$ & $\mathbb{S}^1\sqcup\mathbb{S}^1\sqcup\mathbb{S}^1$ & \cite{NR25, RT25}\\ \hline
    $\Sp(4,\C)$ & $\CP^{3}$ & $\mathrm{Lag}(\C^4)$ & $\CP^2\# \overline{\CP^2}$ & \cite{AMTW25, Har26} \\ \hline 
    $\PSL(2n,\R)$ & $\Gr_{n}(\R^{2n})$ & $\RP^{2n-1}$ & $\T^1\RP^{n-1}$ & \cite{ADL24} \\ \hline 
    $\PSL(2n,\C)$ & $\Gr_{n}(\C^{2n})$ & $\CP^{2n-1}$ & $\T^1\sphere^{2n-1}/\mathrm{U}(1)$ & \cite{ADL24} \\ \hline 
    $\SO_0(2,n+1)$ & $\Ein^{1,n}$ & $\Pho(\R^{2,n+1})$ & $\Pho(\R^{2,n})$ & \cite{CTT19} \\ \hline 
    $\SO_0(2,3)$ & $\Pho(\R^{2,3})$ & $\Ein^{1,2}$ & $\mathbb{S}^1$ & \cite{CTT19} \\ \hline 
    $\PSL(4,\R)$ & $\Gr_2(\R^4)$ & $\RP^3$ & $\mathbb{S}^1\sqcup\mathbb{S}^1$ & \cite{GW08}\\ \hline
    \end{tabular}
    \caption{The topological type of the fiber $\mathfrak{F}$ for ($\frac{\pi}{2}$-Tits metric thickening) domains of discontinuity $\Omega \subset G/Q$ when $\rho$ is a Hitchin or quasi-Hitchin representations in $G$. Here, $P$ is the Anosov parabolic. }
    \label{Table:TopologyFibers}
\end{table}

We now make some further remarks. 

\begin{remark}
\begin{itemize}
    \item Let us fix $\Gamma =\pi_1S, G,P,Q$. The topological type of the fiber $\mathfrak{F}$ of such domains of discontinuity is not independent of $\iota$ in general. For example, as explained in \cite[Example 6.9]{Dav25}, for $\iota_{\mathrm{red}}:\SL(2,\R)\rightarrow \SL(3,\R)$ the reducible representation and $\iota_{\mathrm{pr}}:\SL(2,\R)\rightarrow \SL(3,\R)$ the irreducible representation, with $G/P=\Flag(\R^3)=G/Q$ and $\Omega_{\rho}$ the $\frac{\pi}{2}$-Tits metric thickening domain, the fiber $\mathfrak{F}$ for $\iota_{\mathrm{red}}$-Fuchsians is a circle, while for $\iota_{\mathrm{pr}}$, the fiber $\mathfrak{F}$ is a disjoint union of three circles. 
    \item In the case of $G/Q=\mathrm{Lag}(\C^4)$, we note that \cite{AMTW25} classified the fiber $\mathfrak{F}$ topologically for $\rho$ quasi-Hitchin, and \cite{Har26} extended their work, exhibiting $M_{\rho} \rightarrow S$ as a smooth $\mathfrak{F}$-fiber bundle, where $\mathfrak{F}= \CP^2 \# \overline{\CP^2}$ with its standard smooth structure. 
    \item \cite{CTT19} actually proves more: for $\rho:\pi_1S\rightarrow \SO_0(2,n+1)$ any maximal representation, including those which are not $\iota$-Fuchsian deformations, and $G/Q = \Pho(\R^{2,n+1})$, they exhibit a smooth fiber bundle  $M_{\rho}\rightarrow S$ with fiber $\Pho(\R^{2,n})$. 
    \item In \cite{ADL24}, alongside $\iota_{\mathrm{pr}}$, the authors also consider the `most reducible' representation $\iota_0$, with irreducible subspaces of dimension $2+2+\cdots +2$. The fiber for $\iota_0$-Fuchsians is shown to be the same as for Fuchsian-Hitchin representations.
    \item In \cite{ADL24}, to be more precise, the fiber in $\CP^{2n-1}$ is the total space of the sphere bundle $\sphere(\T\CP^{n-1} \oplus \varepsilon^1_\R)$ over $\CP^{n-1}$, where $\varepsilon^1_\R \rightarrow \CP^{n-1}$ is a trivial real rank one line bundle. This sphere bundle is equivalently a $\mathrm{U}(1)$-quotient of $\T^1\sphere^{2n-1}$.
\end{itemize}
\end{remark}

We note two further related works. In \cite{DS20}, Dumas-Sanders extensively studied the complex-analytic properties of the \cite{GW12, KLP18}-manifolds when $G$ is complex, and proved that for $G=\SL(3,\C)$ and $G/P=G/Q = \Flag(\C^3)$, the quotient $M_{\rho}$ is a fiber bundle over $S$ when $\rho$ is quasi-Hitchin. In \cite{Har26}, Hart studies the topology of domains of discontinuity $\Omega_{\rho}$ for three-dimensional complex flag manifolds $X \in \{\Flag(\C^3), \mathrm{Lag}(\C^4), \CP^3\}$, with respective group $G \in \{\SL(3,\C), \Sp(4,\C), \SL(4,\C)\}$, and $\rho: \pi_1S\rightarrow G$ is any $\iota$-Fuchsian. In this case, $\mathfrak{F}_{\rho}$ is homotopy equivalent to a simply-connected 4-manifold, and Hart's work hinges around the classification of circle actions on such manifolds to determine the topological type of $\mathfrak{F}_{\rho}$. Using the structure theorem of \cite{AMTW25}, a global description of $M_{\rho}$ is achieved as an associated fiber bundle. \medskip 

As already mentioned, \cite{ADL24} obtained results on the  $\RP^{2n-1}$ and $\CP^{2n-1}$-manifolds associated to Hitchin and quasi-Hitchin representations, respectively. This work predated the ideas of \cite{Dav25} used presently. Consequently, the arguments in \cite{ADL24} are rather different from those in the present work. Most notably, in \cite{ADL24} an explicit developing map is verified in the Fuchsian cases, and shown to diffeomorphically map onto the domain $\Omega_{\rho}$. From this argument, the topology of the fiber $\mathfrak{F}$, and the global topology of $M$ are read off from a certain model space $N \rightarrow S$ that is a fiber bundle over the surface. 

Presently, we circumvent a tedious developing map calculation by using \cite{Dav25} and bases of pencils. This leads to shorter proofs. 
Moreover, we find more general results, namely, the homotopy type of $\mathfrak{F}_{\rho}$ and diffeomorphism type of $M_{\rho}$ for all $p$-Anosov deformations of $\iota$-Fuchsians.

\subsection*{Acknowledgments}
The authors thank Daniele Alessandrini for discussions regarding this work and the results of \cite{AMTW25}. 
C. Davalo was funded by the European Union via the ERC 101124349 ``GENERATE''. Views and opinions expressed are however those of the authors only and
do not necessarily reflect those of the European Union or the European Research
Council Executive Agency. Neither the European Union nor the granting authority
can be held responsible for them.

\section{Preliminaries} 
In this section, we describe various background
information needed for the paper. 

\subsection{The Model Space \texorpdfstring{$\Ein^{p-1,p}$}{Ein(p-1,p)} in the Visual Boundary}\label{Subsec:EinPrelims}  

First, we recall details on the group $G =\SO_0(p,p+1)$ and its maximal parabolic subgroup $P_1 <G $ such that $G/P_1 \cong \Ein^{p-1,p}$. 
We describe the model space $\Ein^{p-1,p}$ geometrically and then realize it inside the visual boundary of the $G$-symmetric space $\X$. The latter embedding is essential for everything that follows. 

Recall that $\R^{p,p+1}$ denotes pseudo-Euclidean space $\R^{p,p+1} =(\R^{2p+1},q_{p,p+1})$, namely the vector space $\R^{2p+1}$ equipped with a non-degenerate quadratic form $q=q_{p,p+1}$ of signature $(p,p+1)$. The group $\SO(p,p+1)$ is the stabilizer in $\SL(2p+1,\R)$ of $q_{p,p+1}$ and has two connected components. We denote by $\SO_0(p,p+1)$ the identity component of $\SO(p,p+1)$.

\begin{remark}
Since the vector space $\R^{p,p+1}$ contains both spacelike vectors $(q(x) > 0)$, and timelike vectors $(q(x) <0)$, we frequently denote $Q_+(U)$ and $Q_-(U)$ for the subsets of unit spacelike and unit timelike elements, respectively, in a given subspace $U < \R^{p,p+1}$. 
\end{remark}

The \emph{Einstein universe of signature $(p-1,p)$} is is the projective lightcone in $\R^{p,p+1}$: 
\[ \Ein^{p-1,p} := \{ [x] \in \mathbb{P}(\R^{p,p+1} )\; | \; q_{p,p+1}(x) = 0\}. \] 
The space $\Ein^{p-1,p}$ is also commonly referred to as the \emph{projective null quadric} in $\R^{p,p+1}$.

We shall also need one more flag manifold of $G$, namely
\[ \Iso_p(\R^{p,p+1}) = \{T \in \Gr_p(\R^{p,p+1}) \mid q_{p,p+1}|_T\equiv0\},\]
the Grassmannian of isotropic $p$-planes in $\R^{p,p+1}$. 

To describe some Lie-theoretic preliminaries, pick a basis $B=(e_i)_{i=1}^{2p+1}$ of $\R^{p,p+1}$ such that 
such that the quadratic form $q$ obtains the matrix form $[q]$ as follows : 
\begin{align}\label{Q_ModelBasis}
     [q] = \begin{pmatrix} 
                     & & & & & &1 \\ 
                     & & & &&\reflectbox{$\ddots$}&  \\
                     & & & & 1& &  \\
                     & & &-1 & & & \\     
                     & & 1& & & &\\ 
                     &\reflectbox{$\ddots$}& & & & & \\   
                     1 & & & & & &\\
 \end{pmatrix}. 
 \end{align}
In this basis, a Cartan subalgebra $\mathfrak{a} \subset \g$ is given by:
 \begin{align}\label{ModelCSA}
     \mathfrak{a} = \{ \; \diag(\lambda_1,\, \lambda_2,\, \dots,\, \lambda_p, \,0,\, -\lambda_p,\, \dots, \,-\lambda_2, \,-\lambda_1)\in \gl_{2p+1}(\R) \: | \; \lambda_i \in \R\}.
 \end{align}
We will treat a general element $X \in \mathfrak{a}$ as having these coordinates $(\lambda_i)$. 
A choice of simple roots $\Delta\coloneqq(\alpha_i)_{i=1}^p$ for the real root system $\Sigma(\g, \mathfrak{a})$ is as follows: define $\alpha_i \in \mathfrak{a}^*$ by $\alpha_i=\lambda_{i+1}^*-\lambda_i^*$ for $1 \leq i \leq p-1$ and $\alpha_p = \lambda_p^*$. Note that $\alpha_p$ is the unique short simple root. 

 The $G=\SO_0(p,p+1)$-Riemannian symmetric space $\X$ can and will be equivariantly identified with the model space
$\Gr_{(p,0)}(\R^{p,p+1})$, the spacelike $p$-Grassmannian of $\R^{p,p+1}$. Indeed, a maximal compact subgroup $K<G$ is isomorphic to $\SO(p)\times \SO(p+1)$ and realized as the stabilizer of a splitting $\R^{p,p+1}=\R^{p,0} \oplus \R^{0,p+1}$. Equivalently, $K =\Stab_{G}(P)$ for a point $P \in \X$. The tangent space  $\T_P\X=\T_P\Gr_{(p,0)}(\R^{p,p+1})$ is naturally identified with $\Hom(P, P^\perp)$.
The Riemannian metric of $\X$ can be written in this model as:
 \[ g_{P}(\phi, \psi) = -\tr(\phi^{*q} \circ \psi). \]
We will also identify $\T_P\X$ with the subset of $q$-anti-symmetric endomorphisms of $\R^{p,p+1}$, which can be written $A_\phi \coloneqq \begin{pmatrix} 0 & -\phi^{*q} \\ \phi & 0 \end{pmatrix} \in \mathfrak{so}(p,p+1)$, in block form relative to $\R^{p,p+1} = P\oplus P^\bot$. \medskip

The basis $B$ yields a basepoint $P_0 \in \X$. Indeed, we can set 
\begin{align}\label{X_Basepoint}
  P_0:=\spann\langle e_1+e_{2p+1},e_2+e_{2p}, \cdots, e_{p}+e_{p+2}\rangle.
\end{align}
 We now set $K:= \Stab_{G}(P_0)$. Under the corresponding Cartan decomposition $\g = \frakk \oplus \mathfrak{p}$, the model Cartan subalgebra $\fraka$ in \eqref{ModelCSA} satisfies $\fraka \subset \mathfrak{p}$. Viewing $\T_{P_0}\X = \mathfrak{p}$, then we can treat $\fraka \subset \T_{P_0}\X$. We fix the following (open) model Weyl chamber $\mathfrak{a}^+$: 
 \begin{align}\label{ModelWeylChamber}
    \mathfrak{a}^+ = \{ \; \diag(\lambda_1, \lambda_2,\dots, \lambda_p, 0, -\lambda_p, \dots, -\lambda_2, -\lambda_1)\in \mathfrak{a} \: | \; \lambda_1 > \lambda_2 > \cdots  > \lambda_p > 0\},
 \end{align}
the intersection of the half-planes $\{ t \in \mathfrak{a} \mid \alpha_i (t)>0\}$, for $i \in \{1,2,\dots, p\}$. The \emph{Cartan projection} $\mu: \T\X\rightarrow \overline{\mathfrak{a}}^+$ is the map which takes $X \in \T\X$ to the unique element of its $G$-orbit in $\overline{\mathfrak{a}}^+$. 

The space $\X$ is a \emph{Hadamard manifold}, and admits a compactification $\vis\X$ called the \emph{visual boundary}, described in detail in \cite{Ebe96, BH99}. In particular, $\vis\X$ consists of equivalence classes of parametrized unit speed geodesics rays $\gamma: [0,\infty)\rightarrow \X$ up to equivalence of being at bounded Hausdorff distance. As $G$ acts on $\X$ by isometries, this induces an action of $G$ on $\vis\X$. We say that a non-zero vector $v \in \T_P\X$ \emph{points towards} $x\in \vis \X$ if the geodesic ray with initial velocity $v$, denoted $\gamma_{P,v}$, is in the class of $x$ in $\vis\X$.

We now consider the point $\ell_0\in \vis \X$ corresponding to the following ray $\gamma_t:[0,\infty) \rightarrow \X$:
\begin{align}\label{ModelGeodesic}
     \gamma_t = \left(\diag(e^t,0, \cdots, 0, e^{-t})\cdot P_0\right). 
\end{align}

The stabilizer of $\ell_0$ in $G$ is exactly the stabilizer of the isotropic line $\langle e_1\rangle$, which we will denote by $P_1$ (cf. \cite[Proposition 10.64]{BH99}). Hence the $G$-orbit of $\ell_0$ in $\vis \X$ is naturally identified with the space of isotropic lines $G/P_1=\Ein^{p-1,p}$. From now on we will therefore view $\Ein^{p-1,p}$ as a subset of the visual boundary.

We use the terminology that a non-zero vector $v \in \T_P\X$ \emph{points towards} $\Ein^{p-1,p}$ if the geodesic ray $\gamma_{P,v}$ has the property that $[\gamma_{P,v}]$ lies in the $G$-orbit of $\ell_0$ in $\vis\X$. We can make this abstract property of interest completely concrete with a simple geometric criterion. 

\begin{proposition}[Pointing Towards $\Ein^{p-1,p}$] \label{Prop:PointingTowardsEinstein}
Let $\phi \in \T_P\X\simeq \Hom(P, P^\perp)$. Viewed as a map $\phi: P \rightarrow P^\bot$, then $\phi$ points towards $\Ein^{p-1,p}$ if and only if $\mathrm{rank}(\phi) = 1$. Moreover, in this case $\phi$ points towards $\ell =\graph(\phi_{\mid L})$ where $L\subset P$ is the orthogonal to the kernel of $\phi$. 
\end{proposition}

\begin{proof}
Suppose $\phi \in \T_P\X$ has rank one. 
Up to the action of $G=\SO_0(p,p+1)$, we can assume that $P$ is the basepoint of $\X$ from \eqref{X_Basepoint}. First note that $\Stab_{G}(P)\cong\SO(p)\times \SO(p+1)$, which preserves $P\in \X$, acts transitively on the space of rank one elements $\phi\in \Hom(P, P^\perp)$, up to positive scalars. Hence up to positive scalars and the $K_P$-action, we can assume that $\phi=\dot{\gamma}(0)$ is the derivative of the geodesic ray \eqref{ModelGeodesic} pointing towards $\ell_0$, which has rank one and by definition points towards $\Ein^{p-1,p}$. Thus, if $\phi$ has rank one, we conclude it does point towards $\Ein^{p-1,p}$. 
Note finally that for the model geodesic, $\ell_0=\graph(\phi_{\mid L})$ where $L\subset P$ is the line orthogonal to the kernel of $\phi$. By $G$-equivariance, the `moreover' statement follows.

We now prove the `only if' statement. Since $K_P$ acts transitively on $\Ein^{p-1,p}$, for every element $\ell\in \Ein^{p-1,p}$, there is a rank one tangent vector $\phi \in \T_P\X$ that points towards $\ell$. As a general fact, for each point $P\in \X$ in the symmetric space $\X$ and all point $\ell \in \vis \X$ there is a unique unit vector in $\T_P\X$ pointing towards $\ell$ \cite{Ebe96}. Hence if $\phi\in \Hom(P,P^\perp)$ does not have rank $1$, it does not point towards $\Ein^{p-1,p}$.\end{proof}

There is an important geometric consequence of Proposition \ref{Prop:PointingTowardsEinstein}:  a realization of $\Ein^{p-1,p}$ as a fiber bundle, which describes the embedding of this flag manifold in the visual boundary.

\begin{proposition}[Fiber Bundle for $\Ein^{p-1,p}$]\label{Prop:EinFiberBundle}
Fix $P \in \X$. The orthogonal projection map $\pi_P:\Ein^{p-1,p} \rightarrow \Gr_{1}(P)$ defines an $\mathbb{S}^{p}$-fiber bundle. 
\end{proposition}

\begin{proof}
Fix $\ell \in \Ein^{p-1,p}$. Choose $u \in Q_+(\pi_P(\ell))$. Then we may write $\ell = [u+z]$ uniquely for some $z \in Q_-(P^\bot)$. Working backwards, this means 
\[ \pi^{-1}([u])=\Ein([u]\oplus P^\bot)\cong \Ein^{0,p} \cong \sphere^{p}. \]
With $u$ fixed, the sphere $Q_-(P^\bot)\cong \mathbb{S}^{p}$ is identified with the fiber $\pi^{-1}([u])$ by the map $z\mapsto [u+z]$. The map $\pi$ is a surjective submersion, as a $K_P$-equivariant map for $K_P=\Stab_{G}(P)$. By compactness and the Ehresmann fibration lemma, the result follows. 
\end{proof}

Proposition \ref{Prop:EinFiberBundle} shows $\Ein^{p-1,p}$ is an $\mathbb{S}^p$-fiber bundle over $\RP^{p-1}$. 
We shall see the fiber $\mathfrak{F}_{\rho}$ of the cocompact quotients $M_{\rho}:= \rho(\pi_1S)\backslash \Omega_{\rho}$ of $p$-Anosov $\iota$-Fuchsian deformations $\rho:\pi_1S\rightarrow G$ interact nicely with this fibration $\pi_P$. Indeed, in Section \ref{Sec:StructuralResult}, we show $\mathfrak{F}_{\rho}$ is realized as an $\mathbb{S}^{p-2}$-fiber subbundle of the fibration $\pi_P$ for appropriate $P \in \X$. 

Beyond the relation with domains of discontinuity, the fibration in Proposition \ref{Prop:EinFiberBundle} is intimately linked with Proposition \ref{Prop:PointingTowardsEinstein}. The following remark clarifies this point.
\begin{remark}
Let $\ell \in \Ein^{p-1,p}$ and $P \in \X$. The unique tangent vector $\phi \in \T^1_P\X$, up to positive scalars, that points points towards $\ell$ is exactly described by Proposition \ref{Prop:EinFiberBundle}. Indeed, write $\ell =[u+z]$ for $u \in Q_+(\pi_P(\ell))$ and $ z\in Q_-(P^\bot)$, unique up to replacing $(u,z)$ by $(-u,-z)$. The tangent vector $\phi:P\rightarrow P^\bot$ is the unique rank one map satisfying $\phi(u)=z$. 
\end{remark}

\subsection{Hitchin Representations}\label{Subsec:Hitchin}

Let us first define Fuchsian-Hitchin representations. Let $\Gamma=\pi_1S$ be the fundamental group of a closed surface $S$ of genus $g \geq 2$, and recall that a Fuchsian representation is a discrete and faithful representation of $\Gamma$ into $\PSL(2,\R)$.

\begin{definition}
   A Fuchsian-Hitchin representation $\rho:\Gamma\to \SO_0(p,p+1)$ is the composition of a Fuchsian representation $\rho_0:\Gamma\to \PSL(2,\R)$ through the unique irreducible representation $\eta:\PSL(2,\R)\to \SO_0(p,p+1)$ up to conjugation. 
\end{definition}

The representation $\eta$ is well-known, and is explicitly described in (the proof of) \cite[Lemma 5.8]{CTT19}. More generally Hitchin representations can be defined as follows:

\begin{definition}
   A representation $\rho:\Gamma\to \SO_0(p,p+1)$ is Hitchin if it can continuously be deformed to a Fuchsian-Hitchin representation.
\end{definition}

Hitchin representations are particular examples of an interesting class of representations: they are \emph{Anosov}. Labourie proved that they are Anosov with respect to every parabolic subgroup, so in particular they are $p$-Anosov \cite{Lab06}. We give here a short definition of the Anosov property for completeness, but we will not work with it directly. For the definition, given a transformation $g \in \SO_0(p,p+1)$, we denote 
$\sigma_1(g)\geq  \sigma_2(g) \geq \cdots \geq \sigma_{2p+1}(g)$ as the singular values of $g$,  i.e., the eigenvalues of $(gg^t)^{1/2}$.

\begin{definition}[{\cite{KLP17,BPS}}]\label{Defn:Anosov}
A representation $\rho:\Gamma\to \SO(p,p+1)$ is \emph{$P_p$-Anosov}, or simply $p$-Anosov, if there exist some constants $c,d>0$ such that for all $\gamma\in \Gamma$,
\[ \frac{\sigma_p}{\sigma_{p+1}}(\rho(\gamma))\geq ce^{d\lvert\gamma\rvert}, \]
where $\lvert \gamma\rvert$ is the norm of $\gamma$ for the word metric relative to a finite generating set of $\Gamma$. 
\end{definition}

The notion of singular values depends on the choice of a (Euclidean) scalar product on $\R^{p,p+1}$. Definition \ref{Defn:Anosov} is independent of this choice, and of choice of generating set of $\Gamma$ used to define the word metric. Note also that if the scalar product is well-chosen the singular values of elements of $\SO_0(p,p+1)$ have the property that $\sigma_i=-\sigma_{2p+1-i}$, and $\sigma_{p+1}=1$. \medskip

In the paper, we will only use the fact that Anosov representations admit \emph{limit maps}. Let $\partial \Gamma$ be the Gromov boundary of $\Gamma$, which is topologically a circle and carries a natural $\Gamma$-action, and let $ \Iso_p(\R^{p,p-1})$ be the space of isotropic $p$ planes in $\R^{p,p+1}$.

\begin{theorem}[{\cite{GGKW17}}]\label{thm:AnosovlimitMap}
Let $\rho:\pi_1S\rightarrow \SO_0(p,p+1)$ be a $p$-Anosov representation. There exists a unique continuous $\rho$-equivariant map $\xi^p: \partial \Gamma \rightarrow \Iso_p(\R^{p,p+1})$ that sends the unique attracting fixed point $\gamma^+ \in \partial\Gamma$ of $\gamma\in \Gamma$ to the unique attracting fixed isotropic $p$-plane of $\rho(\gamma)$ in $\Iso_p(\R^{p,p+1})$.
\end{theorem}

\subsection{Domains of Discontinuity via Tits Metric Thickening}\label{Subsec:MetricThickening} 

Next, we recall how the relevant cocompact domain of discontinuity $\Omega \subset \Ein^{p-1,p}$ is defined. This domain is built by the general construction by Kapovich-Leeb-Porti \cite{KLP18}, by \emph{Tits metric thickening}, which was first described in \cite{GW12} in this case.  

For each point $T \in \Iso_p(\R^{p,p-1})$, we define the  \emph{thickening} $K_T \subset \Ein^{p-1,p}$ as in \cite{GW12} by:
\begin{align}\label{Thickening_Containment}
    K_T= \{\ell \in \Ein^{p-1,p} \mid \ell \subset T\}.
\end{align}
Note that $\RP^{p-1} \cong K_T \subset \Ein^{p-1,p}$ is a projective $(p-1)$-plane in $\Ein^{p-1,p}$.

Now, the domain of discontinuity of interest, denoted $\Omega^{\text{Thick}}_{\rho} \subset \mathcal{F}_0$, is obtained by removing the thickening of the entire limit set $\Lambda = \image(\xi^p)$:
\begin{align}\label{Omega_Thick}
    \Omega_{\rho}^{\text{Thick}}:= \Ein^{p-1,p} \backslash \bigcup_{x \in \partial \Gamma}K_{\xi^p(x)}.
\end{align}
We may write $\Omega_{\rho}$ for  $\Omega^{\Thick}_{\rho}$. This domain interacts pleasantly with the $\rho(\Gamma)$-action:

\begin{theorem}
\label{thm:GuichardWienhardEin}
Let $\Gamma$ be a word hyperbolic group and $\rho: \Gamma \rightarrow \SO_0(p,p+1)$ a $p$-Anosov representation. The domain \eqref{Omega_Thick} is a cocompact domain of proper discontinuity for $\rho(\Gamma)$. 
\end{theorem}

This result was originally proven in {\cite[Proposition 8.1, Theorem 8.6]{GW12}}, and also follows by \cite[Theorem 1.8]{KLP18}. Note that a direct dimension count shows $\Omega_{\rho}$ is non-empty in the case $\Gamma = \pi_1S$ is a surface group.

\subsection{Domains of Discontinuity via Bases of Pencils}  

We now describe a seemingly different way to build domains of discontinuity: with \emph{bases of pencils of tangent vectors}. We will then show this approach yields a more geometric realization of the domain \eqref{Omega_Thick} from the previous subsection.

\begin{definition}[Pencil]
For any $x\in \X$, we call a ($k$-)plane $\mathcal{P}\subset \T_x\X$ a ($k$-)\textbf{pencil of tangent vectors} or ($k$-)\textbf{pencil} for short. 
\end{definition}

When specification is not provided, we shall always write `pencil' to mean a 2-pencil. For most of the work, we consider only this case. 

A pencil $\mathcal{P}$ defines naturally a subset of the flag manifold $\Ein^{p-1,p}$ of expected codimension two, that we call the \emph{$\tau$-base}. 
For the following definition, recall that for a tangent vector $v \in \T^1_x\X$, we use the notation $\gamma_{x,v}(\infty)$ to denote the class $[\gamma_{x,v}] \in \vis\X$ of $\gamma_{x,v}$.

\begin{definition}[Base of Pencil]\label{Defn:TauBasePencil}
Let $\mathcal{P} \subset \T_{x} \mathbb{X}$ be a pencil. Then the \textbf{base of} $\boldsymbol{\mathcal{P}}$, denoted $\mathcal{B}(\mathcal{P})$, is given by 
\[ \mathcal{B}(\mathcal{P}) = \{ \gamma_{x, v}(\infty) \in \Ein^{p-1,p}\subset \vis \X\mid v \in \T_x\X, \; v\, \bot \,\mathcal{P} \}\]
\end{definition}

In other words, the base $\mathcal{B}(\mathcal{P})$ contains the points in $\Ein^{p-1,p}$ that can be reached in $\vis\mathbb{X}$ by traveling from $x$ via directions orthogonal to $\mathcal{P}$ in $\T_x\X$.\medskip 
\begin{remark}
One can more generally define bases of pencils for other flag manifolds, viewed as orbits in the visual boundary. In the present paper, the bases of pencil considered will always be in the Einstein universe $\Ein^{p-1,p}$.
\end{remark}

As in \cite{Dav25}, it is useful to distinguish a notion of \emph{regularity} of a pencil. Presently, we consider only a single notion of regularity, related to the $G$-orbit $\Ein^{p-1,p}\subset \vis \X$.

\begin{definition}[Ein-Regularity]\label{Defn:EinRegular}
(i) A pencil $\mathcal{P} \subset \T_P\X$ is \textbf{$\mathbf{Ein}$-regular}, or just regular, when all non-zero elements $\phi\in \mathcal{P}$, viewed as elements of $\Hom(P,P^\perp)$, have rank $p$.

(ii) A (non-trivial) representation $\iota:\SL(2,\R)\rightarrow \SO_0(p,p+1)$ is \textbf{$\mathbf{Ein}$-regular} when its sub-symmetric space $\Ha^2_{\iota}$ has tangent pencil $\T_P\Ha^2_{\iota}$ that is $\Ein$-regular for every $P \in \Ha^2_{\iota}$. 
\end{definition}

In fact, Definition \ref{Defn:EinRegular} (i) is just a slight modification of \cite[{Definition 5.6}]{Dav25}, as clarified by the following proposition. 

\begin{proposition}[$\Ein$-regularity, Lie-theoretically]
\label{Prop:EinRegularity}
A tangent vector $\phi \in \T_P\X$ is $\Ein$-regular if and only if it is $\tau$-regular, where $\tau=\diag(1,0,\dots,0,-1) \in \overline{\mathfrak{a}}^+$ in the sense that its Cartan projection $\mu(\phi)$ satisfies $\langle\mu(\phi), w\cdot \tau\rangle\neq 0$ for all $w \in W$ in the Weyl group $W$.
\end{proposition}

\begin{proof}
Regard $\phi$ as of type $\phi: P \rightarrow P^\bot$. Form the matrix $A_\phi \coloneqq \begin{pmatrix} 0 & -\phi^{*q} \\ \phi & 0 \end{pmatrix} \in \mathfrak{so}(p,p+1)$, in block form relative to $\R^{p,p+1} = P\oplus P^\bot$. Observe that $2\,\rank(\phi) = \rank(A_\phi)$. 
The conclusion follows from the fact that $\phi \in \T\X$ has $\rank(A_\phi) = 2p$ if and only if the Cartan projection
\begin{align}\label{CartanProj}
    \mu(\phi) = (\mu_1, \dots, \mu_p, 0, -\mu_p, \dots, -\mu_1),
\end{align}
where $\mu_1\geq \mu_2\geq \cdots\geq \mu_p\geq 0$, satisfies $\mu_p > 0$. This is equivalent to having $\langle\mu(\phi), w\cdot \tau\rangle\neq 0$ for all $w$ in the Weyl group.
\end{proof}

The point of regularity is to ensure that the base of pencil is well-behaved. 
\begin{proposition}\label{Prop:SmoothBaseOfPencil}
If $\mathcal{P} \subset \T_P\X$ is an $\Ein$-regular pencil, then the base $\mathcal{B}(\mathcal{P})\subset \Ein^{p-1,p}$ is a smooth codimension two submanifold. 
\end{proposition}

A proof of Proposition \ref{Prop:SmoothBaseOfPencil} is given in greater
generality in \cite[Lemma 6.7]{Dav25}. We re-prove the result in the present context later on in Lemma \ref{Lem:EinBaseTopology} in a more direct way.

A fact we shall vitally rely on is that the base of a pencil is invariant under regular deformation. 
\begin{lemma}[Base of Pencil Invariance]\label{Lem:RegularPencilDeformation}
Suppose that $(\mathcal{P}_t)_{t \in [0,1]} \subset \T_P\X$ is a smooth family of pencils that is $\Ein$-regular for all $t$. Then $\mathcal{B}(\mathcal{P}_0)$ and $\mathcal{B}(\mathcal{P}_1)$ are diffeomorphic. 
\end{lemma}

Lemma \ref{Lem:RegularPencilDeformation} is proven in \cite[Corollary 6.8]{Dav25}, though perhaps it is useful to say here that the proof is essentially  due to the Ehresmann fibration lemma and Proposition \ref{Prop:SmoothBaseOfPencil}.\medskip 

Next, we recall how the notion of bases of pencils relates to fibrations of cocompact domains of discontinuity $\Omega^{\tau}_{\rho}$. 

Let $f: \tilde{S} \rightarrow \mathbb{X}$ be a totally geodesic embedding that is $\Ein$-regular.
Fixing an arbitrary basepoint $o \in \X$, we can define a domain $\Omega_{f}\subset \Ein^{p-1,p}$ using Busemann functions by 
\begin{align}\label{Omega_Busemann}
    \Omega_{f} := \{a \in \mathcal{F}_{\tau} \mid b_{a,o} \circ f\; \text{is proper, bounded below} \}.
\end{align}
Here, the \emph{Busemann function} $b_{a,o}$ measures the relative distance of points $x \in \X$ to $a \in \vis\X$ from the point of view of $o$ by 
\[ b_{a,o}(x) := \lim_{t \rightarrow \infty} d_{\X}(\gamma_{o,a}(t),x)-t.\]
By the triangle inequality, the definition of $b_{a,o}$ is well-defined. Busemann functions $b_{a,o}$ are well-known to be smooth when $\X$ is a symmetric space. 

There is a natural projection from the Busemann domain $\Omega_{f}^{\tau}$ to the universal cover $\tilde{S}$ of $S$ as follows. Here, the projection $\pi: \Omega_f \rightarrow \tilde{S}$ maps $a$ to the unique point $x \in \tilde{S}$ such that $b_{a,o} \circ f$ has a critical point at $f(x)$. This critical point is unique by \cite[Lemma 7.2]{Dav25}. 

\begin{lemma}[Nearest Point Projection]\label{Lem:NearestPointProjection}
Let $f: \tilde{S} \rightarrow \mathbb{X}$ be totally geodesic and $\Ein$-regular. Then: 
\begin{enumerate}[noitemsep, label=(\roman*)]\item $\Omega_{f}$ is open, 
    \item $\pi$ is a fibration.
    \item $(\Omega_f)|_x = \mathcal{B}(\mathcal{P}_x)$, where $\mathcal{P}_x \subset \T_{f(x)}\mathbb{X}$ is the pencil $ df(\T_{x}\tilde{S})$. 
\end{enumerate}
\end{lemma} 

\begin{proof}
The definition of $\pi$ is well-defined and $\Omega_f$ is open by \cite[Lemma 7.2]{Dav25}. Then \cite[Theorem 7.3]{Dav25} settles points (2) and (3).   
\end{proof}

In fact, the domain \eqref{Omega_Busemann} defined via Busemann functions is the same as the domain \eqref{Omega_Thick} defined via Tits metric thickening. It is through this link that we can find the fibers of interest as a base of pencil. 
For this connection to be possible, we must first relate Einstein regularity to the Anosov property. The idea is well-known (cf. \cite[Proposition 4.7]{GW12}). 
\begin{proposition}
\label{prop:Anosov for Ein Regular representations}
Let $\iota:\SL(2,\R) \rightarrow \SO_0(p,p+1)$ be a representation. If $\iota$ is $\Ein$-regular, then for any Fuchsian representation $\rho:\pi_1S\rightarrow \SL(2,\R)$, the composition $\iota\circ \rho$ is $p$-Anosov. 
\end{proposition}

We now apply the nearest point projection Lemma to $\Ein$-regular Fuchsian representations, and state the link between the two domains of discontinuity of interest. 

\begin{proposition}[{\cite[Theorem 7.11]{Dav25}}]
Let $\rho: \pi_1S\rightarrow \SL(2,\R) \stackrel{\iota}{\hookrightarrow} \SO_0(p,p+1)$ be an $\Ein$-regular $\iota$-Fuchsian representation and $f: \tilde{S} \hookrightarrow \X$ the corresponding totally geodesic $\rho$-equivariant map. The fibration $\pi:\Omega_f\to \widetilde{S}$ is $\rho$-equivariant and thus descends to a fiber bundle projection $\rho(\pi_1S)\backslash \Omega_{f} \rightarrow S$. The domains $\Omega_\rho^{\Thick}$ in \eqref{Omega_Thick} and $\Omega_f$ in \eqref{Omega_Busemann} coincide. 
\end{proposition}

The diffeomorphism type of the quotient $M_{\rho} = \rho(\pi_1S)\backslash \Omega_{\rho}$ is locally constant in deformations that remain $p$-Anosov. In particular, all Hitchin representations have the same smooth quotients $M_{\rho}$. This notion of invariance of topology originates in \cite[Theorem 9.12]{GW12}. An appropriate version applying in the present context is given in \cite[Corollary 7.15]{Dav25}. 

\begin{corollary}[Topological Invariance]\label{cor:FiberInvariance}
Let $\rho:\pi_1S\rightarrow \SO_0(p,p+1)$ be a $p$-Anosov deformation of an $\iota$-Fuchsian. The diffeomorphism type of the quotient $M_{\rho}=\rho(\pi_1S)\backslash \Omega_{\rho}$ and the homotopy type of the fiber $\mathfrak{F}_{\rho}$ depend only on $\iota$.  
\end{corollary}

\subsection{\texorpdfstring{$\Ein$}{Ein}-regular \texorpdfstring{$\sllie_2\R$}{sl(2,R)}-subalgebras in \texorpdfstring{$\solie(p,p+1)$}{so(p,p+1)}}

We now describe all representations $\sllie_2\R\rightarrow \solie(p,p+1)$ that give rise to $\Ein$-regular representations $\SL(2,\R)\rightarrow\SO_0(p,p+1)$ in the sense of Definition \ref{Defn:EinRegular} (ii).

Recall that an $\sllie_2\C$-representation in $\sllie_{2p+1}\C$ associated to an integer partition of $2p+1$ factors through $\solie_{2p+1}\C$ if and only if the even parts occur with even multiplicity \cite[{Chapter 5.1}]{CM17}. In the real case, additional information is needed \cite[Chapter 9]{CM17}, however, for the integer partitions of interest, this information is essentially determined.

\begin{proposition}[$\Ein$-regular $\sllie_2\R \hookrightarrow \solie(p,p+1)$]
\label{Prop:EinRegularSL2}
Let $I$ be an integer partition of $2p+1$ with one odd part and all even parts of even multiplicity.
\begin{enumerate}[noitemsep, label=(\roman*)]
    \item For some positive integers $q, r_1, \dots, r_k$, the partition $I$ can be written as 
    \begin{align}\label{SpecialPartition}
     2p+1=(2q+1)+ 2r_1 + 2r_1+ \cdots + 2r_k + 2r_k. 
    \end{align}
    \item There is a unique corresponding $\Ein$-regular representation $\iota_I:\sllie_2\R \hookrightarrow \solie(p,p+1)$ up to conjugation in $\mathrm{O}(p,p+1)$, with an indecomposable invariant splitting 
        \[ \R^{p,p+1} = \R^{q,q+1}\oplus V_1 \oplus V_1' \oplus \cdots \oplus V_k \oplus V_k',\]
    where $V_i$ and $V_i'$ are dual isotropic $2r_i$-planes, and all other subspaces are  orthogonal. 
    \item Every $\Ein$-regular representation $\SL(2,\R)\rightarrow \SO_0(p,p+1)$ has an indecomposable invariant splitting of the form in (ii). 
\end{enumerate}
\end{proposition}

\begin{proof}
Note that (i) is obvious. 

(ii) Following \cite[{Theorem 9.34}]{CM17}, given a partition $I$ with an odd part and all even parts with even multiplicity, there is a unique representation $\iota_I$ up to conjugation in $\solie(p,p+1)$. More precisely such representation are parametrized by \emph{signed} Young diagrams, and the Young diagram associated with these partitions admit a unique such decoration. We will show this unique representation satisfies the desired conditions by constructing it directly.
\medskip

Write $I$ as in (i). 
The representation $\iota_I$ can be seen as built from the unique $\sllie_2\R$-representation in $\sllie_{2p+1}$ associated to this integer partition. We show how to place an invariant quadratic form of signature $(p,p+1)$, so that the representation factors through $\solie(p,p+1)$, is $\Ein$-regular, and has the desired invariant subrepresentations. 

We split $\R^{2p+1} = \R^{2q+1} \oplus \bigoplus_{i=1}^k(V_{i} \oplus V_{i}')$, where $V_{i}$ and $V_{i}'$ are $2r_i$-dimensional subspaces, such that $\iota_I $ has $\R^{2q+1}, V_{1}, V_1', \dots, V_{k},V_{k}'$ as irreducible sub-representations. By \cite[$\S$ 9.3]{CM17}, the representation $\iota_I$ preserves a quadratic form of signature $(q,q+1)$ on $\R^{2q+1}$. We can also define an $\iota_I$-invariant non-degenerate pairing of $V_{i}, V_{i}'$ so that $(V_{i} \oplus V_{i}') \cong \R^{2r_i,2r_i}$ is a sub-representation. It follows that $\iota$ takes values in $\solie(p,p+1)$. 

Now, to see that $\iota_1$ is $\Ein$-regular, we take the Cartan decomposition of $\g = \image(\iota)$. Write $\g = \frakk \oplus \mathfrak{p}$. By Proposition \ref{Prop:EinRegularity} then suffices to see that every element $X \in \mathfrak{p}$ has rank $2p$. However, examining the $\sllie_2\R$-triple $\{E,F,H\}$ associated to $\image(\iota_I)$, the rank condition is verified  by checking the sub-representations $\R^{q,q+1}$ and $(V_{i} \oplus V_{i}')$. For the sub-representation  $\iota_q:\sllie_2\R \rightarrow \solie(q,q+1)$, which is principal, each such tangent vector $X|_{\R^{q,q+1}}$ has rank $2q$ since $\iota_q$ is $\Ein$-regular. For the sub-representations in $\solie(2r_i,2r_i)$, one verifies each tangent vector $X$ has $X|_{\R^{2r_i,2r_i}}$ with full rank equal to $4r_i$. In total, this verifies $\Ein$-regularity.

(iii) Let $\iota_I: \SL(2,\R)\rightarrow \SO_0(p,p+1)$ be $\Ein$-regular. 
Note that for every odd part of the integer partition $I$, the semisimple element $H$ of the $\sllie_2$-triple has a one dimensional kernel on the corresponding odd dimensional sub-representation. Hence, $H$ has rank at most $2p-k$ if the partition $I$ has $k$ odd parts. We conclude $\iota_I$ is $\Ein$-regular if and only if it has exactly one odd part and the rest even (with even multiplicity). The claim then follows by (ii). 
\end{proof}

We finish this subsection with some examples of integer partitions of interest.

\begin{example}\label{Ex:p=5}
Let $p= 5$. The integer partitions associated to $\Ein$-regular $\sllie_2\R$-subalgebras in $\solie(5,6)$ are as follows:
\begin{itemize}[noitemsep]
    \item $11$,
    \item $7+2+2$,
    \item $4+4+3$,
    \item $3+2+2+2+2$.
\end{itemize}
\end{example}

\begin{example}\label{Ex:p=6}
Let $p= 6$. The integer partitions associated to $\Ein$-regular $\sllie_2\R$-subalgebras in $\solie(6,7)$ are as follows:
\begin{itemize}[noitemsep]
    \item $13$,
    \item $9+2+2$,
    \item $6+6+1$,
    \item $5+4+4$,
    \item $5+2+2+2+2$,
    \item $4+4+2+2+1$,
    \item $2+2+2+2+2+2+1$. 
\end{itemize}
\end{example}

\section{Uniformizing Higgs Bundles}\label{Sec:Uniformizing}

The goal of this section is to describe Higgs bundles that uniformize $\Ein$-regular $\iota$-Fuchsian representations $\rho:\pi_1S\rightarrow \SO_0(p,p+1)$. 
By Proposition \ref{Prop:EinRegularSL2}, such a representation $\rho$ is reductive and decomposes 
into sub-representations of two types: Fuchsian-Hitchin in $\SO_0(q,q+1)$ for $q < p$ and 
reducible $r+r$ representations in $\SO_0(r,r)$ for  $r$ even. To uniformize the desired representation $\rho$, we
then need only describe how to uniformize these two types of sub-representations, which we do in two separate subsections below.

We will use this framework solely to understand the associated totally geodesic $\rho$-equivariant map into the symmetric space $f:\widetilde{\Sigma}=\widetilde{S}\to \X$.
The material in this section requires knowledge of $\SO(p,q)$-Higgs bundles. We include the bare essentials and refer the reader to \cite{AA09, ABCBGPO19} for further details on the subject. We also assume basic familiarity with the non-abelian Hodge (NAH) correspondence, as detailed in \cite{Col19, Li19} (among many such survey articles). 

\subsection{Fuchsian-Hitchin Representations in \texorpdfstring{$\SO_0(q,q+1)$}{SO(q,q+1)}}\label{Subsec:FuchsianHitchin}

We now describe uniformizations of Fuchsian-Hitchin representations $\rho:\pi_1S\rightarrow \SO_0(q,q+1)$. 

To describe a Higgs bundle associated to $\rho$ via the non-abelian Hodge correspondence, fix a Riemann surface $\Sigma= (S,J)$ on $S$ and form the holomorphic rank $(2p+1)$-vector bundle 
\begin{align}\label{Ep_Fuchsian}
    \mathcal{E}_p \coloneqq \bigoplus_{i=-p}^{p} \K^i,
\end{align} 
where $\K =\K_{\Sigma}= (\T^{1,0}\Sigma)^*$ is the holomorphic cotangent line bundle. We define a holomorphic endomorphism valued-one form $\Phi \in H^0(\End(\V)\otimes \K)$ as follows:
\[ \Phi   = \left( \K^{p} \stackrel{1}{\longrightarrow} \mathcal{K}^{p-1} \stackrel{1}{\longrightarrow}\cdots \stackrel{1}{\longrightarrow} \mathcal{O} \stackrel{1}{\longrightarrow}  \cdots \stackrel{1}{\longrightarrow} \K^{1-p}  \stackrel{1}{\longrightarrow} \K^{-p} \right ). \]
Here, in this diagram, each element $1$ is a sub-tensor of $\Phi$, namely some holomorphic endomorphism valued one-form. For example, $1 \in H^0(\Hom(\K^{p}, \K^{p-1})\otimes \K)$ makes sense because $\Hom(\K^{p}, \K^{p-1})\otimes \K\cong \mathcal{O}$ is a holomorphically trivial line bundle.

The pair $(\V,\Phi)$ will be the $\SO_0(p,p+1)$ Higgs bundle of interest, once we endow it with further structure $(\mathcal{U},\mathcal{V},Q, \omega)$. 
To this end, we first split $\V=\mathcal{U}\oplus \mathcal{V}$ into two parts:
\begin{align}
    \mathcal{U}&=\mathcal{K}^{p-1}\oplus \mathcal{K}^{p-3}\oplus  \cdots \oplus \K^{3-p}  \oplus \mathcal{K}^{1-p}. \label{U_Subbundle} \\
 \mathcal{V}&= \mathcal{K}^{p}\oplus \mathcal{K}^{p-2}\oplus \cdots  \oplus \K^{2-p}\oplus \mathcal{K}^{-p}. \label{V_Subbundle}
\end{align}
Note that $\rank(\mathcal{U})=p$ and $\rank(\mathcal{V})=p+1$. 
We then define a holomorphic symmetric bilinear form $Q= Q_{\mathcal{U}}\oplus (-Q_{\mathcal{V}})$ on $\V$ respecting this splitting $\mathcal{U} \oplus \mathcal{V}$ by letting each of $Q_{\mathcal{U}}$ and $Q_{\mathcal{V}}$ be the natural dual pairings. Explicitly, 
\[ Q_{\mathcal{V}}= Q_{p,-p}+Q_{p-2,2-p}+\cdots + Q_{p,-p},\] 
where each sub-tensor $Q_{i,-i}$ is the dual pairing $Q_{i,-i}:\K^{i} \otimes \K^{-i}\rightarrow \mathcal{O}$. Then $Q_{\mathcal{U}}$ is defined completely analogously. Finally, we may denote $\omega =\underline{1}\in \det(\V) \cong \mathcal{O}$ as the `obvious' volume form on $\V$. In fact, $\omega = \omega_{\mathcal{U}}\wedge \omega_\mathcal{V}$, for $\omega_{\mathcal{U}}, \, \omega_{\mathcal{V}}$ the natural volume forms on  $\mathcal{U}, \mathcal{V}$. 

The non-degenerate bilinear form $Q$, and the volume form $\omega$ together reduce the structure group of $\V$ to $\SO(2p+1,\C)$. Furthermore, the splitting  $\mathcal{E}=\mathcal{U} \oplus \mathcal{V}$ along with $(\omega_{\mathcal{U}}, \omega_{\mathcal{V}})$ reduce the structure group further to $\SO(p,\C)\times \SO(p+1,\C) = K^\C$, where $K \cong \SO(p)\times \SO(p+1)$ is the maximal compact subgroup of $G$.

The Higgs field $\Phi$ is compatible with all the structures imposed. Indeed, we may write $\Phi = \varphi -\varphi^{*Q}$, where $\varphi \in H^0(\Hom(\mathcal{U},\mathcal{V})\otimes \K)$ is given by $\varphi = \Phi|_{\mathcal{U}}$, to see that $\Phi$ is traceless and satisfies $\Phi \in \Omega^0(\solie(Q)\otimes \K)$. As explained in \cite[Proposition 3.10]{Col19}, the tuple $(\V, \Phi, \mathcal{U}, \mathcal{V}, Q, \omega) $ defines an $\SO_0(p,p+1)$-Higgs bundle. Note also that \cite[Section 8.5]{AA09} shows this Higgs bundle corresponds under NAH to a Fuchsian-Hitchin representation $\rho:\pi_1S\rightarrow \SO_0(p,p+1)$. \medskip 

This Higgs bundle carries a distinguished Hermitian metric $h$. The condition distinguishing $h$ is the following: we demand $\nabla := \nabla_{\delbar, h} + \Phi+\Phi^{*h}$ is \emph{flat}, where $\delbar$ is the Dolbeault operator on the underlying smooth bundle $E$, $\nabla_{\delbar, h} $ is the Chern connection of the Hermitian holomorphic bundle $(E,\delbar,h)$, and $\Phi^{*h} \in \Omega^{0,1}(\End(\V))$ is the adjoint of $\Phi$ with respect to $h$. Such a Hermitian metric $h$ is unique in this case, which follows from \emph{stability} of the Higgs bundle. The metric $h$ is \emph{diagonal} under the splitting $\V =\bigoplus_{i=p}^{-p}\K^i$, and can be even written down explicitly, though this is not needed presently.  

Now, the connection $\nabla$ has holonomy in $\SO_0(p,p+1)$ due to the compatibility of $\Phi$. This entails that $\nabla$ preserves a real sub-bundle $\V^{\R}$, with fibers pointwise isomorphic to $\R^{p,p+1}$, which is the fixed point set of an anti-holomorphic involution $\lambda: \V \rightarrow \V$. The involution $\lambda$ relates $Q$ and $h$. Indeed, $h(\cdot , \cdot) = Q(\cdot, \lambda(\cdot))$ (see e.g. \cite[Section 2.3]{CTT19}). Thus, on the real locus $\V^\R$, we simply have $h|_{\mathcal{E}^\R}= 
Q|_{\mathcal{U}^\R} \oplus (-Q|_{\mathcal{V}^\R})$. Concretely, the conjugation $\lambda$ is described in a basis $(e_k)_{k=p}^{-p}$ that is $h$-unitary with $e_k\in \mathcal{K}^k$ by:
\begin{equation}
    \label{eq:RealLocusHitchin}
    \lambda=\big[(z_p, \,z_{p-1},\cdots,\, z_{1-p},\,z_{-p})\longmapsto (\overline{z_{-p}},\,\overline{z_{1-p}},\,\cdots ,\,\overline{z_{p-1}},\,\overline{z_{p}} )\big].
\end{equation}  

The $\rho$-equivariant totally geodesic map into the symmetric space $f:\widetilde{\Sigma}=\widetilde{S}\to \X$ can be understood through the Higgs bundle as follows. The pullback bundle $\pi^*\mathcal{E}^\R$, under the universal covering $\pi:\widetilde{\Sigma} \rightarrow \Sigma$, can be identified using the flat connection with $\widetilde{\Sigma} \times \R^{p,p+1}$. With this identification, $f(x)=P$ where 
\begin{align}
    P&={\left(\mathcal{E}^\R\cap \mathcal{U}\right)}_{\mid x}\subset \R^{p,p+1} \\
    P^\perp&={\left(\mathcal{E}^\R\cap \mathcal{V}\right)}_{\mid x}\subset \R^{p,p+1}. 
\end{align}

Let $\varphi$ be the section of $\mathcal{K}\otimes \Hom(\mathcal{U},\mathcal{V})$ such that $\varphi-\varphi^{*Q}=\Phi$. Explicitly, 

\[\varphi: \bigoplus_{0\leq j\leq p-1} \mathcal{K}^{p-1-2j}\stackrel{1}{\longrightarrow}\mathcal{K}^{p-2-2j}. \]

Since $h|_{\V^{\R}} = Q|_{\mathcal{U}^\R} \oplus (-Q|_{\mathcal{V}^\R})$ and $\Psi:=\Phi + \Phi^{*h}$ is $h$-self adjoint, we find $\Psi$ is $Q$-anti-self-adjoint, meaning $\Psi \in \Omega^{1}(\solie(Q))$. In fact, $\Psi$ is \emph{real}: it preserves the real locus $\V^\R$. Through this identification, actually $\Psi$ corresponds to the differential $df$, up to a constant multiplicative factor (cf. \cite{Gui18, Li19}). 
More precisely, take $x_0\in \widetilde{S}$ and $P_0=f(x_0)$, and $v \in \T_{x_0}\tilde{S}$. 
The differential $\mathrm{d}f(v)\in T_{P_0}\X$ is identified with $(\Phi+\Phi^{*h})(\pi(v))$. We shall use this dictionary between the Higgs bundle and the harmonic map repeatedly. Hence, we record the following remark. 
\begin{remark}[Pencils in Higgs Bundles]\label{Remk:HiggsPencils}
Let $f: \widetilde{\Sigma} \rightarrow\X$ be the unique $\rho$-equivariant harmonic map associated to a Higgs bundle $(\V,\Phi)$. Let $x_0 \in \tilde{S}$, denote $p_0 = \pi(x_0) \in S$, $P_0=f(x_0)$.
\begin{enumerate}[label=(\roman*)]
    \item (Tangent Pencil) The tangent plane $df(\T_{x_0}\tilde{S}) \subset \T_{P_0}\X$ can be identified with the plane $\mathscr{P}$ of endomorphisms of the fiber $\V^\R|_{p_0}$ given by: 
\begin{align}\label{HiggsPencil}
 \mathscr{P} =\{ (\Phi+\Phi^{*h})(v) \in \End(\V^\R|_{p_0}) \mid v \in \T_{p_0}S  \}.
\end{align}
The endomorphisms $(\Phi+\Phi^{*h})(v)$ each obtain the form $\eta-\eta^{*Q}$ for some map $\eta \in \Hom(\mathcal{U}^\R, \mathcal{V}^\R)$ and thus the pencil $\mathscr{P}$ can be treated as an object of type $\mathfrak{E} \in \Gr_2\big(\Hom_{\R}(\mathcal{U}^\R|_{p_0}, \mathcal{V}^\R|_{p_0})\big)$:

\begin{align}\label{HiggsPencilbis}
\mathfrak{E}=\{ \varphi-(\varphi^{*Q})^{*h}(v) \in \Hom_{\R}(\mathcal{U}^\R|_{p_0}, \mathcal{V}^\R|_{p_0})  \mid v \in \T_{p_0}S  \}.
\end{align}

\item (General Pencil) In fact, if $\mathcal{P} \subset \T_{P_0}\X$ is any pencil, then there is a real two-plane $\mathfrak{E}(\mathcal{P}) \in \Gr_2\big(\Hom_{\R}(\mathcal{U}^\R|_{p_0}, \mathcal{V}^\R|_{p_0})\big)$ corresponding to $\mathcal{P}$. 
\end{enumerate}
\end{remark}

\subsection{\texorpdfstring{$r+r$}{r+r} Fuchsians in \texorpdfstring{$\SO_0(r,r)$}{SO(r,r)} }\label{Subsec:Uniformize_p+p}

We now describe the uniformizing Higgs bundle $(\V_r, \Phi, Q)$ for a representation $\rho:\pi_1S\rightarrow \SO_0(r,r)$ that factors through a pair of Fuchsian-Hitchin sub-representations in $\SL(r,\R)$. The descriptions below work regardless of parity of $r$, but later on we shall consider only $r$ even. 

Now, consider the following holomorphic vector bundle $\V_{r,r}$ and a holomorphic bilinear form $Q$ on $\V_{r,r}$: 
\begin{align}
    \V_{r,r} &= \bigoplus_{n=-r+1}^{r-1} (\K^{n/2} \oplus \K^{-n/2}) \label{SO(r,r)Bundle}\\ 
    Q &= \bigoplus_{n=-r+1}^{r-1} (-1)^{n-1}\begin{pmatrix} 0 & 1 \\ 1& 0 \end{pmatrix}
\end{align}
As in the Fuchsian-Hitchin case one can equip this bundle with a Hermitian metric $h$ such that $\nabla := \nabla_{\delbar, h} + \Phi+\Phi^{*h}$ is a flat connection. 

Here, we split $\V_{r,r} = \mathcal{U} \oplus \mathcal{V}$ as follows: 
\begin{align}\label{SO(p,p)UV}
    \mathcal{U} = \bigoplus_{\stackrel{n \in [-r+1, r-1]}{n \in (2\Z+1)}}(\K^{n/2} \oplus \K^{-n/2}) \\
    \mathcal{V} = \bigoplus_{\stackrel{n \in [-r+1, r-1]}{n \in 2\Z}}(\K^{n/2} \oplus \K^{-n/2})
\end{align}
That is, the 2-dimensional sub-bundles $\K^{n/2} \oplus \K^{-n/2}$ alternate between $\mathcal{U}$ and $\mathcal{V}$. This choice forces the following Higgs field $\Phi_{r,r} \in H^0(\End(\V_{r,r})\otimes \K)$ to be compatible with the real form $\lambda$ determined by $Q$ and $h$ as in the Fuchsian-Hitchin case:
\begin{align*}
   \K^{\frac{-1+r}{2}} &\stackrel{1}{\longrightarrow} \K^{\frac{-3+r}{2}} \stackrel{1}{\longrightarrow} \cdots \stackrel{1}{\longrightarrow} \K^{\frac{+1-r}{2}}\\  
   \K^{\frac{+1-r}{2}} 
   &\stackrel{1}{\longleftarrow} \K^{\frac{+3-r}{2}} \stackrel{1}{\longleftarrow} \cdots \stackrel{1}{\longleftarrow} \K^{\frac{-1+r}{2}}.
\end{align*}
Here, each vertical pair of dual line bundles comprises one of the two-dimensional sub-bundles in \eqref{SO(r,r)Bundle}. 
Note that this Higgs bundle $(\V, \Phi, Q)$ is polystable as an $\SL(2n,\C)$-Higgs bundle, and decomposes as the sum of two uniformizing Fuchsian $\SL(r,\R)$-Higgs bundles. By the choice of the real structure encoded by $Q$, we see that $\Phi \in \solie(\End(\V_r),Q)$ is compatible with the real form, i.e. lies in $\p^\C$, and hence we have the desired $\SO_0(r,r)$-Higgs bundle. 

\section{\texorpdfstring{$(\SO_0(p, p+1), \Ein^{p-1,p})$}{(SO(p,p+1),Ein(p-1,p))}-Geometric Structures}\label{Sec:EinGeometricStructures}

In this section, we prove the main results: the topology of the fibers and the global topology of the $(\SO_0(p,p+1),\Ein^{p-1,p})$-manifolds $M_{\rho}$ for $\rho:\pi_1S\rightarrow \SO_0(p,p+1)$ a $p$-Anosov deformation of an $\iota$-Fuchsian. 

We now outline the strategy for these proofs. 
Fix $p \geq 3$. 
By Corollary \ref{cor:FiberInvariance}, it suffices to consider when $\rho$ is an $\iota$-Fuchsian. In this case, let 
$\Ha^{2}_{\iota}$ denote the sub-symmetric space of the associated $(\mathrm{P})\SL(2,\R)$-subgroup defined by $\iota$. Choose any $x \in \Ha^2_{\iota}$. 
By Lemma \ref{Lem:NearestPointProjection}, the fiber $\mathfrak{F}_{\rho}$ of the bundle $\pi:M_\rho \rightarrow S$ defined via nearest point projection is the base of pencil  
$\mathcal{B}(\mathcal{P}) \subset \Ein^{p-1,p}$ for the tangent pencil $
\mathcal{P} = \T_x\Ha^2_{\iota}$.
Moreover, by \cite{Dav25}, a global version of this base of pencil construction describes the global quotient $M_{\rho}$. 
To determine the smooth isomorphism type of $\pi:M_{\rho} \rightarrow S$, there are four main steps.  

\begin{enumerate}[label=(\roman*)]
    \item (Structure Lemma) First, we prove a general structural for about $\Ein$-regular pencils. Namely, if $\mathcal{P} \subset \T_x\X$ is $\Ein$-regular, then its base $\mathcal{B}(\mathcal{P})$ is (a $\Z_2$-quotient of) a sphere bundle $\sphere(E)$ of a rank $(p-1)$ vector bundle $E \rightarrow \mathbb{S}^{p-1}$. Classifying the topology of the fiber is then reduced to classifying the vector bundle $E=E(\mathcal{P})$. 
    \item (Deformation) Next, we prove an invariance result: the associated vector bundle $E(\mathcal{P})$ is the same for all tangent pencils $\mathcal{P}_{\iota} = \T_x\Ha^2_{\iota}$, independent of the $\Ein$-regular representation $\iota$. 
    We prove this by showing all such pencils $\mathcal{P}_{\iota}$ can be deformed in an $\Ein$-regular fashion to the simplest possible $\Ein$-regular pencil. In fact, this simplified pencil $\mathcal{P}_{\iota_0}$ is the tangent pencil to $\Ha^{2}_{\iota_0}$ for $\iota_0$ a specific $(\mathrm{P})\SL(2,\R)$-subgroup, which happens depends on parity of $p$. Indeed, $\iota_0$ is associated by Proposition \ref{Prop:EinRegularSL2} to the integer partition $2+2+\cdots+2+1$ when $p$ is even and to $3+2+2+\cdots +2$ when $p$ is odd. 
    \item (Fiber) Using some elementary geometric reasoning, we determine that $E=E(\mathcal{P}_{\iota_0})$ is the trivial vector bundle $E=\sphere^{p-1}\times \R^p$ when $p$ is odd and $E = \T\sphere^{p-1}$ when $p$ is even. This leads to the fiber $\Ein^{p-1,p-2}$ for $p$ odd and $\T^1\RP^{p-1}$ for $p$ even. 
    \item (Global Topology) Finally, by viewing the deformation in (ii) as happening fiberwise inside the Higgs bundle associated to the representation, we determine the global topology of a global ``Higgs bundle base of pencil'' inside the flat $\Ein^{p-1,p}$-fiber bundle $\tilde{S}\times_{\rho} \Ein^{p-1,p}$, which serves as a diffeomorphic model for the quotient $M_{\rho}$ as a fiber bundle over $S$. 

\end{enumerate} 
In order to execute this plan, the steps are presented slightly out of order, for reasons we now explain. In $\S$\ref{Sec:StructuralResult}, we prove the structure lemma (i). Next, we describe the simplified fibers $\mathcal{B}(\mathcal{P}_{\iota_0})$ in $\S$\ref{Sec:Fiber_Odd}, \ref{Sec:Fiber_Even}, completing (iii). Vitally, $\S$\ref{Sec:Fiber_Odd}, \ref{Sec:Fiber_Even} set up all necessary input data needed to functorially identify the fiber, which later aids in our execution of (iv). 
In $\S$\ref{Sec:HiggsPencilDeformation}, we explain the deformation of the tangent pencils from $\T_{x}\Ha^2_{\iota}$ to $\T_x\Ha^2_{\iota_0}$, using Higgs bundles, handling (ii). This deformation section is a bridge. We introduce bundle versions of the auxiliary data from $\S$\ref{Sec:Fiber_Odd}, \ref{Sec:Fiber_Even}, in preparation to identify moving bases of pencils. Finally, in  $\S$\ref{Sec:GlobalTopology_Odd}, \ref{Sec:GlobalTopology_Even}, we use the identifications of the fiber in a smoothly varying way to study the global Higgs bundle base of pencil and determine the isomorphism type of $M_{\rho}$ as a smooth fiber bundle over $S$, finishing (iv).

\subsection{Geometry of the \texorpdfstring{$\Ein$}{Ein}-Base}\label{Sec:StructuralResult}

In this section, we prove a key structural result:
base $\mathcal{B}(\mathcal{P})$ of an $\Ein$-regular pencil $\mathcal{P}\subset \T_x\X$ is the sphere bundle of an associated rank $(p-1)$-vector bundle $E\rightarrow \RP^{p-1}$. This construction occurs for a fixed pencil $\mathcal{P}$ at a fixed point $x \in \X$. We then extend this identification to the bundle setting, for a moving base of pencil.  

\subsubsection{A Single Base of Pencil}
The core idea to determine the base $\mathcal{B}(\mathcal{P})$ of a pencil $\mathcal{P} \subset \T_{P} \mathbb{X}$ is given by the following lemma. 

\begin{lemma}[Einstein Base of Pencil]\label{Lem:EinBaseTopology}
Let $\mathcal{P} \subset \T_{P} \mathbb{X}$ be an $\Ein$-regular $d$-pencil for $2 \leq d \leq p$, i.e., a $d$-dimensional subspace whose non-zero elements are $\Ein$-regular. Write $\underline{P^\bot} \rightarrow Q_+(P)$ for the trivial vector bundle $Q_+(P) \times P^\bot$ over $Q_+(P)$. Then:
\begin{enumerate}[noitemsep,label=(\roman*)]
    \item $\mathcal{P}$ yields a trivial $d$-dimensional vector sub-bundle $\mathcal{R} \rightarrow Q_+(P)$ of  $\underline{P^\bot}$ with fiber 
    \[\mathcal{R}_u \coloneqq \{ \psi(u) \mid \psi \in \mathcal{P} \} .\] 
    \item The base $\mathcal{B}(\mathcal{P})$ is diffeomorphic to $Q_-(\mathcal{R}^\bot)/\sim,$ 
where $\mathcal{R}^\bot \rightarrow Q_+(P)$ is the orthogonal complement of $\mathcal{R}$ in $\underline{P}^\bot$ and $(u,v) \sim (-u,-v)$. 
    \item The projection map $\pi_P:\Ein^{p-1,p} \rightarrow \mathbb{P}(P)$ restricts to $\mathcal{B}(\mathcal{P})$ as a surjective submersion that realizes $\mathcal{B}(\mathcal{P})$ as an $\mathbb{S}^{p-2}$ fiber bundle over $\mathbb{P}(P)$. 
\end{enumerate}
\end{lemma}

\begin{proof}
(i) Suppose that $\mathcal{P} \subset \T_P\X$ is an $\Ein$-regular $d$-pencil. By Proposition \ref{Prop:EinRegularity}, we have that $\dim \mathcal{R}_u = \dim \mathcal{P}=d$. Any basis $(\psi_i)_{i=1}^d$ for $\mathcal{P}$ yields a global frame $(s_i)_{i=1}^d$ for $\mathcal{R}$ given by $u \stackrel{s_i}{\longmapsto} \psi_i(u)$. 

(ii) The idea rests entirely on Proposition \ref{Prop:PointingTowardsEinstein}. Any line $\ell \in \Ein^{p-1,p}$ obtains the form $\ell= [u+v]$ for $u \in Q_+(P), \, v \in Q_-(P^\bot)$ for a unique pair of elements $(u,v), (-u,-v)$. The antipodal pair $\pm (u,v)$ determines the unique rank one linear map $X_{\ell, P}: P \rightarrow P^\bot$ such that $u \mapsto v$ and $\ker(X_{\ell, P})\perp u$. By Proposition \ref{Prop:PointingTowardsEinstein}, the unique geodesic $\gamma: [0,\infty) \rightarrow \X$ with $\gamma(0)=P$, $\dot{\gamma}(0) =X_{\ell,P}$ has $\gamma(\infty)=\ell \in \vis\X$. Observe that if $X_{\ell, P}$ is such a map and $\psi \in \T_P\X$, then $X_{\ell, P} \,\bot\, \psi$ if and only if $\psi(u) \bot v$, because of the shape of the Riemannian metric on $\Hom(P,P^\perp)\simeq \T_P\X$.

We conclude by the previous argument that $\ell \in \mathcal{B}(\mathcal{P})$ if and only if $v \bot \psi(u)$ for all $\psi \in \mathcal{P}$. The desired claim (ii) follows. 

(iii) For every $[u] \in \mathbb{P}(P)$, there is a point $\ell = [u+z] \in \mathcal{B}(\mathcal{P})$ by selecting $z \in \mathcal{R}^\bot|_{u}$. Hence, $\pi_P$ surjects. One easily verifies $\pi_P$ is a submersion at $\ell$ using (ii). Then (iii) follows by the Ehresmann fibration lemma. 
\end{proof}

Lemma \ref{Lem:EinBaseTopology} (ii) says that $\mathcal{B}(\mathcal{P})$ is, up to a $\mathbb{Z}_2$-quotient, the total space $\sphere(E)$ of a sphere bundle of a rank $(p-1)$ vector bundle $E \rightarrow \mathbb{S}^{p-1}$. Alternatively, (iii) says $\mathcal{B}(\mathcal{P})$ is an $\mathbb{S}^{p-2}$-fiber-subbundle of the $\mathbb{S}^{p}$-fiber bundle realization of $\Ein^{p-1,p}$ over $\RP^{p-1}=\mathbb{P}(P)$ from Proposition \ref{Prop:EinFiberBundle}. 

The following corollary emphasizes the generality of the result of Lemma \ref{Lem:EinBaseTopology}. 

\begin{corollary}
Let $\iota:\SL(2,\R) \rightarrow \SO_0(p,p+1)$ be any $\Ein$-regular representation. Let $\rho:\pi_1S\rightarrow \SO_0(p,p+1)$ an $\iota$-Fuchsian representation. Then $\rho$ is $p$-Anosov and the fibers of $M_\rho = \rho(\pi_1S)\backslash \Omega_{\rho}$ are diffeomorphic to an $\mathbb{S}^{p-2}$-fiber bundle over $\RP^{p-1}$. 
\end{corollary}

\begin{proof}
Note that $\rho$ is $p$-Anosov by Proposition \ref{prop:Anosov for Ein Regular representations}.
By Lemma \ref{Lem:NearestPointProjection}, the fibers of $M_\rho$ are diffeomorphic to $\mathcal{B}(\mathcal{P})$, where $\mathcal{P}$ is any tangent pencil to the sub-symmetric space of $\image(\iota)$. By Lemma \ref{Lem:EinBaseTopology}, the result follows. 
\end{proof}

\subsubsection{Moving Bases of Pencils}
We now describe a bundle analogue of Lemma \ref{Lem:EinBaseTopology}. To this end, we introduce some notation for the statement. 
We write the discussion below in general for maps from a $d$-manifold $N$. When equivariance is needed, we reintroduce the hypothesis that $N=\widetilde{\Sigma}$.

Let $f: N \rightarrow \X$ be a smooth map and $\mathscr{P} \in \Omega^0(N, \Gr_2(f^*\T\X))$ a smooth pencil along $f$ that is pointwise $\Ein$-regular. We consider the moving bases of pencils defined by $\mathscr{P}$ assembled in to a smooth manifold: 
\begin{align}\label{MovingBaseOfPencil}
    \mathcal{B}({\mathscr{P}}) =\{(p, \ell) \in N \times \Ein^{p-1,p} \mid \ell \in \mathcal{B}(\mathscr{P}_p)\}. 
\end{align}
We introduce the bundle analogue of $\mathcal{R}$ as we vary the point $P \in \X$ using the map $f$. 

Let $\mathscr{T} \rightarrow \X$ denote the tautological vector bundle, with fiber $\mathscr{T}|_{P} = P \subset \R^{p,p+1}$. We form the $\mathbb{S}^{p-1}$-bundle $f^*(Q_+(\mathscr{T}))$. We may then define a vector bundle $\mathscr{R}^\bot \rightarrow f^*(Q_+(\mathscr{T}))$ by 
\begin{align}\label{RBotBundle}
    \mathscr{R}^\bot = \{ (p,u, z) \in N\times Q_+(\R^{p,p+1}) \times \R^{p,p+1} \mid u \in Q_+(f(p)), \,z \in f(p)^\bot, \, z\in \mathscr{P}|_p(u)^\bot\}. 
\end{align}
Here, we view $\mathscr{P}_p(u) = \{\psi(u) \in f(p)^\bot \mid \psi \in \mathscr{P}_p\}$. 
In other words, in \eqref{RBotBundle}, for every choice of $P = f(p)$, we form the same construction as in Lemma \ref{Lem:EinBaseTopology} fiberwise. 

\begin{lemma}[Moving Ein-Base]\label{Lem:MovingBase}
Let $f: N\rightarrow \X$ and $\mathscr{P} \in \Omega^0(N, f^*\T\X))$ be smooth with $\mathscr{P}$ also $\Ein$-regular. Form the bundle \eqref{RBotBundle}. The manifold $\mathcal{B}(\mathscr{P})$ in \eqref{MovingBaseOfPencil} is diffeomorphic to the $\Z_2$-quotient of $Q_-(\mathscr{R}^\bot)$ by the antipodal map $a(p,u,z) = (p,-u,-z)$. 
\end{lemma}

\begin{proof}
Consider the map $F: Q_-(\mathscr{R}^\bot) \rightarrow \mathcal{B}(\mathscr{P})$ by 
$F(p,u,z) = (p, [u+z])$, well-defined by Lemma \ref{Lem:EinBaseTopology}. $F$ is evidently a fiber bundle map lifting $\id:N\rightarrow N$. The map $F$ is a smooth 2-1 covering by Lemma \ref{Lem:EinBaseTopology}, which after quotienting by $a$, becomes a diffeomorphism. 
\end{proof}

We now show that if we form a homotopy of pencils $(\mathscr{P}_t)_{t \in [0, 1]}$ along the same map $f: N \rightarrow \X$, which is $\Ein$-regular for all $t$, then the associated moving bases of pencils $\mathcal{B}(\mathscr{P}_t)$ are diffeomorphic as fiber bundles.

\begin{corollary}[Homotopy of Moving Bases of Pencils]\label{Cor:HomotopyBasesPencil}
Let $f:N \rightarrow \X$ be smooth and $(\mathscr{P}_t)_{t \in [0, 1]}$ a smoothly varying 2-pencil along $f$ that is pointwise $\Ein$-regular for all $t$. Then $\mathcal{B}(\mathscr{P}_1)$ is isomorphic to $\mathcal{B}(\mathscr{P}_0)$ as a smooth fiber bundle over $N$.
\end{corollary}

In fact, we show the stronger fact that $\mathcal{B}(\mathscr{P}_1)$ is isomorphic to $\mathcal{B}(\mathscr{P}_0)$ as a smooth $\mathbb{S}^{p-1}$-fiber bundle over $f^*(\mathbb{P}\mathscr{T})$.
 
\begin{proof}
The result is an immediate application of the ideas of Lemma \ref{Lem:MovingBase}, after we introduce some more notation. 

Let $\underline{\R^{p,p+1}} = \X \times \R^{p,p+1}$ be the trivial bundle. Denote $\mathscr{T} \rightarrow \X$ as the tautological vector bundle, i.e. $\mathscr{T}|_P =P$, such that $\underline{\R^{p,p+1}}= \mathscr{T} \oplus \mathscr{T}^\bot$. 
We denote $E = f^*Q_+(\mathscr{T})$ for the $\mathbb{S}^{p-1}$-bundle over $N$. 
Let $\pi: E\rightarrow N$ denote the fiber bundle projection. We consider $(\pi \circ f)^*(\mathscr{T}^\bot)$. Then we have vector bundles
$\mathscr{R}_t\rightarrow E $ given by 
\[ \mathscr{R}_t = \{ (p,u,z) \in N \times Q_+(\R^{p,p+1}) \times\R^{p,p+1} \mid u \in Q_+(f(p)) ,\, z \in \mathscr{P}_t|_p(u) \}.\]
Let $t$ be fixed. Then $\mathscr{P}_t|_p \in \Gr_2(\Hom(f(p), f(p)^\bot))$. Observe the vector bundle splitting $(\pi \circ f)^*(\mathscr{T}^\bot) =\mathscr{R}_t\oplus \mathscr{R}_t^\bot $ over $E$. Concretely, for $(p,u)\in N\times Q_+(f(p))$, we split 
\[ f(p)^\bot = \mathscr{R}_t|_{(p,u)}\oplus \mathscr{R}_t^\bot|_{(p,u)}.\] 

Now, we view the bundles $(\mathscr{R}_t)$ as sub-bundles of a fixed bundle $(\pi \circ f)^*(\mathscr{T}^\bot)$, varying smoothly in $t$. Hence $\mathscr{R}_0 \cong \mathscr{R}_t$ for all $t$, as smooth vector bundles over $E$, see \cite[{Proposition 1.7}]{Hat17}. 

Hence, $\mathscr{R}_0^\bot \cong \mathscr{R}^\bot_t$ for all $t$, which induces smooth fiber bundle isomorphisms $Q_-(\mathscr{R}_0^\bot) \cong Q_-(\mathscr{R}_t^\bot)$ as $\mathbb{S}^{p-1}$-fiber bundles over $E$. Now, the $\Z_2$-quotient of $Q_-(\mathscr{R}^\bot_t)$ by the antipodal map $a$ as realizes $\mathcal{B}(\mathscr{P}_t)$ by Lemma \ref{Lem:MovingBase}. Under the quotient by $a$, the fiber bundle isomorphism $Q_-(\mathscr{R}_0^\bot) \cong Q_-(\mathscr{R}_t^\bot)$ descends to a fiber bundle isomorphism $\mathcal{B}(\mathscr{P}_t) \cong \mathcal{B}(\mathscr{P}_0)$ as $\mathbb{S}^{p-1}$-fiber bundles over $\mathbb{P}(E) = f^*(\mathbb{P}(\mathscr{T}))$. The claim follows immediately.  
\end{proof}

In the case of interest, $N = \widetilde{\Sigma}$. We now briefly remark on equivariant deformations.

\begin{remark}\label{Remk:EquivariantPencil}
Following the proof of Corollary \ref{Cor:HomotopyBasesPencil}, there is no issue to upgrade the identifications above to be $\pi_1S$-equivariant as long as the pencils $\mathscr{P}_t$ are each $\rho$-equivariant for the same representation $\rho: \pi_1S \rightarrow \SO_0(p,p+1)$. 
\end{remark} \medskip

We now consider the problem of quotients. 
Suppose $f:\widetilde{\Sigma}\rightarrow \X$ is a $\rho$-equivariant harmonic map. We will use the Higgs bundle $(\V,\Phi)$ associated to $\rho$ via the non-abelian Hodge correspondence to describe the quotient $\rho(\pi_1S)\backslash \mathcal{B}(\mathscr{P})$ for $\mathscr{P}$ a smooth $\rho$-equivariant $\Ein$-regular pencil along $f$. By Corollary \ref{Cor:HomotopyBasesPencil}, it does not matter which such pencil $\mathscr{P}$ is used.

To this end, we introduce a Higgs bundle version of base of pencil 
inside the flat $\Ein^{p-1,p}$ fiber bundle $\Ein(\V^\R) \rightarrow \Sigma$. 
As in Remark \ref{Remk:HiggsPencils}, associated to $\mathscr{P}$ is an object $\mathfrak{E}$ of type $\mathfrak{E}\in\Omega^0(\Sigma, \Gr_2( \Hom(\mathcal{U}^\R, \mathcal{V}^\R))$. 
We define its associated base of pencil in analogous fashion to Lemma \ref{Lem:EinBaseTopology} as follows: 
\begin{align}\label{HiggsBundleBasePencil}
 \mathcal{B}(\mathfrak{E})  \coloneq \{ (p,\ell) \in \Ein(\V^\R) \mid\ell \in \mathcal{B}(\mathscr{P}|_p) \subset \Ein(\V^\R|_p)\}.
\end{align}
In other words, the fiber of this fiber bundle at $p \in S$ is given by 
\[ \mathcal{B}(\mathfrak{E})|_p= \{ [u+v] \in \Ein(\V^\R|_p)\mid u \in Q_+(\mathcal{U}^\R), \; v \in Q_-(\mathcal{V}^\R)_, \;v \in \mathscr{P}|_p(u)^\bot\}.\] 
By the same circle of ideas, we note the following identification: the Higgs bundle pencil \eqref{HiggsBundleBasePencil} is a geometric representative for the $\pi_1S$-quotient of the moving base of pencil \eqref{MovingBaseOfPencil}.

\begin{remark}[Pencils: Upstairs+Downstairs]\label{remk:PencilUp+Down}
By the definition of $\mathcal{B}(\mathfrak{E})$, we have an isomorphism of smooth fiber bundles between $\mathcal{B}(\mathfrak{E})$ and $\rho(\pi_1S)\backslash \mathcal{B}(\mathscr{P})$.
\end{remark}

We can use these Higgs bundle versions of bases of pencils for explicit descriptions of the quotient of the domain of discontinuity. 
\begin{corollary}\label{Cor:HiggsBundlePencil}
Let $f: \widetilde{\Sigma} \rightarrow \X$ be a $\rho$-equivariant harmonic map for $\rho:\pi_1S\rightarrow \SO_0(p,p+1)$ and $(\mathscr{P}_t)_{t \in [0,1]}$ a smooth $\rho$-equivariant pencil along $f$ that is $\Ein$-regular for all $t$. Then 
\begin{enumerate}[label=(\roman*), noitemsep] 
    \item \emph{(Invariance of Moving Pencils)}. $\pi_1S\backslash \mathcal{B}(\mathscr{P}_t)$ is isomorphic to $\mathcal{B}(\mathfrak{E}_0)$ as a smooth fiber bundle over $S$ for all $t$, where $\mathfrak{E}_0$ is the Higgs bundle pencil associated to $\mathscr{P}_0$.
    \item \emph{(Higgs Pencil = Global Quotient)}. Suppose that $\rho$ is a $P_p$-Anosov $\iota$-Fuchsian representation and $\mathscr{P}_0$ is the tangent pencil to $f$. Then the quotient $M_{\rho} =\rho(\pi_1S)\backslash \Omega_{\rho}$ is isomorphic to $\mathcal{B}(\mathfrak{E}_1)$ as a smooth fiber bundle.
\end{enumerate}
\end{corollary}

\begin{proof}
Note that (i) is a direct consequence of Remark \ref{Remk:EquivariantPencil}. 

(ii) By Lemma \ref{Lem:NearestPointProjection} and Remark \ref{remk:PencilUp+Down}, we have $M_{\rho} \cong \rho(\pi_1S)\backslash \Omega_{\rho} \cong \mathcal{B}(\mathfrak{E}_0)$. Now, (ii) follows immediately from (i).
\end{proof}

\subsection{Description of the fiber, \texorpdfstring{$p$}{p} odd}\label{Sec:Fiber_Odd}

Let $p =2k+1$ be odd. 
In this section, we define a certain kind of \emph{reducible pencil} $\mathcal{P}_0\subset \T_P\X$, and then parametrize the base of this pencil, showing it is diffeomorphic to $\Ein^{p-1,p}$. These reducible pencils $\mathcal{P}_0$ are generalizations of the tangent pencil $\T_P\Ha^2_{\iota_0}$ for $\iota_0: \SL(2,\R)\rightarrow \SO_0(p,p+1)$ the representation associated to the integer partition $2p+1=3+\underbrace{2+2+\dots +2}_{p-1}$. 
This section contains the heart of the proof that the fiber of $M_{\rho}\rightarrow S$ for $\rho$ a deformation of (any) $\iota$-Fuchsian is $\Ein^{p-1,p}$. 
\medskip

\begin{figure}[ht]
    \centering
    \includegraphics[width=0.55\linewidth]{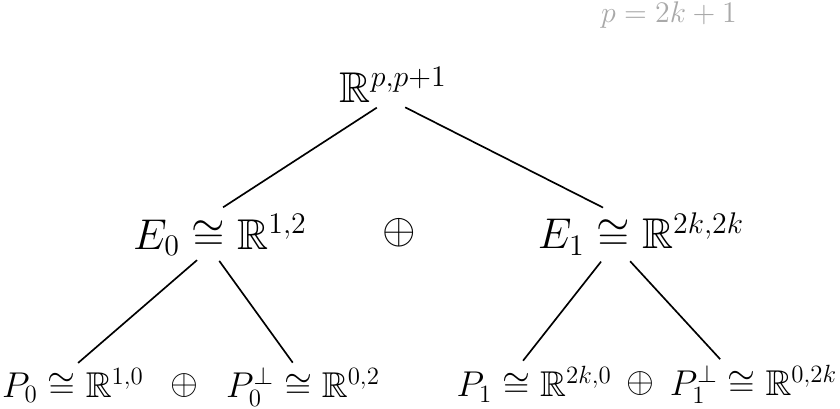}
    \caption{Schematic summary of the relevant splittings for $\S$4.2.}
    \label{Fig:OddSplitting}
\end{figure}

We now introduce various decompositions of the ambient vector space $\R^{p,p+1}$ that are essential to the definition of a reducible pencil.
\begin{itemize}
    \item Let $E_0\oplus E_1=\R^{p,p+1}$ be an orthogonal splitting with $E_0\cong \R^{1,2}$ and $E_1\cong \R^{p-1,p-1}$.
    \item Let $P\in \X$ be a spacelike subspace of $\R^{p,p+1}$ of dimension $p$ that decomposes as $P=P_0\oplus P_1$, where $P_0\subset E_0$ is a spacelike line and $P_1\subset E_1$ is a spacelike $(p-1)$-plane. We decompose similarly $P^\perp=P_0^\perp\oplus P_1^\perp$ where $P^\perp_0\subset E_0$ and $P^\perp_1\subset E_1$.
    \item  Finally fix complex structures $J_1:P_1^\perp\to P_1^\perp$, $J'_1:P_1\to P_1$, and $J_0:P_0^\perp \rightarrow P_0^\bot$ that are $q$-orthogonal. 
\end{itemize}
Hence, we may identify $P_1$ and $P_1^\perp$ with $\C^k$ for $k=\frac{p-1}{2}$. 
For simplicity, we write $J=J_0\oplus J_1\oplus J'_1$ for the complex structure on $P_0^\perp\oplus E_1$. See Figure \ref{Fig:OddSplitting}.

We now introduce the technical definition of interest for pencils in the odd case. 

\begin{definition}\label{Defn:ReduciblePencilOdd}
A pencil $\mathcal{P}_0\subset \T_P\X\simeq \Hom_{\R}(P,P^\perp)$ is \textbf{reducible with data} $\mathbf{(E_0,J)}$ if it lies in $\Hom(P_0,P_0^\perp)\oplus \Hom(P_1,P_1^\perp)$, is $\Ein$-regular, and satisfies the following properties:
\begin{itemize}
    \item[(i)] $(\psi_0+\psi_1)\in \mathcal{P}_0\iff (J_0\circ\psi_0+J_1\circ\psi_1) \in \mathcal{P}_0$,
    \item[(ii)] for all $(\psi_0+\psi_1)\in \mathcal{P}_0$, $\psi_1:P_1\to P_1^\perp$ is holomorphic; i.e. $\psi_1\circ J_1'=J_1\circ\psi_1$,
\end{itemize}
\end{definition}
Here we denote $\psi_0\in \Hom(P_0,P_0^\perp)$ and $\psi_1\in \Hom(P_1,P_1^\perp)$. Note that condition (i) asserts $\mathcal{P}_0$ is a complex line in $\Hom(P,P^\bot)$, i.e., it is closed under post-composition $J_{P^\bot}$. \medskip 

With all the background data fixed, we show a reducible pencil has a canonically identified base of pencil. Below, associated to $P \in \X$, we write $h =h(P)$ for the Euclidean metric $h=q|_P+(-q|_{P^\bot})$ on $\R^{2p+1}$. 
\begin{lemma}[Geometry of Simplified Pencil, Odd Case]\label{lem:FiberOdd}
Let $p =(2k+1) \geq 3$ be odd and $\mathcal{P}_0 \subset \T_P\X$ be a reducible pencil with data $(E_0,J)$. Fix any unit spacelike element $e_0 \in Q_+(P_0)$. 
\begin{enumerate}[label=(\roman*)]
    \item There is a natural real vector bundle isomorphism $\Theta: \underline{\mathcal{P}_0} \otimes_{\C} \underline{P_1} \rightarrow \mathcal{R}^\bot$, depending only on $(E_0,J, \mathcal{P}_0, e_0)$, given by: 
    \[ \Theta: \underline{\mathcal{P}_0} \otimes_{\C} \underline{P_1} \rightarrow \mathcal{R}^\bot, \;\;\; \big(u,\psi \otimes_{\C} v\big)\longmapsto \bigg(u, \,\psi(v)-\frac{\langle \psi(v),\psi(u)\rangle_h}{\langle \psi(u)- 2\psi(e_0),\psi(u) \rangle_h }(\psi(u)-2\psi(e_0))\bigg).\] 
Here, $\underline{\mathcal{P}_0}\otimes_{\C} \underline{P_1}=Q_+(P)\times (\mathcal{P}_0\otimes_\C P_1)$ is a trivial bundle and $\mathcal{R}^\perp \rightarrow Q_+(P)$ is as in Lemma \ref{Lem:EinBaseTopology}. 
\item Moreover, under the isomorphism $\underline{\mathcal{P}_0} \otimes_{\C} \underline{P_1} \cong \underline{P_1^\perp}$ via the map $(\psi,u)\mapsto \psi(u)$, we obtain a natural $\R$-linear vector bundle isomorphism depending only on $(E_0,J,\mathcal{P}_0,e_0)$:
\[ \widehat{\Theta}: \underline{P_1^\perp} \rightarrow \mathcal{R}^\bot.\]
\item For each tuple $(E_0,J,\mathcal{P}_0,e_0)$, there is a natural diffeomorphism the between base of pencil $\mathcal{B}(\mathcal{P}_0)$ and $\Ein(P \oplus P_1^\bot) \cong \Ein^{p-1,p-2}$.
\end{enumerate}
\end{lemma}

\begin{remark}
Here what we mean by natural is that the diffeomorphisms in (i)-(iii) are defined for every quadruple $(E_0,J,\mathcal{P}_0,e_0)$, and depend smoothly on this input data. In particular, we can consider bundle versions of these maps. 
\end{remark}

 \begin{proof}
(i) We now check that $\Theta$ is well-defined and produces the desired vector bundle isomorphism. First, note that the fiberwise formula for $\Theta$ is $\R$-linear in $\psi$ and $\R$-linear in $v$. Next, we find that 
the formula for $\Theta(u,\psi, v)$ agrees on the inputs $(u, J\psi, v)$ and $(u, \psi, Jv)$.
Indeed, this follows from holomorphicity of the pencil, namely properties (i) and (ii) from Definition \ref{Defn:ReduciblePencilOdd}, which imply 
$\psi(J(v)) = J(\psi(v)) $ for $\psi \in \mathcal{P}_0$ and $v \in P_1$, as well as $(J\psi)(v)=J( \psi(v))$. Hence, the map $\Theta$ is well-defined on the domain $\underline{\mathcal{P}_0}\otimes_{\C} \underline{P_1}$ as desired. 
Next, we observe by direct calculation:
\[ \langle \Theta(u, \psi\otimes_{\C} v), \psi(u)\rangle_h = 0.\] 
 This  implies that $\Theta(u, \psi \otimes_{\C} v) \in \mathcal{R}^\bot|_u$. Indeed, 
 $\mathcal{R}^\bot|_u = P^\bot \cap \mathcal{P}_0(u)^\bot = P^\bot \cap \C \{ \psi(u) \}^\bot$ for any $\psi \in \mathcal{P}_0$ by (i) of Definition \ref{Defn:ReduciblePencilEven}. We conclude that $\Theta$ is a well-defined $\R$-linear vector bundle morphism. 

 By dimension count, it suffices to show that $\Theta$ is injective. For $\psi \neq 0$, note that $\psi(u)-2\psi(e_0)\notin P^\perp_1$. Indeed, we may write $u \in Q_+(P)$ as $u=\lambda e_0+u_1$ with $\lambda\in \R$ and $u_1\in P_1$. Since $e_0\in Q_+(P)$, one has $|\lambda|\leq 1$, which implies $\psi(u)-2\psi(e_0)$ has a nonzero component in $P_0^\perp$. Since $\psi(v)\in P_1^\perp$, then $\Theta(u, \psi\otimes_{\C} v)=0$ if and only if $\psi(v)=0$, hence if and only if $v=0$. \medskip
 
 (ii) The $J$-holomorphicity of $\mathcal{P}_0$ from Definition \ref{Defn:ReduciblePencilOdd}(ii) implies that each $\psi \in \mathcal{P}_0$ gives a $\C$-linear identification of $P_1$ and $P_1^\bot$. The remaining claims of (ii) follow from (i).\medskip 

  (iii) This is a direct consequence of (ii) and of Lemma \ref{Lem:EinBaseTopology}. Indeed, since $\mathcal{R}^\bot \rightarrow Q_+(P)$ is a trivial vector bundle, the total space of $\sphere(\mathcal{R}^\bot) \rightarrow Q_+(P)$ is diffeomorphic to $\sphere(P)\times \sphere(P_1^\perp)$. Hence, under the quotient $(u,v)\sim (-u,-v)$, we see $\mathcal{B}(\mathcal{P}_0) \cong \sphere(\mathcal{R}^\bot)/\sim$ is diffeomorphic to $\Ein(P\oplus P_1^\perp)\cong \Ein^{p-1,p-2}$. 
 \end{proof}

\subsection{Description of the fiber, \texorpdfstring{$p$}{p} even}
\label{Sec:Fiber_Even}

We proceed in a very similar manner as in the odd case. We first define a notion a \emph{reducible pencil} $\mathcal{P}_0\subset \T\X$, and then we parametrize the base of this pencil. This time, the pencils $\mathcal{P}_0$ defined are, up to the $\SO_0(p,p+1)$-action, exactly the tangent pencil to a copy of $\Ha^2_{\iota_0}$ for $\iota_0$ the $\SL(2,\R)$-subgroup associated to the integer partition $2p+1=\underbrace{2+2+\dots+2}_p+1.$ 
\medskip

The background data is as follows. 
\begin{itemize}
    \item Let $E_0\oplus E_1=\R^{p,p+1}$ be a orthogonal decomposition with $E_0\simeq \R^{0,1}$ and $E_1\simeq \R^{p,p}$.
    \item Let $P\in \X$ be of the form $P \in \Gr_{(p,0)}(E_1)$. Thus, $P^\bot$ decomposes as $P^\perp= E_0\oplus P_1^\perp$, where $P^\perp_1\subset E_1$.
    \item Fix complex structures $J_1:P_1^\perp\to P_1^\perp$ and $J'_1:P\to P$ that are $q$-orthogonal.
\end{itemize}
This allows us to identify both $P$ and $P_1^\perp$ with $\C^k$ for $k= \frac{p}{2}$. We will write $J$ for the complex structure obtained on $E_1= P \oplus P_1^\bot$.
\medskip

\begin{definition}[Reducible Pencil, Even Case]\label{Defn:ReduciblePencilEven}
Let $p =2k$ be even. 
We call a pencil $\mathcal{P}_0\subset \T_P\X\cong \Hom_{\R}(P,P^\perp)$ \textbf{reducible with data} $\mathbf{(E_0,J)}$ if it is $\Ein$-regular and is a complex line in $\Hom_{\C}(P,P_1^\perp)\subset \Hom_{\R}(P,P^\perp)$. 
\end{definition}

In other words, the above definition asserts $\mathcal{P}_0 $ admits the description $\mathcal{P}_0 =\spann_{\C} \langle \psi \rangle $, where $\psi:P\rightarrow P_1^\bot$ is complex-linear and of full rank. 

\begin{remark}
The set of reducible pencils $\mathcal{P}_0\subset \T\X$ in the tangent bundle $\T\X$ comprises exactly one $\SO_0(p,p+1)$-orbit, namely the orbit of any tangent pencil to $\Ha^{2}_{\iota_0}$. 
\end{remark}

We can now describe geometrically the base of a reducible pencil. 
For the proof, recall that if $V$ is a real vector space with a Euclidean metric, $\T_vQ_+(V)$ is naturally identified with $\langle v\rangle ^\perp\subset V$. If moreover $V$ admits a complex structure $J$ preserving the metric, $\R\{v \}^\perp$ naturally decomposes as $ \R \{Jv \} \oplus\C \{v\} ^\perp $. 

\begin{lemma}[Geometry of Simplified Pencil, Even Case]\label{Lem:FiberEvenQuasiNatural}
Let $p \geq 4$ be even and $\mathcal{P}_0 \subset \T_P\X$ a reducible pencil with data $(E_0, J)$. 
\begin{enumerate}[label=(\roman*)]
    \item Fix a non-zero element $\psi\in \mathcal{P}_0$ and a norm one vector $e_0\in Q_+(E_0)$. There is a natural vector bundle isomorphism 
\[ \widehat{\Theta}_{\psi}: \T Q_+(P) \rightarrow \mathcal{R}^\bot, \; \; \; \;\]
where $\widehat{\Theta}_{\psi}(u, Ju) = (u, e_0)$ and $\widehat{\Theta}_{\psi}(u, v) = \left(u, \psi(v)\right)$ for $v \in \C\{u\}^\bot$. 
\item In particular, $\widehat{\Theta}_{\psi}$ induces a a natural diffeomorphism $\Theta_{\psi}:\T^1\mathbb{P}(P) \rightarrow \mathcal{B}(\mathcal{P}_0)$ depending only on $(E_0,\mathcal{P}_0, J, e_0)$. Hence, $\mathcal{B}(\mathcal{P}_0)\cong \T^1\RP^{p-1}$.
\end{enumerate}
\end{lemma}

\begin{proof}
(i) Let $u \in Q_+(P)$ be a fixed element and recall that $\mathcal{R}_u=\lbrace \psi(u)\mid \psi\in \mathcal{P}_0\rbrace$. Observe that $\mathcal{R}_u^\perp$ decomposes as $E_0\oplus \C\{ \psi(u) \}^\perp$. Since $\psi$ defines an isometry between $P$ and $P_1^\perp$, we see 
$\psi$ maps $\C \{u\}^\bot$ to $\C\{ \psi(u)\}^\bot$. 
Hence, there is a unique vector bundle morphism $\widehat{\Theta}_{\psi}: \T Q_+(P) \rightarrow \mathcal{R}^\bot$ satisfying $\widehat{\Theta}_{\psi}(u, Ju) = (u, e_0)$ and 
$\widehat{\Theta}_{\psi}(u, v) = (u, \psi (v))$. The map $\widehat{\Theta}_{\psi}$ is injective and hence an isomorphism by dimension count. 

(ii) This is a consequence of (i) and Lemma \ref{Lem:EinBaseTopology}, and the fact that the quotient $\T^1Q_+(P)/\sim$, where $(u,v)\sim (-u,-v)$, is naturally identified with $\T^1\mathbb{P}(P)$.
\medskip
\end{proof}

Recall $p=2k$. Given an orthogonal complex structure $J$ on $\R^{2k}$, we define an action $\eta_0$ of $\mathrm{U}(1)$ on $\T^1\RP^{p-1}$.
To this end, recall for $\ell \in \RP^{p-1}$ that there is a natural identification $\T_{\ell}\RP^{p-1} \cong \Hom_{\R}(\ell, \ell^\bot)$. 
Consequently, we can decompose the tangent space fiberwise by
\[ \T_{\ell}\RP^{p-1} \cong \Hom_{\R}(\ell, \ell^\bot) = \Hom_{\R}(\ell,J\ell) \oplus \Hom_{\R}(\ell, \C \{ \ell\}^\perp ).\]
Let $\mathscr{L}\rightarrow \RP^{p-1}$ be the tautological line bundle. The above pointwise splitting yields a bundle splitting $\T\RP^{p-1} = J\mathscr{L} \oplus (\mathscr{L}^\C)^\bot$. Since $\T\RP^{n} =\T\sphere^n/(-\id,-\id)$, we may denote $[u,v] = [ \pm (u,v)] \in \T_{[u]}\RP^n$ for $(u,v) \in \T_u\sphere^n$. 

\begin{definition}[The $\eta_0$-action]\label{Defn:eta0Action}
Decompose $[u,x] \in \T_{[u]}\RP^{p-1}$ as $[u,x]=[u,v+w]$ with $v \in J\mathscr{L}|_{u}$ and $w \in (\mathscr{L}^\C)^\bot|_{u}$. We define the action $\eta_0:\U(1)\rightarrow \Diff(\T^1\RP^{p-1})$ by 
\begin{align}\label{eta0}
    \eta_0(\lambda) \cdot [u,v+w] =[u,v+\lambda w] .
\end{align} 
\end{definition}

That is, $\eta_0$ acts by fiber bundle automorphisms lifting the identity map on the base, fixing the line bundle $J\mathscr{L}$, and acting by rotation on the complex vector bundle $(\mathscr{L}^\C)^\bot$. 

\begin{remark}
Given $P \in \X$ and $\mathcal{P}_0\subset \T_P\X$ a reducible pencil with data $(E_0,J)$, 
we can perform the same construction as in Definition \ref{Defn:eta0Action}, replacing $\T^1\RP^{p-1}$ with $\T^1\mathbb{P}(P)$.
\end{remark}

Now, note that $\U(1)$ also acts naturally on the right on the pencil $\mathcal{P}_0$ by rotation, i.e., so that $\psi\cdot \lambda=\psi\circ \lambda$.

Now, let $\mathcal{P}_0 \subset \T_P\X$ be a reducible pencil with data $(E_0,J)$. 
The next result provides need a more precise understanding of how the identifications in Lemma \ref{Lem:FiberEvenQuasiNatural} depend on the choice of $\psi$. Here, we use the left and right action of $\U(1)$ on $\T^1\mathbb{P}(P)$ and $\mathcal{P}_0$, respectively.

\begin{theorem}[Canonical Description of Fiber, $p$ Even]\label{thm:FiberEvenNatural}
Let $p=2k$ be an even integer and $\mathcal{P}_0 \subset \T_P\X$ be a reducible pencil with data $(E_0,J)$. 
Fix $e_0\in E_0$ of norm one and $\psi \in \sphere(\mathcal{P}_0)$.
\begin{enumerate}[label=(\roman*)]
    \item The diffeomorphism 
$\Theta_{\psi}:\T^1\mathbb{P}(P)\rightarrow \mathcal{B}(\mathcal{P}_0)$ from Lemma \ref{Lem:FiberEvenQuasiNatural}(ii) satisfies 
\[ \Theta_{\psi\cdot \lambda}= \Theta_{\psi}\circ \eta_0(\lambda).\]
    \item The base of pencil $\mathcal{B}(\mathcal{P}_0)$ is canonically identified with $\big(\mathbb{S}(\mathcal{P}_0)\times \T^1\mathbb{P}(P) \big)/\sim$, under the quotient $(\psi \cdot \lambda ,x)\simeq (\psi,\eta_0(\lambda)\cdot x)$ for all $\lambda \in \mathrm{U}(1)$.
\end{enumerate}

\end{theorem}

\begin{proof}
    (i) For any non-zero element $\psi\in \mathcal{P}_0$, Lemma \ref{Lem:FiberEvenQuasiNatural} provides the diffeomorphism $\Theta_\psi:\T^1\mathbb{P}(P)\rightarrow \mathcal{B}(\mathcal{P}_0)$, which is induced by a vector bundle isomorphism $\widehat\Theta_\psi:\T Q_+(P)\to \mathcal{R}^\perp$. The equivariance property for $\Theta_{\psi}$ follows if we prove the property for $\widehat{\Theta}_{\psi}$. We now undertake this proof.

     By definition of $\widehat{\Theta}_{\psi}$, we first note that 
    \[ \widehat{\Theta}_{\psi \cdot \lambda }(u, Ju)=\widehat{\Theta}_{\psi}(u, Ju) = (u, e_0) = \widehat{\Theta}_{\psi}\big( \eta(\lambda)\cdot (u,Ju) \big).\]
     This verifies the desired equivariance when $v \in J\mathscr{L}$. 

     Otherwise, we may suppose $v \in (\mathscr{L}^\C)^\bot|_u$. In this case, unraveling definition shows 
     \[ \widehat{\Theta}_{\psi \cdot \lambda}(u, v)= (u, (\psi \circ \lambda)(v))= (u, \psi(\lambda v)) = \widehat{\Theta}_{\psi}(u,\lambda v)= \widehat{\Theta}_{\psi}\circ \eta_0(\lambda)  (u,v).\]
     The equivariance property of $\widehat{\Theta}$ follows.

     (ii) This is an immediate consequence of (i). 
\end{proof}

\subsection{Deformations of Tangent Pencils}\label{Sec:HiggsPencilDeformation}

Let $\iota:\SL(2,\R)\rightarrow \SO_0(p,p+1)$ be an $\Ein$-regular representation and $\Ha^2_{\iota} \subset \X$ the associated sub-symmetric space. In this section, we produce a smooth deformation $(\mathscr{P}_t)_{t \in [0,1]}$ from the tangent pencil
$\mathscr{P}_1 \in \Omega^0(\Ha^2_{\iota}, \Gr_2(\T \Ha^2_{\iota}))$ to a `simplified' pencil $\mathscr{P}_0 \in \Omega^0(\Ha^2_{\iota}, \Gr_2(\T\Ha^{2}_{\iota}))$ that is pointwise a \emph{reducible pencil} in the sense of Definitions \ref{Defn:ReduciblePencilOdd}, \ref{Defn:ReduciblePencilEven}, depending on parity of $p$. We do this using Higgs bundles, with the ideas introduced in Remark \ref{Remk:HiggsPencils}. As a consequence of the deformation and the geometry in $\S$\ref{Sec:Fiber_Odd}, \ref{Sec:Fiber_Even}, we find the fiber of $M_{\rho}\rightarrow S$ for $\rho$ an $\iota$-Fuchsian deformation.  

We now describe the strategy to build the deformation $(\mathscr{P}_t)$. 
By Proposition \ref{Prop:EinRegularSL2}, the integer partition associated to $\iota$ obtains the form 
$(2p+1)= (2q+1)+2r_1+2r_1+\cdots +2r_k+2r_k$. The desired deformation occurs in blocks, using the corresponding deformation of $\R^{p,p+1}$ into $\iota$-invariant sub-representations. Consequently, we then consider deformations of pencils for $\SO_0(q,q+1)$-Hitchin Higgs bundles, and for the sum of $2r_i+2r_i$ reducibles in $\SO_0(2r_i,2r_i)$. Taken in conjunction, we obtain the desired deformation $(\mathscr{P}_t)_{t \in[0,1]}$ block-by-block. 

In this section, we use all the Higgs bundle notation from $\S$\ref{Subsec:FuchsianHitchin} and $\S$\ref{Subsec:Uniformize_p+p}. 

\subsubsection{Fuchsian-Hitchin in $\SO(p,p+1)$ for $p$ odd.}\label{Sec:GlobalPencilHitchinOdd}

Let $p = (2k+1) \geq 3$ be an odd positive integer and $\rho:\pi_1S\rightarrow \SO_0(p,p+1)$ a Fuchsian-Hitchin representation. By \cite{Hit92}, for an appropriate Riemann surface $\Sigma = (S,J)$, the representation $\rho$ is associated to the Higgs bundle described in Section \ref{Subsec:FuchsianHitchin}:

\begin{equation}\label{eq:FuchsianHitchinBundleOdd} \begin{cases}  \V &= \bigoplus_{i=-p}^{p} \K^{i} ,  \\
\Phi   &= \left( \K^{p} \stackrel{1}{\longrightarrow} \mathcal{K}^{p-1} \stackrel{1}{\longrightarrow}\cdots \stackrel{1}{\longrightarrow} \mathcal{O} \stackrel{1}{\longrightarrow}  \cdots \stackrel{1}{\longrightarrow} \K^{1-p}  \stackrel{1}{\longrightarrow} \K^{-p} \right ).
\end{cases} \end{equation}
Note that for $p$ odd, the bundle $\mathcal{U}$ from \eqref{U_Subbundle} is the sum of even powers of $\mathcal{K}$ and $\mathcal{V}$ from \eqref{V_Subbundle} is the sum of odd powers of $\K$. In particular, $\mathcal{O}\in \mathcal{U}$.

\medskip

We can make sense of a \emph{reducible pencil} in the sense of Definition \ref{Defn:ReduciblePencilOdd} pointwise for $x \in S$, replacing $\T_P\X$ with the vector space $\Hom_{\R}(\mathcal{U}^\R, \mathcal{V}^{\R})|_x$ (here, recall Remark \ref{Remk:HiggsPencils}). To this end, we  introduce bundle versions of the auxiliary data associated in Definition \ref{Defn:ReduciblePencilOdd}: 

\begin{itemize}
    \item We define a bundle decomposition $\mathcal{E}=\mathcal{E}_0\oplus \mathcal{E}_1$ such that $\mathcal{E}_0^\R$ and $\mathcal{E}_1^{\R}$ are fiberwise isomorphic to $\R^{1,2}$ and $\R^{2k,2k}$.
    \item We decompose further as $\mathcal{U}_i = \mathcal{U} \cap \mathcal{E}_i$ and $\mathcal{V}_i = \mathcal{U} \cap \mathcal{E}_i$. 
    \item We place a complex structure $\mathcal{J}$ on $\mathcal{U}_0\oplus \mathcal{V}_0 \oplus \mathcal{V}_1$ that preserves $\mathcal{U}_0^\R, \mathcal{V}_0^{\R}, \mathcal{V}_{1}^{\R}$ (which mirror $P_0, P_0^\bot, P_1^\bot$, respectively) and is $Q$-orthogonal on these real loci. 
\end{itemize}

We first define the bundle decomposition as follows:
\begin{align}
    \mathcal{E}_0 &=\mathcal{K}\oplus \mathcal{O}\oplus \mathcal{K}^{-1},\\ 
\mathcal{E}_1&=\bigoplus_{i=2}^{p} (\K^{i}\oplus\K^{-i}) .
\end{align}

We define $\mathcal{J}: \mathcal{V}_0 \oplus \mathcal{U}_1 \oplus \mathcal{V}_1 \rightarrow \mathcal{V}_0 \oplus \mathcal{U}_1 \oplus \mathcal{V}_1$ as follows:
\begin{align}
   \mathcal{J}(v)= \begin{cases}
         -\mathbf{i}v & v \in \K^j, \;\;\;j > 0\\
         \mathbf{i}v & v \in \K^{-j}, \,j > 0 \\
    \end{cases}
\end{align}
By the equation \eqref{eq:RealLocusHitchin}, the complex structure $\mathcal{J}$ satisfies the desired properties. 

Finally, we define holomorphic endomorphism-valued one forms $\varphi_t$, for $t\in [0,1]$, by:
 \begin{align}
        \varphi_{t} &= \left( \bigoplus_{0< j\leq k} \K^{2j} \stackrel{t\cdot 1}{\longrightarrow} \mathcal{K}^{2j-1}\right )\oplus \left(\bigoplus_{-k\leq j\leq 0} \K^{2j} \stackrel{1}{\longrightarrow} \mathcal{K}^{2j-1}\right )\in H^0(\Hom(\mathcal{U},\mathcal{V})\otimes \K), \label{DeformationOdd1}\\
    -\varphi_{t}^{*Q} &=\left(\bigoplus_{0\leq j\leq k} \K^{2j+1} \stackrel{1}{\longrightarrow} \mathcal{K}^{2j}\right )\oplus \left( \bigoplus_{-k\leq j<0} \K^{2j+1} \stackrel{t\cdot 1}{\longrightarrow} \mathcal{K}^{2j}\right ) \in H^0(\Hom(\mathcal{V},\mathcal{U})\otimes \K). \label{DeformationOdd2}
    \end{align}

Let us write 
\[ \Phi_t=(\varphi_{t}-\varphi_{t}^{*Q}) \in H^0(\Sigma, \End(\V)\otimes \K). \]
Associated to this object $\Phi_t$ is a Higgs bundle pencil $\mathfrak{E}_t\in \Omega^0(\Sigma,\Gr_2(\Hom(\mathcal{U}^\R, \mathcal{V}^\R)))$. In short, $\Phi_t+\Phi_t^{*h}$ consists of four sub-tensors, and we take the parts of appropriate type. Precisely, $\mathfrak{E}_t$ is given fiberwise by 
\begin{align}\label{HiggsPencilDefn}
 \mathfrak{E}_t\mid _x=\lbrace(\varphi_t-(\varphi_t^{*Q})^{*h})(v)\mid v\in \T_x\Sigma \rbrace \in \Gr_{2}(\Hom(\mathcal{U}^\R,\mathcal{V}^\R)_{\mid x}). 
\end{align}
Observe that by construction, $(\Phi_t+\Phi_t^{*h})(v)$ preserves the real locus $\V^\R=\mathcal{U}^\R \oplus \mathcal{V}^\R$ and moreover $[\varphi_t-(\varphi_t^{*Q})^{*h}](v)$ is the restriction of the former transformation to a map $\mathcal{U}^{\R}\rightarrow \mathcal{V}^{\R}$. 

\begin{figure}[ht]
\centering
\includegraphics[width=.75\textwidth]{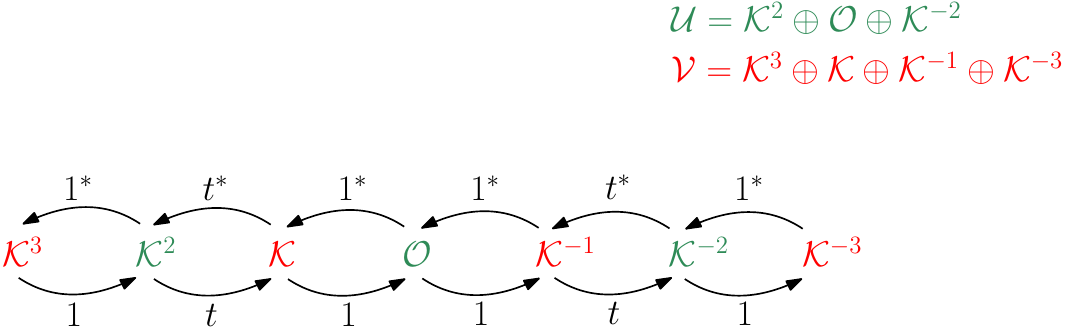}
\caption{An illustration of the deformation of pencils $\mathfrak{E}_t$. The forwards arrows come from $\Phi_t$ and the backwards arrows from $\Phi_t^{*h}$. The Higgs bundle pencil $\mathfrak{E}_t =\big[\varphi-(\varphi^{*Q})^{*h}\big]$ consists of all arrows from $\mathcal{U}$ to $\mathcal{V}$.}
\label{Fig:DeformationOdd}
\end{figure}

Note that $\Phi=\Phi_1$ is the original Higgs field. The Higgs bundle pencils $(\mathfrak{E}_t)_{t\in [0,1]}$ serve the purpose of deforming the tangent pencil to $\Ha^{2}_{\Delta}$, 
the sub-symmetric space of the principal $\PSL(2,\R)$-subgroup, to a simpler pencil. The following remark clarifies this point. 

\begin{remark}
Note that here we are not changing the ambient Higgs bundle $(\V,\Phi)$, the corresponding harmonic metric or the flat connection; we are only defining these new sections $\Phi_t \in H^0(\End(\V)\otimes\K)$ of the same type as the Higgs field.

Fix $\tilde{x} \in \widetilde{\Sigma}, \, x=\pi(\tilde{x}) \in \Sigma$ and $P = f(x) \in \X$. 
By Remark \ref{Remk:HiggsPencils}, the tangent pencil $\mathcal{P} =\T_P\Ha^2_{\Delta}\subset \T_{P}\X$ is identified with ${\mathfrak{E}_1}_{\mid x}$. Moreover, the object $\mathfrak{E}_t|_x$ is identified with a pencil $\mathcal{P}_t \subset \T_P\X$, thereby defining a deformation from $\mathcal{P}=\mathcal{P}_1$ to $\mathcal{P}_0$ in $\T_P\X$. 
\end{remark}

Using the identifications made in the above remark, we show $\mathfrak{E}_t|_x$ is pointwise an $\Ein$-regular a pencil for all $t\in [0,1]$. 
\begin{lemma}[Regularity of Pencils, Hitchin Odd Case]\label{Lem:RegularPencilOddCase}
For any point $x \in \Sigma$ and any $t \in [0,1]$, the pencil $ {\mathfrak{E}_t}|_x$ is $\Ein$-regular. 
\end{lemma} 

\begin{proof}
By Proposition \ref{Prop:EinRegularity}, we need to show that for all non-zero $v\in \T_{x}S$ the endomorphism $(\Phi_t+\Phi_t^{*h})(v)$ has the largest possible rank, namely $2p$. Recall that the metric $h$ on $(\V, \Phi)$ is diagonal. For $t\neq 0$ the restriction of $(\Phi_t+\Phi_t^{*h})(v)$ from $\bigoplus_{j=-p+1}^p\mathcal{K}^j$ to $\bigoplus_{j=-p}^{p-1}\mathcal{K}^j$ is upper triangular with non-zero diagonal coefficients, and hence has rank $2p$. For $t=0$, $(\Phi_t+\Phi_t^{*h})(v)$ is block diagonal with $p$ invertible blocks of size $2$, and hence has rank $2p$.
\end{proof}

We now explain why $\mathfrak{E}_0$ is pointwise a reducible pencil in the sense of Definition \ref{Defn:ReduciblePencilOdd}.

Let $x \in S$ be fixed and set $\mathcal{P}_0 =\mathfrak{E}_0|_x$. 
Define $E_0 = \mathcal{E}^\R_0|_x$ and $E_1 = \mathcal{E}_1^{\R}|_x$. 
We write $P\coloneq\mathcal{U}^\R|_x$. 
As in Subsection \ref{Sec:Fiber_Odd}, we obtain a decomposition $P=P_0\oplus P_1$ and $P^\perp=P_0^\perp\oplus P_1^\perp$. That is, we define $P_i \coloneq P \cap \mathcal{E}_i$ and $P^\bot_i \coloneq P^\bot \cap \mathcal{E}_i$. 
Let $J$ be the restriction of $\mathcal{J}_x$ to $P_0^\perp \oplus E_1$. 

\begin{lemma}
\label{lem:ReduciblePencilOdd}
    For all $x\in \Sigma$ the pencil $\mathcal{P}_0={\mathfrak{E}_0}{\mid_x}$ is a reducible pencil with data $(E_0,J)$.
\end{lemma}

\begin{proof}

All that remains to prove are the two properties (i) and (ii) of Definition \ref{Defn:ReduciblePencilOdd}.

Recall that $\varphi_0: \T^{1,0}\Sigma \rightarrow \End_{\C}(\V)$ and $(\varphi_0^{*Q})^{*h}:\T^{0,1}\Sigma \rightarrow \End_{\C}(\V)$ are holomorphic bundle maps.
Thus, for $v\in \T_x\Sigma$ and $u\in P_0$ or $u \in P_1$, we find 
 \begin{align*}
    \varphi_0(\mathbf{i}v)(u)&= \varphi_0(\mathbf{i}v^{1,0})(u)=\mathbf{i}\varphi_0(v)(u)=J\big(\varphi_0(v)(u)\big), \\ 
    (\varphi_0^{*Q})^{*h}(\mathbf{i}v)(u)&=(\varphi_0^{*Q})^{*h}(-\mathbf{i}v^{0,1})(u)= -\mathbf{i}(\varphi_0^{*Q})^{*h}(u)=J\big((\varphi_0^{*Q})^{*h}(v)(u)\big).
\end{align*} 
Here, we use that the images of $\varphi_0$ and $(\varphi_0^{*Q})^{*h}$, respectively, are in the $\mathbf{i}$ and $-\mathbf{i}$ eigenspaces of $J$. See Figure \ref{Fig:DeformationOdd}. 
Point (i) of Definition \ref{Defn:ReduciblePencilOdd} follows. 
Then note that $J$ commutes with $\varphi_0$ and $(\varphi_0^{*Q})^{*h}$, which implies (ii).
\end{proof}

\subsubsection{Fuchsian-Hitchin in $\SO(p,p+1)$ for $p$ even.}\label{Sec:GlobalPencilHitchinEven}
We proceed in a very similar manner. Let $p = 2k \geq 4$ be an even positive integer and $\rho:\pi_1S\rightarrow \SO_0(p,p+1)$ a Fuchsian-Hitchin representation. For an appropriate Riemann surface $\Sigma = (S,J)$, the Higgs bundle associated to $\rho$ is as follows:

\begin{equation}\label{eq:FuchsianHitchinBundleEven} \begin{cases}  \V &= \bigoplus_{i=-p}^{p} \K^{i} ,  \\
\Phi   &= \left( \K^{p} \stackrel{1}{\longrightarrow} \mathcal{K}^{p-1} \stackrel{1}{\longrightarrow}\cdots \stackrel{1}{\longrightarrow} \mathcal{O} \stackrel{1}{\longrightarrow}  \cdots \stackrel{1}{\longrightarrow} \K^{1-p}  \stackrel{1}{\longrightarrow} \K^{-p} \right ).
\end{cases} \end{equation}
Note that for $p$ even, the bundle $\mathcal{U}$ from \eqref{U_Subbundle} is the sum of the odd powers of $\mathcal{K}$ and $\mathcal{V}$ from \eqref{V_Subbundle} is the sum of the even powers of $\K$. In particular, $\mathcal{O}\in \mathcal{V}$ in this case. \medskip

Again, we introduce bundle analogues of the auxiliary data of a reducible pencil. Namely,
\begin{itemize}
    \item We fix a bundle splitting $\V = \mathcal{E}_0 \oplus \mathcal{E}_1$ with $\mathcal{E}_0^\R$ and $\mathcal{E}_1^\R$ fiberwise isomorphic to $\R^{0,1}$ and $\R^{2k,2k}$, respectively. 
    \item We place a complex structure $\mathcal{J}$ on $\mathcal{E}_1$ such that 
    $\mathcal{U}_1^\R= \mathcal{U}^{\R}\cap \mathcal{E}_1$ and $\mathcal{V}_1^\R= \mathcal{V}^{\R}\cap \mathcal{E}_1$ are preserved and $\mathcal{J}$ is $Q$-orthogonal on these real loci. 
\end{itemize}
This time, the relevant bundle decomposition $\mathcal{E}=\mathcal{E}_0\oplus \mathcal{E}_1$ given by 
\begin{align*}
    \mathcal{E}_0 &= \mathcal{O}, \\
    \mathcal{E}_1 &=\bigoplus_{i=1}^{p} (\K^{i}\oplus\K^{-i}) .
\end{align*}
We define the complex structure $\mathcal{J}$ on $\mathcal{E}_1$ as before:  
\begin{align}
    \mathcal{J}(v)=\begin{cases}
                \mathcal{J}_x(v)= -\mathbf{i}v & v\in \mathcal{K}^j, \;\;\;\;j > 0 \\
                \mathcal{J}_x(v)= \mathbf{i}v & v\in \mathcal{K}^{-j}, \;\;j > 0 .\end{cases} 
\end{align}

We then build new endomorphism-valued one forms $\varphi_t$, for $t\in [0,1]$, by:
\begin{align}
    \varphi_{t} &= \left(\bigoplus_{-k\leq j<0} \K^{2j+1} \stackrel{1}{\longrightarrow} \mathcal{K}^{2j}\right )\oplus \left(\bigoplus_{0\leq j\leq k} \K^{2j+1} \stackrel{t\cdot 1}{\longrightarrow} \mathcal{K}^{2j}\right)\in H^0(\Hom(\mathcal{U},\mathcal{V})\otimes \K), \label{DeformationEven1} \\
    -\varphi_{t}^{*Q} &= \left(\bigoplus_{0< j< k} \K^{2j} \stackrel{1}{\longrightarrow} \mathcal{K}^{2j-1}\right )\oplus \left( \bigoplus_{-k\leq j\leq 0} \K^{2j} \stackrel{t\cdot 1}{\longrightarrow} \mathcal{K}^{2j-1}\right )\in H^0(\Hom(\mathcal{V},\mathcal{U}) \otimes\K). \label{DeformationEven2}
\end{align}

Once again, we set $\Phi_t \coloneqq\varphi_{t}-\varphi_{t}^{*Q}$ and define the pencil $\mathfrak{E}_t \in \Omega^0(\Sigma, \Gr_2(\Hom(\mathcal{U}^\R, \mathcal{V}^\R)))$ just as in \eqref{HiggsPencilDefn}. 
Note that $\varphi_1=\varphi$ from Section \ref{Subsec:FuchsianHitchin} and $\Phi=\Phi_1$ is the original Higgs field. 

\begin{figure}[ht]
\centering
\includegraphics[width=.80\textwidth]{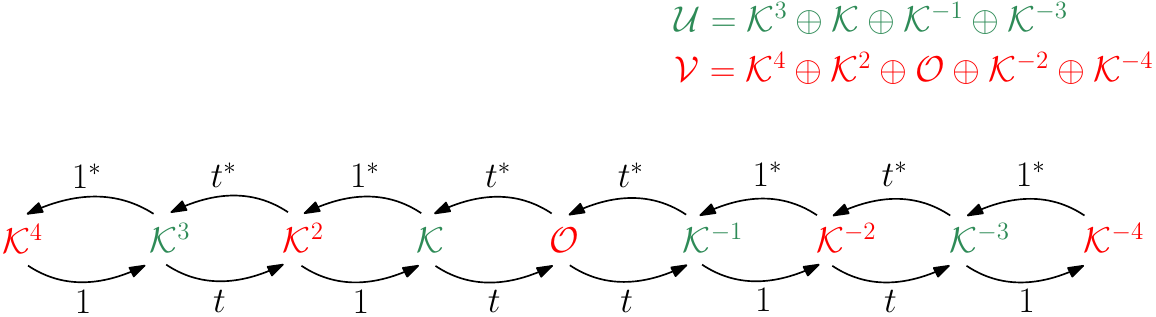}
\caption{An illustration of the deformation of pencils $\mathfrak{E}_t$. The forwards arrows come from $\Phi_t$ and the backwards arrows from $\Phi_t^{*h}$. The Higgs bundle pencil $\mathfrak{E}_t =\big[\varphi-(\varphi^{*Q})^{*h}\big]$ consists of all arrows from $\mathcal{U}$ to $\mathcal{V}$.}
\label{Fig:DeformationEven}
\end{figure}

The Higgs bundle pencil $\mathfrak{E}_t$ again maintains its $\Ein$-regular nature through deformation. 

\begin{lemma}[Regularity of Pencils, Hitchin Even Case]\label{Lem:RegularPencilEvenCase}
Fix $x \in \Sigma$. For all $t \in [0,1]$, the pencil $\mathfrak{E}_t|_x$ is $\Ein$-regular. 
\end{lemma}

\begin{proof}
 Just as in Lemma \ref{Lem:RegularPencilOddCase}, we need only show that each endomorphism $\psi:=(\Phi_{t}+\Phi_{t}^{*h})(v)$ has rank $2p$, for $v \in \T_{x}S$. By the equation \eqref{CartanProj}, we know a priori $\psi$ has rank $\leq 2p$. 
 
 Now, recall the harmonic metric $h$ is diagonal in the splitting $\bigoplus_{j=p}^{-p}\K^j$. For $t\neq 0$ and any non-zero $v\in \T_{p_0}S$ the block of $\psi =(\Phi_t+\Phi_t^{*h})(v)$ from $\bigoplus_{j=-p+1}^p\mathcal{K}^j$ to $\bigoplus_{j=-p}^{p-1}\mathcal{K}^j$ is upper triangular with non-zero diagonal coefficients, and hence has rank $\geq 2p$. For $t=0$, $\psi$ is block diagonal with $p$ invertible $(2\times 2)$ blocks and a $(1\times 1)$ zero block and hence has rank $2p$.
\end{proof}

\medskip

We now check that the `simplified' Higgs bundle pencil $\mathfrak{E}_0$ is pointwise a reducible pencil in the sense of Definition \ref{Defn:ReduciblePencilOdd}. 

Fix $x \in \Sigma$ and write $P = \mathcal{U}^\R|_x$. 
Let $E_0$ and $E_1$ be respectively the real loci of ${\mathcal{E}_0}|_{x}$ and ${\mathcal{E}_1}|_{x}$. We obtain  decompositions $P=P_0\oplus P_1$ and $P^\perp=P_0^\perp\oplus P_1^\perp$ by defining $P_i =P \cap \V_i$ and $P^\bot_i = P^\bot \cap \V_i$. Note that $\mathcal{J}_x$ preserves the real locus $\V^\R$ inside $\mathcal{E}$ as well as $P_0$ and $E_1$. Hence, we let let $J$ be the restriction of $\mathcal{J}_x$ to $ E_1$.

\begin{lemma}
\label{lem:ReduciblePencilEven}
    For all $x\in \Sigma$ the pencil $\mathcal{P}_0={\mathfrak{E}_0}{\mid_x}$ is a reducible pencil with data $(E_0,J)$.
\end{lemma}

\begin{proof}

The proof is completely analogous to the odd case.
\end{proof}

\subsubsection{Reducible $2r+2r$ Fuchsian in $\SO_0(2r,2r)$.}

In this section, we consider a notion of reducible pencils for $\R^{2r,2r}$. We then construct a Higgs bundle version of such a pencil for the $\SO_0(2r,2r)$-Higgs bundle from Section \ref{Subsec:Uniformize_p+p} associated to a sum of dual $\SL(2r, \R)$-Fuchsian-Hitchin representations $\rho:\pi_1S\rightarrow\SO_0(2r,2r)$. 

Recall that such Higgs bundles obtain the form 
\begin{align}
    \V_{2r,2r} &= \bigoplus_{n=-2r+1}^{2r-1} (\K^{n/2} \oplus \K^{-n/2}) \label{SO(r,r)VectorBundle}\\ 
    Q &= \bigoplus_{n=-2r+1}^{2r-1} (-1)^{n-1}\begin{pmatrix} 0 & 1 \\ 1& 0 \end{pmatrix}
\end{align}

We do not define $\mathcal{E}_0$ and $\mathcal{E}_1$ in this case. As we will see in the next section, $\mathcal{E}_0= \underline{\R^0}$ is appropriate here. 

The relevant complex structure $\mathcal{J}$ in this case is as follows: 
\[ \mathcal{J}=\bigoplus_{n=-2r+1}^{2r-1}\begin{pmatrix} \mathbf{i} & 0 \\ 0 & -\mathbf{i} \end{pmatrix},\] 
in the above splitting \eqref{SO(r,r)VectorBundle}. 
\medskip

For $t\in[0,1]$, consider the following endomorphism-valued one forms $\Phi_t\in H^0(\End(\mathcal{E})\otimes \mathcal{K})$: 
\begin{align*}
   \K^{\frac{-1+2r}{2}} &\stackrel{1}{\longrightarrow} \K^{\frac{-3+2r}{2}} \stackrel{t}{\longrightarrow} \cdots \stackrel{t}{\longrightarrow} \K^{\frac{+3-2r}{2}}\stackrel{1}{\longrightarrow} \K^{\frac{+1-2r}{2}}\\  
   \K^{\frac{+1-2r}{2}} 
   &\stackrel{1}{\longleftarrow} \K^{\frac{+3-2r}{2}} \stackrel{t}{\longleftarrow} \cdots \stackrel{t}{\longleftarrow} \K^{\frac{-3+2r}{2}}\stackrel{1}{\longleftarrow} \K^{\frac{-1+2r}{2}}.
\end{align*}
In short, $\Phi_{t}$ alternates between $1$ and $t$ for the maps between consecutive two dimensional sub-bundles of $\V_{2r,2r}$. See Figure \ref{fig:SO(p,p)Deformation} for an example with $2r=4$. These sections can also be written  $\Phi_t=\varphi_t-\varphi_t^{*Q}$, where $\varphi_t \in H^0(\Hom(\mathcal{U},\mathcal{V})\otimes \K)$. 

Again, we set $\mathfrak{E}_t$ to be the Higgs bundle pencil defined by \eqref{HiggsPencilDefn}. At time $t=1$, we have $\Phi_r^1 =\Phi_r$ as the original Higgs field. By a reasoning completely similar to Lemma 
\ref{Lem:RegularPencilOddCase}, we obtain the following result.  

\begin{lemma}[Regularity of Deformation, $2r+2r$ case]\label{lem:Regularity_r+r}
Fix $x \in \Sigma$. The pencil $\mathfrak{E}_t|_x$ is $\Ein$-regular in the sense that for every non-zero $v\in \T\Sigma$, the endomorphism 
$(\Phi_t+\Phi_t^{*h})(v) \in \End(\mathcal{E}_{2r,2r}|_{x})$ has full rank $4r$. 
\end{lemma}

The complex structure $\mathcal{J}$ has the following properties, which are slightly stronger than in the Fuchsian-Hitchin cases.
\begin{lemma}[Holomorphic Pencil, $2r+2r$ case]\label{lem:Holomorphic r+r}
Let $v\in \T_x\Sigma$ and $u\in \mathcal{U}^\R|_x$. Then 
the pencil is dually holomorphic in the following sense: for any transformation $\theta \in \mathfrak{E}_t|_x$, 
\[ \theta(\mathbf{i}v)(u) =\theta(v)(Ju)=J\big( \theta\big(v)(u)\big). \] 
\end{lemma}
\begin{proof}
Both equalities hold for $\theta=\varphi_0|_x$ and $\theta=(\varphi_0^{*Q})^{*h}|_x$, similar to Lemma \ref{Lem:RegularPencilOddCase}. 
\end{proof}

\begin{figure}[ht]
    \centering
    \includegraphics[width=.45\textwidth]{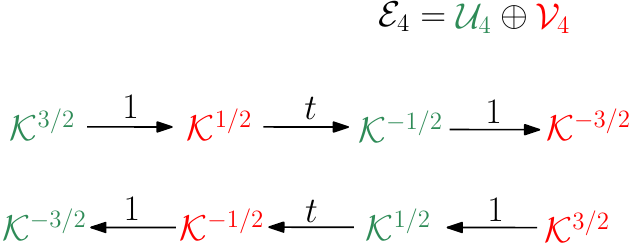}
    \caption{An illustration of the deformation of pencils $\mathfrak{E}_t$ for $2r=4$. The object $\Phi_t$ consists of all arrows in the diagram. The object $\varphi_t$ consists of all arrows from $\mathcal{U}_4$ to $\mathcal{V}_4$. The arrows of $\Phi_t^{*h}$ are not displayed.}
\label{Fig:Deformation_r+r}
    \label{fig:SO(p,p)Deformation}
\end{figure}

\subsubsection{Sum of an irreducible and $2r_i+2r_i$ Fuchsians in $\SO(p,p+1)$.}
\label{subsec:Sum of irreducibles}
We now finally arrive at the case of interest. 
Let $I$ be a partition of the form $2p+1=(2q+1)+2r_1+2r_1+\cdots +2r_k+2r_k$ 
and $\iota_I:\SL(2,\R)\rightarrow \SO_0(p,p+1)$ the  representation associated by Proposition \ref{Prop:EinRegularSL2}. 
Let $\rho: \pi_1S\rightarrow \SO_0(p,p+1)$ be an $\iota_I$-Fuchsian. We will deform the tangent pencil to $\Ha^2_{\iota_I}$ to a reducible pencil and then use this to deduce the topology of the fibers of $M_{\rho}\rightarrow S$. 

Let us consider the uniformizing Higgs bundle associated to $\rho$, denoted $(\mathcal{E},\Phi)$, with its subbundles $\mathcal{U}, \mathcal{V}$ and real locus $\V^\R$. We now introduce the analogue of a \emph{reducible pencil} for the Higgs bundle $(\V,\Phi)$. In essence, such a pencil is just a combination of all the reducible pencils for each Higgs sub-bundle according to the splitting (using notation from $\S$\ref{Sec:Uniformizing})
\[ (\V,\Phi) = (\V_{q}, \Phi_q)\oplus \bigoplus_{i=1}^{n} (\mathcal{E}_{2r_i, 2r_i}, \Phi_{2r_i, 2r_i}).\]

We now introduce further splittings and data:
\begin{itemize}
    \item Write  $\mathcal{E}=\mathcal{E}^\text{odd}\oplus \mathcal{E}^\text{even}$ where $\mathcal{E}^\text{odd}= \V_q$ and $\mathcal{E}^\text{even} = \bigoplus_{i=1}^n \V_{2r_i,2r_i}$. 
    \item Construct the bundle splitting $ \mathcal{E}=\mathcal{E}_0\oplus \mathcal{E}_1$, where $\mathcal{E}_0=\mathcal{E}_0^\text{odd}$ and $\mathcal{E}_1=\mathcal{E}^\text{even}\oplus \mathcal{E}_1^\text{odd}$.
    \item Define a complex structure $\mathcal{J}$ on $(\mathcal{E}_0 \cap \mathcal{V}) \oplus \mathcal{E}_1$ of the form $\mathcal{J}=\mathcal{J}^\text{odd}\oplus \mathcal{J}^\text{even}$, which is a further direct sum of the complex structures previously defined.   
    \item For all $t\in[0,1]$, we define $\phi_t,\Phi_t$ as previously described for $\mathcal{E}^\text{even}$ and $\mathcal{E}^\text{odd}$. 
\end{itemize}
Finally, we once more define the Higgs bundle base of pencil $\mathfrak{E}_t$, given pointwise by: 
\[ \mathfrak{E}_t|_x=\lbrace (\varphi_t-(\varphi_t^{*Q})^{*h})(v)\in \Hom(\mathcal{U}^\R,\mathcal{V}^\R)\mid v\in \T_x\Sigma\rbrace. \]
In this expression, $\Phi_t$ is the direct sum of the holomorphic endomorphism-valued one-forms defined in previous subsections in each block.

This deformation $\mathfrak{E}_t$ is $\Ein$-regular. 

\begin{lemma}[Regularity of Deformation, General Case]\label{Lem:PencilDeformationGeneralCase}
Fix $x \in \Sigma$. The pencil $\mathfrak{E}_t|_x$ is $\Ein$-regular for all $t$. 
\end{lemma}

\begin{proof}
Fix $v \in \T_x\Sigma$. We need to show $(\Phi_t+\Phi_t^{*h})(v)$ has the maximum possible rank, namely $2p$. 
Each endomorphism $(\Phi_t+\Phi_t^{*h})(v)$ is the sum of a map of corank $1$, on the odd block, and a sum of maps of full rank, on each even block, so of corank $1$ in total. Indeed, in the odd block the pencils are $\Ein$-regular because of Lemmas \ref{Lem:RegularPencilOddCase} and \ref{Lem:RegularPencilEvenCase} respectively for $q$ odd and $q$ even, and on the remaining even blocks we apply Lemma \ref{lem:Regularity_r+r}.
\end{proof}

Finally, we deduce the topology of the fiber of $M_{\rho}$ in both the even and odd cases. 
\begin{theorem}[Topology of Fibers]
Let $p \geq 3$ be an integer, $\rho:\pi_1S\rightarrow \SO_0(p,p+1)$ a deformation of an $\iota$-Fuchsian in the space of $p$-Anosov representations and $M_{\rho} = \rho(\pi_1S)\backslash \Omega_{\rho}$. Then $M_{\rho}$ is diffeomorphic to a smooth $\mathfrak{F}_{\rho}$-fiber bundle over $S$, where:  
\begin{enumerate}[noitemsep, label=(\roman*)]
    \item $\mathfrak{F}_{\rho} = \Ein^{p-1,p}$ when $p$ is odd. 
    \item $\mathfrak{F}_{\rho} = \T^1\RP^{p-1}$ when $p$ is even. 
\end{enumerate}
\end{theorem}

\begin{proof}
Fix any point $P \in \Ha^2_{\iota}$. 
By Lemma \ref{Lem:PencilDeformationGeneralCase}, there is a family of pencils $\mathcal{P}_t \subset \T_P\X$, for $t \in [0,1]$, such that $\mathcal{P}= \mathcal{P}_1$ and $\mathcal{P}_0$ is a reducible pencil in the sense of Definition \ref{Defn:ReduciblePencilOdd} or \ref{Defn:ReduciblePencilEven}, depending on whether $p$ is odd or even, respectively. 

Now, Lemma \ref{lem:FiberOdd} says that $\mathcal{B}(\mathcal{P}_0)$ is diffeomorphic to $\Ein^{p-1,p}$ when $p$ is odd and Lemma \ref{Lem:FiberEvenQuasiNatural} says $\mathcal{B}(\mathcal{P}_0)$ is diffeomorphic to $\T^1\RP^{p-1}$ when $p$ is even. The result then follows from Lemma \ref{Lem:PencilDeformationGeneralCase} and Lemma \ref{Lem:RegularPencilDeformation}. 
\end{proof}

We now move on to the problem of the global topology of $M_{\rho}$ as a fiber bundle. To this end, we will use the Higgs bundle base of pencil $\mathcal{B}(\mathfrak{E}_0)$, in the sense of \eqref{HiggsBundleBasePencil}, to serve as a model for the quotients of domains of discontinuity. 

\begin{theorem}\label{thm:Deformation regular pencil blocks}
Let $\rho: \pi_1S\rightarrow \SO_0(p,p+1)$ be an $\iota_I$-Fuchsian representation. 
\begin{enumerate}[noitemsep, label=(\roman*)]
    \item The bundle of pencils $\mathfrak{E}_0\rightarrow S$ consists of reducible pencils with data $(\mathcal{E}_0^\R,\mathcal{J})$.
    \item The quotient of the domain of discontinuity $M_{\rho} = \rho(\pi_1S)\backslash \Omega_\rho$ is diffeomorphic to $\mathcal{B}(\mathfrak{E}_0)$ as a smooth fiber bundle over $S$. 
\end{enumerate}
\end{theorem}

\begin{proof}
(i) Fix $x \in S$. The claim is that $\mathfrak{E}_0|_x$ is a reducible pencil with data $(E_0,J)$, for $E_0 = \mathcal{E}_0^{\R}|_x$ and $J = \mathcal{J}|_x$, in the sense of Definitions \ref{Defn:ReduciblePencilEven}, \ref{Defn:ReduciblePencilOdd} in the respective cases of $p$ even and odd. 
The claim holds as a consequence of Lemma \ref{lem:ReduciblePencilOdd} and \ref{lem:ReduciblePencilEven} in the odd block, respectively for $q$ odd and even, and Lemma \ref{lem:Holomorphic r+r} for the other blocks. More precisely, for any $\theta \in \mathfrak{E}_t|_x$, we verify that for any $u \in \mathcal{U}^{\R}$ and $v \in \T_x\Sigma$, 
\[ \theta(\mathbf{i}v)(u)=(J\circ\theta)(v)(u). \]
If additionally $u \in \mathcal{U}^\R \cap \mathcal{E}_1$, then we also have $\theta(v)(J(u)) = (J\circ \theta)(v)(u)$. This proves (i). 
\medskip

(ii) Theorem \ref{Lem:NearestPointProjection} implies that the fiber bundle $\pi: M_\rho \rightarrow S$ is isomorphic to $\mathcal{B}(\mathfrak{E}_1)$ as smooth fiber bundles over $S$. Since the bundle of pencils $\mathfrak{E}_1$ can be smoothly deformed to $\mathfrak{E}_0$ by a path of $\Ein$-regular pencils, $M_\rho$ is diffeomorphic as a fiber bundle to $\mathcal{B}(\mathfrak{E}_0)$ by Corollary \ref{Cor:HiggsBundlePencil}.
\end{proof}

The remainder of the paper consists of applying Theorem \ref{thm:Deformation regular pencil blocks} to determine the global topology of $M_{\rho}$ in the odd and even cases.

\subsection{Global topology, \texorpdfstring{$p$}{p} odd}
\label{Sec:GlobalTopology_Odd}

We now give an initial description of the global topology in the odd case. The following definition will appear in the description.

\begin{definition}
For a $\Z_2$-graded vector bundles $E_+\oplus E_-$ over $S$, we shall write  $\Ein(E_+ \oplus E_-)$ for the associated fiber bundle $\sphere(E_+)\oplus \sphere(E_-)/\sim$, under the quotient $(p,u,v)\sim (p,-u,-v)$.
\end{definition}

Recall that $\mathcal{E}_1$ was defined in Section \ref{subsec:Sum of irreducibles}.

\begin{theorem}[Global Topology, Geometric]\label{Thm:GlobalTopology_Odd_First}
Let $p \geq 3$ be odd and $\rho: \pi_1S\rightarrow \SO_0(p,p+1)$ a deformation of an $\iota_I$-Fuchsian in the space of $p$-Anosov representations.
Then $M_{\rho} = \rho(\pi_1S)\backslash \Omega_{\rho}$ is diffeomorphic to the fiber bundle  $\Ein\left(\mathcal{U}^\R\oplus \left(\mathcal{V}^\R\cap\mathcal{E}_1\right)\right)$. 
\end{theorem}

\begin{proof}

By Corollary \ref{cor:FiberInvariance}, we can assume that $\rho$ is an $\iota_I$-Fuchsian as the diffeomorphism type of $M_\rho$ is invariant by continuous deformations in the space of $p$-Anosov representations.

    Using Theorem \ref{thm:Deformation regular pencil blocks}, there is a fiber bundle diffeomorphism between $M_{\rho}$ and $\mathcal{B}(\mathfrak{E}_0)$, where the Higgs bundle pencil $\mathfrak{E}_0$ is reducible with data $(\mathcal{E}_0^\R,\mathcal{J})$. Note that the bundle $\mathcal{E}_0^\R\cap \mathcal{U}=\mathcal{O}$ is trivial, so we can fix a global non-vanishing section, denoted $e_0 \in \Omega^0(S,\mathcal{O})$. Using this section $e_0$, we apply Lemma \ref{lem:FiberOdd} to construct a diffeomorphism for each $x\in S$ between $\mathcal{B}(\mathfrak{E}_0)|_x$ and $\Ein\left(\mathcal{U}^\R\oplus \left(\mathcal{V}^\R\cap\mathcal{E}_1\right)\right)|_x$, that varies smoothly in $x$. This defines therefore a diffeomorphism between $\mathcal{B}(\mathfrak{E}_0)$ and $\Ein\left(\mathcal{U}^\R\oplus \left(\mathcal{V}^\R\cap\mathcal{E}_1\right)\right)$. This concludes the proof.
\end{proof}

Using the previous result, we can now describe the topology of the quotients in a more direct way. We begin with the Hitchin case. 

\begin{theorem}[Global Topology in Hitchin case, $p$ odd]\label{thm:HitchinOddGlobal}
Let $\rho:\pi_1S \rightarrow \SO_0(p,p+1)$ be a Hitchin representation, with $p \geq 3$ odd. The manifold $M_{\rho} =\rho(\pi_1S)\backslash \Omega_{\rho}$ is diffeomorphic to: 
\begin{enumerate}[noitemsep, label=(\roman*)]
    \item $\Ein(\varepsilon^3_{\R} \oplus \K^3)$ when $p = 3$. 
    \item $S \times \Ein^{p-1,p}$ when $p \geq 5$. 
\end{enumerate}
\end{theorem}

\begin{proof}
We apply Theorem \ref{Thm:GlobalTopology_Odd_First}. The topology of $M_\rho$ depends on the real vector bundles $\mathcal{U}^\R$ and $\mathcal{V}^\R\cap \mathcal{E}_1$. More concretely on the Hitchin case we obtain:
\begin{align}\label{ExplicitDescriptionHitchinOdd}
   \sphere(\mathcal{U}^\R)\oplus \sphere(\mathcal{V}^{\R} \cap \V_1) \;\simeq\; \sphere(\K^{p-1} \oplus \K^{p-3} \oplus \cdots \oplus \K^2 \oplus \underline{\R})\oplus \sphere(\K^p \oplus \K^{p-2} \oplus \cdots \oplus \K^3).
\end{align}
The manifold $M_\rho$ is the quotient of this space by  the relation $(p, u,z)\sim (p,-u,-z)$.
However, there is a simpler topological description of this manifold. 

Recall that real vector bundles $E \rightarrow S$ of rank $\ge 3$ over a closed surface $S$ are classified by their Stiefel-Whitney classes $w_1(E), w_2(E)$. When $\rank(E) = 2$, we need the finer information afforded by the degree of $E$, rather than its mod 2 reduction $w_2(E)$. 
The vector bundles $\mathcal{U}^\R\rightarrow S$ and $\mathcal{V}^\R \cap \V_1 \rightarrow S$ are direct sums of complex line bundles of even degree, with a trivial real line bundle added for $\mathcal{U}^\R$. Hence, $w_1,w_2$ are both trivial for both bundles. It follows that  $\mathcal{U}^\R$ and $\mathcal{V}^\R \cap \V_1$ are smoothly trivializable when their ranks are at least three. 

The results of interest now follow. 

(i) The bundle $\mathcal{U}^\R$ has real rank three and $w_1=w_2=0$ and is thus smoothly trivial. When $p=3$, we have $\mathcal{V}^\R \cap \V_1 \cong \K^3$, which is a nontrivial complex line bundle over $S$ of degree $6g-6$. 
The description in (i) follows by \eqref{ExplicitDescriptionHitchinOdd}. 

(ii) Both $\mathcal{U}^\R$ and $\mathcal{V}^\R$ are smoothly trivial real vector bundles over $S$ when $p \geq 5$, as their real rank is at least $3$. Hence, $M_{\rho} \cong \mathbb{S}(\underline{\R^{p}})\oplus \mathbb{S}(\underline{\R^{p-1}})/\sim $ by \eqref{ExplicitDescriptionHitchinOdd}. However, the latter is just $(S \times \mathbb{S}^{p-1}\times \mathbb{S}^{p-2})/\sim$ under the quotient $(p,u,v)\sim(p,-u,-v)$, which is diffeomorphic to $S \times \Ein^{p-1,p-2}$.
\end{proof}

We now provide an analogous of $M_{\rho}$ for $p$-Anosov deformations of $\iota$-Fuchsians. 

To state the theorem, we must recall some topological invariants associated to representations $\rho:\pi_1S\rightarrow \SO_0(p,p+1)$. 
Given such a representation, we obtain an associated flat $\SO_0(p,p+1)$-principal bundle $\mathcal{P} \rightarrow S$ by $\mathcal{P} = \tilde{S}\times_{\rho} \SO_0(p,p+1)$. This principal bundle reduces structure group (uniquely up to isomorphism) to its maximal compact $K = \SO(p)\times \SO(p+1)$. We then obtain
two associated oriented vector bundles $E^p=E^p({\rho}),E^{p+1}= E^{p+1}({\rho})$, whose isomorphism classes are determined by their second Stiefel-Whitney classes
$w_2(E^p), w_2(E^{p+1}) \in H^2(S, \Z_2)$. 

For the representations of interest, $w_2(E^p)=w_2(E^{p+1})$. Thus, we shall write $w_2(\rho)$ instead. The global topology of $M_{\rho}$ is completely determined by this invariant.  

\begin{theorem}[Global Topology for $\iota$-Fuchsians, $p$ odd]\label{Thm:GlobalTopologyGeneral_Odd}
Let $p \geq 5$ be an odd integer, $S=S_g$ a closed surface of genus $g$, and $\rho: \pi_1S\rightarrow \SO_0(p,p+1)$ a deformation of an $\iota_I$-Fuchsian in the space of $p$-Anosov representations, and $2q+1$ the the unique odd part of $I$. 

Then $M_{\rho} = \rho(\pi_1S)\backslash \Omega_{\rho}$ is diffeomorphic to the fiber bundle $\Ein(E_+ \oplus E_-)$, where $E_+, E_-\rightarrow S$ are orientable real vector bundles of ranks $p$ and $p-1$, respectively, with $w_2=w_2(\rho)=(g-1)\frac{p-q}{2}\bmod 2$.
\end{theorem}

\begin{proof}
Observe that $\mathcal{U}^\R \cong E^p(\rho)$ and $\mathcal{V}^\R \cong E^{p+1}(\rho)$. A direct calculation with degrees verifies $w_2(\rho)= (g-1)\frac{p-q}{2} \bmod 2$.
Now, by Theorem \ref{Thm:GlobalTopology_Odd_First}, all that remains is to show $w_2(\mathcal{V}^\R \cap \V_1) = w_2(\mathcal{V}^\R)$. However, this follows by noting that $\deg(\mathcal{V}^\R \cap \V_0) \bmod 2 = \deg(\mathcal{K}) \bmod 2=0$. Indeed, one can easily check that on each block $\mathcal{E}'$, we have $\deg(\mathcal{E}' \cap \mathcal{E}_0\cap \mathcal{V}^{\R}) =0 \bmod 2$.
\end{proof}

\begin{corollary}

\label{cor:ProofTopology1}
When $p\geq 5$ is odd, 
the quotient $M_{\rho}$ satisfies the following: 
\begin{enumerate}[label=(\roman*)]
    \item The global topology of $M_{\rho}$ depends only on $w_2(\rho)$. 
    \item When $\rho$ is an $\iota_I$-Fuchsian, the fiber bundle $\pi: M_{\rho} \rightarrow S$ defined by Lemma \ref{Lem:NearestPointProjection} is non-trivial if and only if $w_2(\rho)\neq 0$.
    \item When $w_2(\rho) = 0$, then $M_{\rho}$ is diffeomorphic to $S \times \Ein^{p-1,p-2}$.
\end{enumerate}
\end{corollary}

\begin{proof}
(i) This is evident by Theorem \ref{Thm:GlobalTopologyGeneral_Odd}. Also, (iii) is immediate from (ii). 

Now, we prove (ii). Let $\rho$ be an $\iota_I$-Fuchsian. By Theorems \ref{thm:Deformation regular pencil blocks}, \ref{Thm:GlobalTopology_Odd_First}, and \ref{Thm:GlobalTopologyGeneral_Odd}, we have 
$M_{\rho} \cong \mathcal{B}(\mathfrak{E}_0) \cong \Ein(E_+ \oplus E_-)$ as smooth fiber bundles over $S$.
When $w_2(\rho)=0$, it is clear that $\Ein(E_+\oplus E_-)$ is a trivial fiber bundle. We now show the converse. 

\textbf{Step 1}: Suppose for contradiction that $\Ein(E_+ \oplus E_-)$ were a trivial fiber bundle but $w_2(E_-)\neq 0$.

(\textbf{1a}) Then its double cover $\sphere(E_+) \oplus \sphere(E_-)$, obtained by taking the universal cover of the fiber is also trivial. 

(\textbf{1b}) Crucially, if $\sphere(E_+) \oplus \sphere(E_-)$ is trivial, then $\sphere(E_+)$ must be homotopy equivalent to $S \times \sphere^{p-1}$. 

To prove this step, we introduce a \emph{trivial} $\sphere^{p-1}$-fiber-subbundle $\mathcal{F}$ of $\sphere(E_+)\oplus \sphere(E_-)$. To start, let $\tau: \underline{\sphere^{p-1} \times \sphere^{p-2}}\rightarrow \sphere(E_+)\oplus \sphere(E_-)$ be a fiber bundle diffeomorphism covering a diffeomorphism $\overline{\tau}:S \rightarrow S$. Fix a vector $v_-\in\mathbb{S}^{p-2}$. Consider $\mathcal{F}\coloneq \tau(\underline{\mathbb{S}^{p-1}\times\lbrace v_-\rbrace})\subset \sphere(E_+)\oplus \sphere(E_-)$, which is indeed trivial. For $v\in \underline{\mathbb{S}^{p-1}}$, we write: 
\[ \tau(v,v_-)=\big(\tau_+(v),\tau_-(v)\big)\in \sphere(E_+)\oplus \sphere(E_-). \]
Since $p\geq 5$, the homotopy group $\pi_{p-1}(\mathbb{S}^{p-1}\times \mathbb{S}^{p-2})$ is isomorphic to $\Z\oplus \Z_2$, and the projection $\mathbb{S}^{p-1}\times \mathbb{S}^{p-2}\to \mathbb{S}^{p-1}$ induces a surjective map $\pi_{p-1}(\mathbb{S}^{p-1}\times \mathbb{S}^{p-2})\to \pi_{p-1}(\mathbb{S}^{p-1})\cong \Z$.
Let $x\in S$ be arbitrary. In the fiber $\underline{\sphere^{p-1}\times \sphere^{p-2}}_{\mid x}$, the sphere $Y=\mathbb{S}^{p-1}\times\lbrace v_-\rbrace$ is a generator of the quotient of $\pi_{p-1}(\mathbb{S}^{p-1}\times \mathbb{S}^{p-2})\cong \Z\oplus \Z_2$ by its torsion subgroup, and hence the bundle diffeomorphism $\tau$ maps $[Y]$ to an element of $\pi_{p-1}(\sphere(E_+) \oplus \sphere(E_-)_{\mid \overline{\tau}(x)})$ that is a generator of the quotient by its torsion subgroup. Now, set $G=\pi_{p-1}(\sphere^{p-1}\times \sphere^{p-2})$, $\Z_2 \cong T<G$ its torsion subgroup, and $q:G\rightarrow \pi_{p-1}(\sphere^{p-1})$ 
the projection. The map $q$ kills the torsion and induces an isomorphism $\overline{q}:G/T \rightarrow \pi_{p-1}(\sphere^{p-1})$. 
Consequently, the map ${\tau_+}_{|x}:\underline{\mathbb{S}^{p-1}}_{\mid x}\to \mathbb{S}(E_+)_{\mid \overline{\tau}(x)}$ must be a generator of $\pi_{p-1}(\mathbb{S}(E_+)_{\mid \overline{\tau}(x)})$. Hence, ${\tau_+}_{\mid x}$ is homotopic to a homeomorphism, and in particular it is a homotopy equivalence.

We now use \cite[Theorem 6.3]{Dol63}, which implies that a continuous bundle map $f: E \rightarrow E'$ between continuous fiber bundles $E,E'$ over the same CW-complex $B$ such that $f$ is a fiberwise homotopy equivalence must actually be a homotopy equivalence of total spaces. We apply this to the bundle map $\tau_+:\underline{\mathbb{S}^{p-1}}\to \mathbb{S}(E_+)$ and conclude that $S \times \sphere^{p-1}$ and $\sphere(E_+)$ homotopy equivalent.

\medskip

\textbf{Step 2}: Now, suppose that $E \rightarrow S$ is a real rank $k \geq 4$ vector bundle with $w_2(E)\neq 0$. We will show that $\sphere(E)$ is not homotopy equivalent to $S \times \sphere^{k-1}$, which provides the desired contradiction.

The key is to prove the following sub-claim: $w_2( \T \sphere(E)) \neq 0$. Now, let $\mathscr{L} \rightarrow \sphere(E)$ be the tautological line bundle, $\mathscr{L}|_{(x,v)}= \R\{v\}$. The bundle $\mathscr{L}$ admits a non-vanishing section $s$ given by $\big[(x,v) \stackrel{s}{\longmapsto} (x,v,v) \big]$, so $\mathscr{L} \cong \varepsilon^1_{\R}$ is trivial. Next, observe that $ \T \sphere (E) \cong \pi^*\T S\oplus (\pi^*E)/\mathscr{L}$, where $\pi:\sphere(E)\rightarrow S$ is the bundle projection.
Hence, 
\[ \T \sphere(E) \oplus \varepsilon^1_{\R} \cong \pi^*\T S \oplus (\pi^*E)/\mathscr{L} \oplus \varepsilon^1_{\R} \cong \pi^*\T S \oplus \pi^* E.\]
Thus, 
\[ w_2(\T \sphere(E)) = w_2(\T \sphere(E)\oplus \varepsilon^1_{\R})=w_2(\pi^*(\T S \oplus E)) =_{(\star)} \pi^* (w_2(\T S \oplus E))= \pi^*(w_2(E)) .\] 
Here, in $(\star)$ we use \emph{naturality} of Stiefel-Whitney classes. 
Now, using the Gysin sequence (cf. \cite[page 438]{Hat01}), we find $\rank(E) \geq 4$ implies the fiber bundle projection $\pi: \sphere(E)\rightarrow S$ induces an isomorphism on second cohomology: $\pi^*: H^2(S) \stackrel{\cong}{\longrightarrow} H^2(\sphere(E))$. Thus, $w_2(\T\sphere(E)) \neq 0$ finally. On the other hand, $w_2( \T (S \times \sphere^{k-1}) )=0$ by direct calculation. We now conclude with the remarkable result of Wu \cite[Lemma 11.13]{MS74}, which implies the Stiefel-Whitney classes of the tangent bundle of a smooth manifold are homotopy invariants. In particular, by the calculations of $w_2$ above, $\sphere(E)$ and $S \times \sphere^{k-1}$ are not homotopy equivalent.

Hence, $\Ein(E_+\oplus E_-)$ is non-trivial when $w_2(E_+)\neq0$. 
\end{proof}

We now address one leftover case from Theorem \ref{Thm:GlobalTopologyGeneral_Odd}.
\begin{remark}
In Theorem \ref{Thm:GlobalTopologyGeneral_Odd}, we exclude the exceptional case $p=3$, as in this case the orientable vector bundle $E_-$ is of real rank two and is not classified by $w_2$. 

There are only two partitions corresponding to $\Ein$-regular representations in this case, namely $7=7$ and $7=3+2+2$. The former is the Hitchin case, addressed in Theorem \ref{thm:HitchinOddGlobal}. 
For $\iota$ associated to the latter partition, the global topology of the quotient $M_{\rho}$ is as follows. $E_+$ is of real rank 3, with  $w_2(E_+)=(g-1)\bmod 2$ and 
$E_- $ is of real rank 2 and degree $1-g$. Indeed, $E_- \cong \mathcal{V}^\R\cap \mathcal{E}_1 \cong \K^{-1/2}$ as a real vector bundle. Hence, $M_{\rho} \cong \Ein(E_+ \oplus \K^{-1/2}).$
\end{remark}

Theorem \ref{Thm:GlobalTopologyGeneral_Odd} says there are exactly two possibilities 
for the fiber bundle structure of $M_{\rho}$, with $p$ and $g$ fixed. If the genus $g$ of the surface is even, both of these possibilities are realized.

\begin{example}
We revisit the cases for $p=5$ as in Example \ref{Ex:p=5}. In this case, if the genus $g$ of $S$ is even, and $\iota$ is associated to the partition $11=7+2+2$, then $E_+$ and $E_-$ are non-trivial vector bundles over $S$. 
On the other hand, for $\iota$ associated to the partition $3+4+4$ or $3+2+2+2+2$, the manifold $M_{\rho}$ is diffeomorphic to $S \times \Ein^{4,3}$, regardless of genus. 

More generally, for the partition $(2p+1)=(2p-3)+2+2$, we always obtain $M_{\rho}$ of the form $\Ein(E_+\oplus E_-)$, where $E_+,E_-$ are nontrivial vector bundles with $w_2 \neq 0$. 
\end{example}

\subsection{Global topology, \texorpdfstring{$p$}{p} even}
\label{Sec:GlobalTopology_Even}

Let $\mathcal{T}$ be the unit tangent bundle to our closed surface $S$, naturally seen as a $\U(1)$-bundle. Let $\eta_0$ be the action of $\U(1)$ from Definition \ref{Defn:eta0Action}. We give a description of the global topology of $M_\rho$ that depends on this action.

\begin{theorem}\label{Thm:GlobalTopologyGeneral_Even}
Let $p \geq 4$ be an even integer and $\rho: \pi_1S\rightarrow \SO_0(p,p+1)$ a deformation of an $\iota_I$-Fuchsian in the space of $p$-Anosov representations. 

Then the manifold $M_{\rho} = \rho(\pi_1S)\backslash \Omega_{\rho}$ is diffeomorphic to the fiber bundle $\mathcal{T}\oplus\T^1\mathbb{P}(\mathcal{U}^\R)/\sim$ where $(x\cdot \lambda,v)\sim (x, \eta_0(\lambda)\cdot v)$ for all $\lambda \in \mathrm{U}(1)$ and $(x,v)\in \mathcal{T}\oplus\T^1\mathbb{P}(\mathcal{U}^\R)$. 
\end{theorem}

Recall that for the fiber bundle direct sum $\mathcal{T}\oplus\T^1\mathbb{P}(\mathcal{U}^\R)$, for all $x\in S$, the fiber is given by $\big(\mathcal{T}\oplus\T^1\mathbb{P}(\mathcal{U}^\R) \big)|_x=\mathcal{T}|_x\times\T^1\mathbb{P}(\mathcal{U}^\R)|_x$. 
\begin{proof}

First note that we can assume that $\rho$ is $\iota_I$-Fuchsian as the topology of $M_\rho$ is invariant by continuous deformations inside the space of $p$-Anosov representations, see Corollary \ref{cor:FiberInvariance}.

 Using Theorem \ref{thm:Deformation regular pencil blocks} once again, we see that $M_{\rho}$ is diffeomorphic to the fiber bundle $\mathcal{B}(\mathscr{P}_0)$ where the pencils $\mathscr{P}_0|_x$ are all reducible pencils with data $(\mathcal{E}_0^\R,\mathcal{J})$. Note that the bundle $\mathcal{E}_0=\mathcal{O}$ is trivial so we can fix a global section of norm one, denoted $e\in \Omega^0(S,\mathcal{O}^\R)$. 

We now apply the isomorphism from 
Theorem \ref{thm:FiberEvenNatural}(i) fiberwise. In doing so, we obtain a surjective submersion $\Theta: \mathcal{T} \oplus \T^1\mathbb{P}(P) \rightarrow \mathcal{B}(\mathscr{P}_0)$ via $\Theta(\psi,[u,v]\,) = \Theta_{\psi}([u,v])$ realizing the domain as a circle bundle over $\mathcal{B}(\mathscr{P}_0)$ by the Ehresmann lemma. 
By Theorem \ref{thm:FiberEvenNatural}(ii), $\Theta$ descends to a fiber bundle diffeomorphism between $\big(\mathcal{T}\oplus \T^1\mathbb{P}(\mathcal{U}^\R)\big)/\sim $ and $\mathcal{B}(\mathscr{P}_0)$.
\end{proof}

\begin{remark}
The description in Theorem \ref{Thm:GlobalTopologyGeneral_Even} is similar to that of an associated fiber bundle $\mathscr{F} = \mathscr{P}\times_{\phi} F$ to a principal bundle $\mathscr{P}$, but slightly different. Rather, in this case, we combine a principal $G$-bundle $\mathcal{P}$ with an $F$-fiber bundle $\mathcal{F}$ to obtain a $(G\times F)$-fiber bundle $\mathcal{P} \oplus \mathcal{F}$, which we quotient down to an $F$-fiber bundle via a diagonal $G$-action. In fact, one can view the construction of $\mathscr{F} = \mathscr{P}\times_{\phi} F$ in the same way, where $\mathcal{F} = \underline{F}$ is trivial.

Of note is the fact that $\mathcal{P}$ and $\eta$ are fixed in the theorem, while $\mathcal{F}$ varies. In the associated bundle picture, usually the action $\phi$ varies instead (e.g. in \cite{AMTW25}).
\end{remark}

Recall that $p =2k$ is even presently. If $p$ is also divisible by four, we can simplify our description of $M_\rho$ substantially. For this purpose, the following technical lemma is the key.\footnote{This argument shows that the complex vector bundles $(\mathscr{L}_{J_0}^\C)^\bot$ and $(\mathscr{L}^{\C}_{-J_0})^\bot$ are isomorphic when $p=4k$. Such a phenomenon is prohibited when $p=4k+2$ by the fact that $c_1((\mathscr{L}_{J_0}^\C)^\bot)\neq 0$. We do not currently know whether $\eta_0^2$ is homotopically trivial in the latter case.}

\begin{lemma}[$\eta_0^2$ Triviality]\label{Lem:Eta0^2NullHomotopy}
Let $p$ be divisible by four. Then the loop $\eta_0^2:\U(1) \rightarrow \Diff(\T^1\RP^{p-1})$ is homotopically trivial. 
\end{lemma}

\begin{proof}
We first remark in general. 
The essence of the proof is to consider the dependence of the action of $\U(1)$ on $\T^1\RP^{2k-1}$ on the background identification $(\R^{2k},J)\cong \C^k$. Recall that the symmetric space $\mathscr{J}_k\coloneqq\mathrm{O}(2k)/\U(k)$ is the space of orthogonal complex structures on $\R^{2k}$. 
We shall study the map $A:\U(1)\times \mathscr{J}_k\rightarrow \Diff(\T^1\RP^{p-1})$
given as follows. For $J \in \mathscr{J}_k$ fixed, we obtain a decomposition 
\[ \T\RP^{p-1} = J\mathscr{L} \oplus (\mathscr{L}^\C_{J})^\bot,\]
where $J\mathscr{L}|_{\ell} = \Hom_{\R}(\ell, J\ell)$ and $(\mathscr{L}^\C_J)^\bot = \Hom_{\R}(\ell, (\ell \oplus J\ell)^\bot).$
We then obtain an action $A(\cdot, J)$ of $\U(1)$ by bundle automorphisms on $\T\RP^{2k-1}$ by letting $\lambda \in \U(1)$, thought of as unit complex numbers, act by scalar multiplication on each fiber. Since $J$ is orthogonal, this action preserves $\T^1\RP^{2k-1}$, and can be written $A(\cdot, J):\U(1)\rightarrow \Diff(\T^1\RP^{p-1})$. 

We first make two observations:
\begin{enumerate}
    \item For a fixed background complex structure $J_0 \in \mathscr{J}_k$ for which we built the map $\eta_0$, we can write $\eta_0= A(\cdot, J_0)$. 
    \item The map $\eta_0^{-1}:\U(1)\rightarrow \Diff(\T^1\RP^{p-1})$, rotating the fibers of $(\mathscr{L}^\C_{J_0})^\bot$ in the opposite direction, is determined by $\eta_0^{-1}(\cdot)=A(\cdot, -J_0)$. 
\end{enumerate}
Now, $\eta_0^2$ is nullhomotopic in $\Diff(\T^1\RP^{p-1})$ if $\eta_0 \sim \eta_0^{-1}$ in $\pi_1(\Diff(\T^1\RP^{p-1}))$. We prove this now. 

Finally, we suppose $k$ is even. To produce such a homotopy, we use that $J$ and $-J$ are in the same connected component of $\mathscr{J}_k$ when $k$ is even. Hence, we can take a smooth path $\gamma_{t}:[0,1]\rightarrow \mathscr{J}_k$ such that $\gamma_t(0) = J_0$ and $\gamma_t(1)=-J_0$. 

Hence, we obtain a homotopy 
$H:\U(1)\times[0,1]\rightarrow \Diff(\T^1(\RP^{2k-1}))$ given by $H(\lambda, s) = A(\lambda, \gamma_s)$.
By the observations above, $H(\cdot, 0) = \eta_0$ and $H(\cdot, 1) = \eta_0^{-1}$. Thus, $\eta_0 \sim \eta_0^{-1}$. 
\end{proof}

We now present the simpler description of the global topology when $p$ is divisible by four.

\begin{theorem}[Global Topology for $\iota$-Fuchsians, $p$ even]\label{Thm:GlobalTopologyEven-Geometric}
Let $p \geq 4$ be a multiple of four and $\rho: \pi_1S\rightarrow \SO_0(p,p+1)$ a deformation of an $\iota_I$-Fuchsian in the space of $p$-Anosov representations. Then $M_{\rho} = \rho(\pi_1S)\backslash \Omega_{\rho}$ is diffeomorphic to the fiber bundle $\T^1\mathbb{P}(\mathcal{U}^\R)$.
\end{theorem}

\begin{proof}
First, we make a small alteration to Theorem \ref{Thm:GlobalTopologyGeneral_Even} by using a square root of $\mathcal{T}$. Let us denote $\mathcal{T}'$ to be the principal circle bundle over $\Sigma$ with Euler number $g-1$. Observe that the principle bundle $\mathcal{T}$ is isomorphic to an associated bundle of $\mathcal{T}'$ as follows: 
\[ \mathcal{T}= \mathcal{T}'\times_{m} \U(1),\]
where $m:\U(1)\rightarrow \U(1)$ is $m(\lambda) = \lambda^2$. There is a natural 2-1 covering map $\pi:\mathcal{T}'\rightarrow \mathcal{T}$ given by $\pi(\psi) = [(\psi,\id)]$. Note the property $\pi(\psi \cdot \lambda) =\pi(\psi)\cdot \lambda^2$. 

Now, we define the bundle map 
\[ \pi\oplus \id\colon\mathcal{T}' \oplus \T^1\mathbb{P}(\mathcal{U}^\R){\longrightarrow} \mathcal{T}\oplus \T^1\mathbb{P}(\mathcal{U}^\R).\]
The given map will descend to a diffeomorphism. Define a quotient on the domain by declaring 
$(\psi\cdot \lambda, x)\sim(\psi, \eta_0^2(\lambda)\cdot x)$ for $(\psi,x) \in \mathcal{T}' \oplus \T^1\mathbb{P}(\mathcal{U}^\R)$. 
The quotient improves the situation: $\pi \oplus \id$ descends to a fiber bundle diffeomorphism 
\[ \big(\mathcal{T}'\oplus \T^1\mathbb{P}(\mathcal{U}^\R)\big)/\sim \;\cong \big(\mathcal{T}\oplus \T^1\mathbb{P}(\mathcal{U}^\R)\big)/\sim \; \cong \; M_{\rho},\]
where the latter isomorphism is from Theorem \ref{Thm:GlobalTopologyGeneral_Even}. 
However, since $\eta_0^2$ is homotopically trivial by Lemma \ref{Lem:Eta0^2NullHomotopy},
we conclude by noting $\big(\mathcal{T}'\oplus \T^1\mathbb{P}(\mathcal{U}^\R)\big)/\sim \;\cong \T^1\mathbb{P}(\mathcal{U}^\R)$.
\end{proof}

\begin{corollary}
When $p$ is a multiple of four, the manifold $M_{\rho} =\rho(\pi_1S)\backslash \Omega_{\rho}$ satisfies the following:
\begin{enumerate}[noitemsep,label=(\roman*)]
    \item The topology depends only on $w_2(\rho) = w_2(\mathcal{U}^\R)$.
    \item If $\rho$ is an $\iota_I$-Fuchsian, then the fibration $\pi:M_{\rho} \rightarrow S$ from Lemma \ref{Lem:NearestPointProjection} is trivial if and only if $w_2(\rho) \neq 0$. 
    \item If $w_2(\rho)=\frac{p-q}{2}(g-1) \bmod 2 = 0$, then $M_{\rho}$ is diffeomorphic to $S \times \T^1\RP^{p-1}$. 
\end{enumerate}
\end{corollary}

\begin{proof}
    Again, we need only prove the claim on non-triviality in (ii). Thus, let us suppose $w_2(\rho)\neq 0$. We shall adapt the argument in Corollary \ref{cor:ProofTopology1} to the present setting. 
    
    \textbf{Step 1}: Suppose, for contradiction, that $\T^1\mathbb{P}(E)$ were a trivial fiber bundle. 
    
    \textbf{(1a)}. If $\T^1\mathbb{P}(E)$ is trivial as a fiber bundle, then $\T^1 \sphere(E)$ must be trivial as a fiber bundle. Here, use that $\T^1\mathbb{S}^{p-1}$ is the universal cover of $\T^1\RP^{p-1}$ for $p \geq 4$. 

    \textbf{(1b)} We show that when $\T^1\sphere(E)$ is a trivial fiber bundle, $\sphere(E)$ must be homotopy equivalent to $S \times \sphere^{p-1}$. 
    This sub-step is the heart of the proof. 
    We prove this now for $p\geq 8$ and deal with the remaining case $p=4$ at the end of the proof. 
    
    Take $\tau:\underline{\T^1\sphere^{p-1}}\rightarrow \T^1\sphere(E)$ a global trivialization. Using a section $s$ of the projection $\T^1\sphere^{p-1}\rightarrow \sphere^{p-1}$, we obtain a fiber subbundle 
    $\tau_+: \underline{\sphere^{p-1}}\hookrightarrow \underline{\T^1\sphere^{p-1}}$ via $(x,v) \mapsto (x, v,s(v))$. We then have a map $\beta: \underline{\sphere^{p-1}} \rightarrow \sphere(E)$ by $\beta = \pi \circ \tau_-$, where $\pi:\T^1\sphere(E) \rightarrow \sphere(E)$. Now, $\beta$ is a fiber bundle map. As in Corollary \ref{cor:ProofTopology1}, it suffices to show $\beta$ is a homotopy equivalence fiberwise. We do this now. 
    
    Observe that the map $\pi_{p-1}(\T^1\mathbb{S}^{p-1})\to \pi_{p-1}(\mathbb{S}^p)$  induced by the projection $\T^1\mathbb{S}^{p-1}\to \mathbb{S}^{p-1}$ is surjective. Indeed, this is true because $\T^1\sphere^{p-1}$ has a section when $p$ is even. By the long exact sequence of homotopy groups, for $p \geq 8$ even, one finds $\pi_{p-1}(\T^1\mathbb{S}^{p-1})\leq \Z\oplus \Z_2$. Arguing as in Corollary \ref{cor:ProofTopology1}, one shows $\beta_{\mid x}: \sphere^{p-1}\rightarrow \sphere(E)_{\mid x}$ is homotopic to a homeomorphism and thus is a homotopy equivalence. In other words, we have found $\beta$ is fiberwise a homotopy equivalence. 

    \textbf{Step 2}. We obtain a contradiction again since $\sphere(E)$ cannot be homotopy equivalent to $S \times \sphere^{p-1}$ when $w_2(E)\neq0$, by the proof of Corollary \ref{cor:ProofTopology1}. This concludes the proof for $p\geq 8$ divisible by four.

\medskip

We now address the leftover case $p=4$. Here, $\T^1\mathbb{S}^{3}\cong \mathbb{S}^3\times \mathbb{S}^2$ and $\pi_3(\mathbb{S}^3\times \mathbb{S}^2)\cong \Z\oplus \Z$ has a normal subgroup associated to $\pi_3(\lbrace *\rbrace\times\mathbb{S}^2)\cong \Z$. In order to apply the same technique as in Corollary \ref{cor:ProofTopology1}, we claim that any homeomorphism $\phi:\mathbb{S}^3\times \mathbb{S}^2\to \mathbb{S}^3\times \mathbb{S}^2$ induces a map $\phi_*$ on the third homotopy groups that descends to a group isomorphism of the quotient $\pi_3(\mathbb{S}^3\times \mathbb{S}^2)/\pi_3(\lbrace *\rbrace\times\mathbb{S}^2)$. Then the same arguments from Corollary \ref{cor:ProofTopology1} apply when replacing the torsion subgroup by $\pi_3(\lbrace *\rbrace\times\mathbb{S}^2)\cong \Z$.

We now check this last claim. All homology groups below have integer coefficients. The Hurewicz map $h:\pi_3(\mathbb{S}^3\times \mathbb{S}^2)\to H_3(\mathbb{S}^3\times \mathbb{S}^2)$ satisfies $\phi_*\circ h=h\circ \phi_*$ where now $\phi_*$ denotes the induced action on $H_3$ or $\pi_3$ as appropriate. Using the Kunneth formula, we see that $h$ is surjective and its kernel is exactly $\pi_3(\lbrace *\rbrace\times\mathbb{S}^2)$. In particular, the quotient $\pi_3(\mathbb{S}^3\times \mathbb{S}^2)/\pi_3(\lbrace *\rbrace\times\mathbb{S}^2)$ is naturally identified with $H_3(\mathbb{S}^3\times \mathbb{S}^2)$ and the homomorphism induced by $\phi_*$ on this quotient is an isomorphism. 

Finally, we conclude as before. If $E \rightarrow S$ is a real rank four vector bundle and $w_2(E) \neq 0$, then $\T^1\sphere(E)$ cannot be trivial as a fiber bundle or else $\sphere(E)$ would be homotopy equivalent to $S \times \sphere^3$. 
\end{proof}

\bibliographystyle{alpha}
\bibliography{bib}

\end{document}